\documentstyle[12pt,ctagsplt]{amsart}

\input cyracc.def
\font\tencyr=wncyr10
\def\cyr{\tencyr\cyracc}
\font\tencyrit=wncyi10
\def\cyrit{\tencyrit\cyracc}

\newtheorem{reidslemma}{Reid's Lemma} 
\newtheorem{definition}{Definition} 
\newtheorem{maintheorem}{Main Theorem} 
\newtheorem{corollary}{Corollary} 
\newtheorem{theorem}{Theorem}
\newtheorem{lemma}{Lemma}
\newtheorem{proposition}{Proposition}
\newtheorem{remark}{Remark}

\newcommand{\ip}[2]{( #1 \, | \, #2 )}
\newcommand{\coloneq}{\mathrel{:=}}

\newcommand{\W}{\frak W}
\newcommand{\w}{\frak W}
\newcommand{\Sym}{\frak S}
\newcommand{\sym}{\frak S}
\newcommand{\isom}{\cong}
\newcommand{\union}{\cup}
\renewcommand{\a}{\alpha}
\renewcommand{\c}{\gamma}
\renewcommand{\d}{\delta}
\newcommand{\e}{\varepsilon}
\newcommand{\f}{\varphi}
\newcommand{\s}{\sigma}
\renewcommand{\l}{\lambda}
\newcommand{\m}{\mu}
\renewcommand{\O}{{\cal O}}
\newcommand{\mm}{\frak m}
\newcommand{\A}{{\Bbb A}}
\newcommand{\C}{{\Bbb C}}
\renewcommand{\H}{\bold{H}}
\renewcommand{\P}{{\Bbb P}}
\newcommand{\R}{{\Bbb R}}
\newcommand{\Z}{{\Bbb Z}}
\newcommand{\Proj}{\text{\bf Proj}}
\newcommand{\Aut}{\operatorname{Aut}}
\newcommand{\Coker}{\operatorname{Coker}}
\newcommand{\Def}{\operatorname{Def}}
\newcommand{\Ker}{\operatorname{Ker}}
\newcommand{\pr}{\operatorname{pr}}
\newcommand{\PRes}{\operatorname{PRes}}
\newcommand{\Symm}{\operatorname{Sym}}

\begin{document}

\title[Singularities with Small Resolutions]{Gorenstein
Threefold Singularities with Small Resolutions
via Invariant Theory for Weyl Groups}
\author[S. Katz]{Sheldon Katz}
\address{Department of Mathematics\\
Oklahoma State University\\
Stillwater, OK  74078
}
\email{katz@@math.okstate.edu}
\author[D. R. Morrison]{David R. Morrison}
\address{Department of Mathematics\\
Duke University\\
Durham,
NC 27706}
\email{drm@@math.duke.edu}
\date{}

\maketitle

\section*{}

A fundamental new type of birational modification which first occurs in
dimension three is the {\em simple flip}.  This is a
birational map $Y \dasharrow Y^+$
 which induces an isomorphism
$(Y - C) \cong (Y^+ - C^+),$
 where $C$ and $C^+$ are smooth
rational curves such that $K_Y \cdot C < 0$ and $K_{Y^+} \cdot C^+ > 0$.
($Y$ and $Y^+$ should be allowed to have ``terminal" singularities.)
Mori's celebrated theorem \cite{[Mr]}
shows that these flips exist when numerically
expected.

A closely related type of modification is the {\em simple flop}.
This has a similar definition, except that  $Y$ and $Y^+$
should be Gorenstein, with $K_Y \cdot C  = K_{Y^+} \cdot C^+ = 0$.
(This is more than an analogy:
every flip has a branched double cover which is a flop,
and this construction was used in Mori's proof.)
For both flips and flops, the curves $C$ and $C^+$ can be contracted
to points (in $Y$ and $Y^+$, respectively), yielding the same normal
variety $X$.  The birational map $Y \dasharrow Y^+$ can thus be
described in terms of the two contraction morphisms
$\pi\colon Y \to X$ and $\pi^+ \colon Y^+ \to X$.

In this paper, we study the case of simple flops
with $Y$  smooth, so that $\pi\colon Y \to X$ is an irreducible small
resolution of a Gorenstein threefold singularity
$P \in X$.  (It is called ``small" because the exceptional set is
a curve rather than a divisor, and ``irreducible" because the curve has
only one component.)
In fact, the study of flops with $Y$ smooth
can be reduced to a study of Gorenstein
threefold singularities with small resolutions,
thanks to a theorem
of Reid \cite{[R]} (cf.\ also Koll\'ar \cite{[Kol]}) which produces a
second small resolution $\pi^+\colon Y^+ \to X$ out of the original
$\pi\colon Y \to X$.
Pinkham \cite{[P]} showed that in this situation, $X$ is Gorenstein if it is
merely assumed to be Cohen-Macaulay.

Early examples of small resolutions were constructed in an {\em ad hoc\/}
manner.
The first class of examples
is given by the singularities $x^2+y^2 + z^2 + t^{2k}=0$:  a particularly
nice description of the associated simple flops
 (with pictures!)\ can be found in a paper of
Reid \cite[\S 5]{[R]}.
These flops are exactly those for which the normal bundle of $C$ in $Y$
 is $\O(-1) \oplus \O(-1)$ or $\O \oplus \O(-2)$.
A second  example of simple flops, in which the normal bundle is
$\O(1) \oplus \O(-3)$, was found by
Laufer \cite{[L]};  variants of this example were also investigated by
the second author, Pinkham \cite{[P]}, and Reid \cite{[R]}.
Some other flops were studied in previous work of the authors
(\cite{[M]} and \cite{[K]}).

In general,
if $Y\to X$ is a small resolution of the isolated Gorenstein threefold
singularity $P\in X$, then by a lemma of Reid \cite{[R]}, the general
hyperplane section of $X$ through $P$ has a rational double point
at $P$, and the proper transform of that surface on $Y$ gives a
``partial resolution" of the rational double point.
Pinkham \cite{[P]} used this observation to give a construction which
includes all possible Gorenstein threefold
singularities with small resolutions.  In the irreducible case, each such
singularity can be described by a map from the disk to a space
$\PRes (S,v)$ which parametrizes deformations which partially resolve.
The natural map $\PRes (S,v) \to \Def (S)$ to the deformation
space then gives a map from the disk to $\Def (S)$ which describes
the space $X$ as the pullback of a (semi-)universal family.

At first glance, Pinkham's construction appears to give a countable number
of families of Gorenstein threefold singularities with irreducible
small resolutions.
The discrete data which appear in the construction are the  type
of the rational double point, together with a choice of component in the
exceptional divisor of the minimal resolution of that point.
However, the construction uses a {\em particular\/} hyperplane section through
$P$, and there is no guarantee that this hyperplane section is ``general".
(Examples for which it is {\em not\/} general were known to Pinkham at the time
he gave the construction.)

We have discovered that there are in fact only six families of
Gorenstein threefold singularities with irreducible
small resolutions.  They can be
distinguished by a very simple invariant (the ``length") which was
introduced by Koll\'ar \cite[pp.~95, 96]{[CKM]} a few years ago.  The precise
statement of our main theorem can be found in section 1.

Our methods do much more than simply characterize the six families:
our techniques can be used to calculate the map
$\PRes (S,v) \to \Def (S)$
 quite explicitly in each of the six cases.
Composing this map with a general map from the disk to
$\PRes (S,v)$ describes the most general Gorenstein threefold with
irreducible small resolution of each type.  For example, when the length
is 1, the
result is precisely the class of examples $x^2+y^2+z^2+t^{2k}=0$ mentioned
above.

\bigskip

In order to prove the main theorem, we need to solve a
 fundamental problem:  compute the singularity type of the
generic hyperplane section for irreducible small resolutions which
have been produced by Pinkham's construction.
The main tools we use to solve this problem are derived from the theory of
simultaneous resolution of rational double points, as developed by
Brieskorn \cite{[Bri0]},
\cite{[Bri]},
\cite{[Bri-Nice]} and Tyurina \cite{[T]}.  In this theory, the
deformation space $\Def (S)$ is identified with a quotient
space $V/\W$, where $V$ is the complex root space and $\W$ is the
Weyl group of a certain root system $R$.
(The Dynkin diagram of this root system coincides with the dual
graph of the minimal resolution of the rational double point.)
 A simultaneous resolution
of the semi-universal deformation over $V/\W$ is possible after making the
base change by $V \to V/\W$.  Pinkham's construction of the space
$\PRes (S,v)$ included an identification of it with the quotient
space $V/\W_0$, where $\W_0$ is the Weyl group of a certain subsystem
of $R$.  Thus, an understanding of the map $\PRes (S,v) \to \Def (S)$
can be obtained if one has an adequate understanding of the invariant
theory of the Weyl groups $\W$ and $\W_0$.

Our first theorem is a slight extension of the fundamental theorem on
simultaneous resolution.  We need some additional information about
the loci in which certain curves in the exceptional set deform,
but most importantly, we need a version of the proof which provides
the simultaneous resolution in an explicitly computable form.  With
a small amount of modification, the approach of Tyurina \cite{[T]} as amplified
by Pinkham \cite{[P-RDP]} provides us with what we need.

The theory of simultaneous resolution implies that the coefficients
of the defining polynomial of a semi-universal deformation serve as generators
of the algebra of $\W$-invariant polynomials on $V$.  (By a theorem
of Coxeter \cite{[Cx]} and Chevalley \cite{[Chv]},
this is a free polynomial algebra.)
Our next theorems show that these generators are in an appropriate sense
unique up to $\C^*$-action,
 and can be computed from the invariant theory alone.  This is
an important step, since it allows us to recover the defining polynomials
of various semi-universal deformations directly from the invariant theory.
It has the effect of reducing our fundamental problem to a problem in
the invariant theory of Weyl groups.

The problem in invariant theory can be solved by hand when the Weyl
groups in question are those associated to the root systems $A_{n-1}$
or $D_n$, but in the cases of $E_6$, $E_7$, and $E_8$,
the invariant theory is not so well
understood.  We computed the invariants in these cases with the
aid of the symbolic computing languages {\sc maple} and {\sc reduce},
running on a Macintosh II and on a Sun 3/60
workstation.  We have also developed a number of tools for manipulating
these invariants, which enabled us to extract enough information
relating $\W$-invariants to $\W_0$-invariants to solve our fundamental
problem.

We describe our algorithm in sufficient detail that it could be
implemented in any symbolic manipulation language.  We have, however,
been somewhat selective in displaying the calculated results in
the paper.
The implementations in {\sc maple} and {\sc reduce} were carried out
independently by the second and first authors, respectively,
and run on different machines.
We are happy to report that the two sets of
calculated results agree in every particular.

One of the by-products of the work described here is a
new explicit set of generators for the invariants of the Weyl groups
of $E_6$, $E_7$, and $E_8$.  Following ideas of Tyurina,
these generators are obtained by computing
the anti-pluricanonical mappings for an appropriate family of del~Pezzo
surfaces, and putting the defining polynomials of the images into a
semi-universal form.  This calculation was attempted in 1918 by
C. C. Bramble \cite{[Bra]} in the case of $E_7$, but our machine-aided
calculation
reveals that Bramble made a few errors.
A corrected version of Bramble's calculation appears in Appendix 2;
the corresponding calculation for $E_6$ is in Appendix 1.

As the final draft of this paper was being prepared, we learned of some
recent work of Shioda \cite{[Shioda]} who has also done an explicit
calculation of simultaneous resolutions using a quite
different approach.

J\'anos Koll\'ar has communicated to us another possible
approach to proving our
main theorem \cite{[Kpvt]}, using techniques of Clemens and Jim\'{e}nez
\cite{[J]} and an
analysis of higher order neighborhoods similar to that used by Mori for
extremal rays in \cite{[Mr]}.

The paper is organized as follows.  Section 1 contains the statement
of the main theorem, and a discussion of all of the ingredients needed
to state it, including  Pinkham's
construction and Koll\'ar's ``length" invariant.
In section 2, we set up the notation for root systems, and
we introduce the notion of {\em distinguished polynomials}, which provide an
efficient means for comparing invariants.  In section 3, we establish
notation for the rational double points, and state the
simultaneous resolution theorem in the form we will need it.  We prove
that theorem in sections 4 and 5, giving explicit constructions of
simultaneous resolutions.
In section 6, we show that the simultaneous resolutions we have
constructed are essentially unique, and we use this uniqueness to
analyze the defining polynomials of simultaneous partial resolutions.
In section 7, we use the distinguished polynomials to compare
invariants in different coordinate systems, and section 8 contains the proof of
the
main theorem, modulo some computer calculations.  The computer-dependent
portions of the paper have been isolated in sections 9 and 10,
which describe the algorithm for explicitly putting the defining polynomials
of the simultaneous
resolutions for $E_6$, $E_7$, and $E_8$ into semi-universal form,
and for manipulating the
invariant polynomials to finish the proof of the main theorem.

We would like to thank Robert Bryant, Jonathan Wahl, and
Klaus Wirthm\"uller for helpful conversations, and Yunliang Yu for
sharing his {\sc maple}-to-\TeX\ package with us.
The first author would also like to thank the University of Bayreuth for a
stimulating working environment while part of this work was done.
We both gratefully acknowledge the facilities of Internet and {\sc bitnet},
which made our collaboration possible.

\section{Reid's lemma, Pinkham's construction, Koll\'ar's invariant, and
the statement of the main theorem.}

The analysis of Gorenstein threefold singularities with small resolutions
begins with the following lemma of Reid \cite[(1.1),(1.14)]{[R]}.

\begin{reidslemma}
Let $\pi\colon Y \to X$ be a resolution of an isolated
 Gorenstein threefold singularity
$P \in X$.  Suppose that the exceptional set of $\pi$ has pure
dimension 1.  Let $X_0$ be a generic hyperplane section of $X$
which passes through $P$.
Then $X_0$ has a rational double point at $P$.

Moreover, if $X_0$ is any hyperplane section through $P$
with a rational double point,
and  $Y_0$ is its proper transform, then $Y_0$ is normal,
and the minimal resolution
$Z_0 \to X_0$ factors through the induced map
$\pi |_{Y_0}\colon Y_0 \to X_0$.
\end{reidslemma}

Following Wahl \cite{[W]}, a map $Y_0 \to X_0$
through which the minimal resolution $Z_0 \to X_0$ factors
is called
a {\em partial resolution\/} of $X_0$ (provided that $Y_0$ is normal).
There is a natural graph
associated to such a map.  Start with the dual graph $\Gamma$ of the
components of the exceptional divisor of the minimal resolution
$Z_0 \to X_0$.
The curves contracted by $Y_0 \to X_0$ correspond to vertices in the
graph which span a subgraph $\Gamma_0$;
we call $\Gamma_0 \subset \Gamma$ the {\em partial resolution graph\/}
of $\pi$.
Figure~\ref{figure1} below and figure~\ref{figure2}
in section 7 show some partial resolution graphs.
In both figures, the vertices
corresponding to $\Gamma_0$ are shown with an open circle $(\circ)$,
while those corresponding to $\Gamma - \Gamma_0$ are shown with a
closed circle $(\bullet)$.

Based on Reid's lemma, Pinkham \cite{[P]} gave a construction which
in principle describes
all Gorenstein threefold singularities with small resolutions.
Let ${\cal Y} \to \Def (Y_0)$ and
${\cal X} \to \Def (X_0)$ be semi-universal deformations of
$Y_0$ and $X_0$, respectively.  By a theorem of Wahl
\cite[Theorem 1.4]{[Wahl]}, all deformations
of $Y_0$ blow down to give deformations of $X_0$, and  there is a
natural map $\tau\colon \Def (Y_0) \to \Def (X_0)$.  If we choose a local
equation $\{f=0\}$ for $X_0$ in a neighborhood $U$ of $P$ in $X$, then
$f$ can be regarded as a map  $f\colon U \to \Delta$, where
$\Delta \subset \C$ is a  small disk, such that
$X_0 \cap U = f^{-1}(0)$.
The fibers of $f$ give a deformation of $X_0$, so there is
 a classifying map $\mu_f\colon \Delta \to \Def (X_0)$
which enables us to recover $U$.  Similarly, the fibers of
$f \circ \pi$ give a deformation of $Y_0$, and  the classifying
map $\mu_{f \circ \pi}\colon \Delta \to \Def (Y_0)$
satisfies  $\tau \circ \mu_{f \circ \pi} = \mu_f$.

The only data which are really necessary to describe the map
$\tau\colon \Def (Y_0) \to \Def (X_0)$ are the type $S$ of the rational
double point ($S$ is one of $A_{n-1}$, $D_n$, $E_6$, $E_7$, or $E_8$)
and the subgraph $\Gamma_0 \subset \Gamma$.  Having fixed such data, we denote
the map by $\tau\colon \PRes (S,\Gamma_0) \to \Def (S)$.
The space $\PRes (S,\Gamma_0)$ represents the functor of deformations of
a singularity of type $S$ which {\em p\/}artially {\em res\/}olve,
with partial resolution type specified by $\Gamma_0$.
When $\Gamma_0$ consists of a single vertex $v$, we abbreviate the
notation to $\PRes (S,v)$.

Pinkham's general construction then goes as follows. Fixing $S$ and
$\Gamma_0$ determines the map
 $\tau\colon \PRes (S,\Gamma_0) \to \Def (S)$.
For any map $\nu\colon \Delta \to \PRes (S,\Gamma_0)$ there is an
induced map $\mu \coloneq \tau \circ \nu$.
Pulling back the semi-universal families by $\nu$ and $\mu$ gives
  threefolds $Y \to X$.  Pinkham shows that if $\nu$ is sufficiently
general, then $Y$ is smooth, $X$ is Gorenstein with an isolated
singular point, and $Y \to X$ is a small resolution.  Reid's lemma
implies that all Gorenstein threefold
 singularities with small resolutions arise in this way.

However, as Pinkham observed \cite[p.~367]{[P]},
 ``Although we have a construction for the singularities that arise, we
do not really have a classification.  For example, the construction does not
explain how to get the generic hyperplane \dots".
In fact, the generic hyperplane section of $X$ can have a much simpler
singularity type than $S$ (the type of the hyperplane section which
was used to construct $Y \to X$).
The main theorem of this
paper will compute the singularity type of the
generic hyperplane section in the case that the
exceptional set of $\pi\colon Y \to X$ is irreducible.

\bigskip

There is a fundamental invariant of  singularities $P \in X$ with an
 irreducible small resolution
which was introduced by
Koll\'ar \cite[pp.~95, 96]{[CKM]}.
Let $\mm_{P,X}$ denote the maximal ideal sheaf of $P$ in
$X$.

\begin{definition}
Let $\pi\colon Y\to X$ be an irreducible small resolution of an isolated
threefold singularity $P$.  (That is, the exceptional set
$\pi^{-1}(P)$ is an irreducible curve $C$.)
  The\/ {\em length\/} of $P$ is the length at the generic point of the
scheme supported on $C$ with structure sheaf ${\cal O}_Y/\pi^{-1}(\mm_{P,X})$.
\end{definition}

\begin{remark}  A different invariant, the\/ {\em normal bundle sequence\/}
of the exceptional curve, was introduced by Pinkham in \cite{[P]}.  This is the
collection of normal bundles of the sequence of curves $C_1,\ldots $ beginning
with the
exceptional curve $C$, with $C_{i+1}$ being the negative section of the
exceptional divisor of the blowup along $C_i$.  The sequence terminates when
the normal bundle becomes $\O(-1)\oplus\O(-1)$.
There are
five possible types of normal bundle sequences,
listed in \cite[p.~367]{[P]}.
The length and the normal bundle sequence measure a kind of ``thickness'' of
the singularity and exceptional curve, respectively.  It can
be shown that these invariants are related as follows.
Length 1 corresponds to Pinkham's cases (1) or (2), length 2 corresponds to
case (3), lengths 3, 4, 6 each correspond to case (4),
and length 5 corresponds to
case (5).
\end{remark}

In the Gorenstein case,
if $X_0$ is any hyperplane section of $X$ with rational
double points, then the length of $X$ can also be computed from the
sheaf ${\cal O}_{Y_0}/\pi^{-1}(\mm_{P,X_0})$.
The length therefore coincides with the multiplicity of the curve $C$ in
the maximal ideal cycle of the rational double point.
(The maximal ideal cycle is the same as Artin's {\em fundamental cycle\/}
 \cite{[A]}.)

This property of the length puts a constraint on the possible deformations of
the hyperplane section.  If we have a family $(X_0)_t$
 of hyperplane sections
through $P$, then we can get different partial resolution graphs for different
values of $t$.  However, the multiplicity of the curve $C$ in the maximal ideal
cycle must be independent of $t$, since it always coincides with the
length.  This leads to the
following definition.

\begin{definition}
A partial resolution $Y_0 \to X_0$ with irreducible
exceptional set $C$ is\/ {\em primitive\/} if for
every nontrivial 1-parameter deformation of $Y_0 \to X_0$ for which
$C$ deforms,  the
multiplicity of $C$ in the maximal ideal cycle at the generic point of the
family is strictly less than the multiplicity at the special point.
\end{definition}

The point of the definition is that if $Y\to X$ is an irreducible
 small resolution with a
hyperplane section $X_0\subset X$ yielding a primitive $Y_0\to X_0$,
 then
a 1-parameter deformation from $X_0$ to a general hyperplane section
through $P$ must be trivial.  Thus,
$X_0$ has
the same singularity type as a general hyperplane section of $X$ through $P$.

The primitive partial resolution
graphs can be computed from the deformation theory of rational
double points.  There are exactly six of them:  the six graphs
shown in figure~\ref{figure1}.
The numbers labeling the vertices in the figure
 are the multiplicities of the corresponding curves in
the maximal ideal cycle.
We omit the proof that these six graphs are the only primitive ones;
the proof is a fairly straightforward computation, but the result
is not logically necessary
for the purposes of this paper.
In fact, it can be obtained as a corollary of our main theorem.

\begin{figure}[t]
\begin{picture}(2,1)(1.9,.5)
\thicklines
\put(2.9,1){\circle{.075}}
\put(2.775,.7){\makebox(.25,.25){\footnotesize 1}}
\end{picture}
\hspace*{\fill}
\begin{picture}(3,1)(1.65,.5)
\thicklines
\put(1.9,1){\circle*{.075}}
\put(1.9,1){\line(1,0){.5}}
\put(2.4,1){\circle*{.075}}
\put(2.4,1){\line(1,0){.4625}}
\put(2.9,1){\circle{.075}}
\put(2.9,.9625){\line(0,-1){.4625}}
\put(2.9,.5){\circle*{.075}}
\put(2.9375,1){\line(1,0){.4625}}
\put(3.4,1){\circle*{.075}}
\put(3.4,1){\line(1,0){.5}}
\put(3.9,1){\circle*{.075}}
\put(3.9,1){\line(1,0){.5}}
\put(4.4,1){\circle*{.075}}
\put(1.775,1.05){\makebox(.25,.25){\footnotesize 2}}
\put(2.275,1.05){\makebox(.25,.25){\footnotesize 3}}
\put(2.775,1.05){\makebox(.25,.25){\footnotesize 4}}
\put(2.925,.375){\makebox(.25,.25){\footnotesize 2}}
\put(3.275,1.05){\makebox(.25,.25){\footnotesize 3}}
\put(3.775,1.05){\makebox(.25,.25){\footnotesize 2}}
\put(4.275,1.05){\makebox(.25,.25){\footnotesize 1}}
\end{picture}

\bigskip

\bigskip

\begin{picture}(2,1)(1.9,.5)
\thicklines
\put(2.4,1){\circle*{.075}}
\put(2.4,1){\line(1,0){.4625}}
\put(2.9,1){\circle{.075}}
\put(2.9,1.0375){\line(0,1){.4625}}
\put(2.9,1.5){\circle*{.075}}
\put(2.9375,1){\line(1,0){.4625}}
\put(3.4,1){\circle*{.075}}
\put(2.275,.7){\makebox(.25,.25){\footnotesize 1}}
\put(2.775,.7){\makebox(.25,.25){\footnotesize 2}}
\put(2.925,1.375){\makebox(.25,.25){\footnotesize 1}}
\put(3.275,.7){\makebox(.25,.25){\footnotesize 1}}
\end{picture}
\hspace*{\fill}
\begin{picture}(3,1)(1.9,.5)
\thicklines
\put(1.9,1){\circle*{.075}}
\put(1.9,1){\line(1,0){.5}}
\put(2.4,1){\circle*{.075}}
\put(2.4,1){\line(1,0){.5}}
\put(2.9,1){\circle*{.075}}
\put(2.9,1){\line(0,-1){.5}}
\put(2.9,.5){\circle*{.075}}
\put(2.9,1){\line(1,0){.4625}}
\put(3.4,1){\circle{.075}}
\put(3.4375,1){\line(1,0){.4625}}
\put(3.9,1){\circle*{.075}}
\put(3.9,1){\line(1,0){.5}}
\put(4.4,1){\circle*{.075}}
\put(4.4,1){\line(1,0){.5}}
\put(4.9,1){\circle*{.075}}
\put(1.775,1.05){\makebox(.25,.25){\footnotesize 2}}
\put(2.275,1.05){\makebox(.25,.25){\footnotesize 4}}
\put(2.775,1.05){\makebox(.25,.25){\footnotesize 6}}
\put(2.925,.375){\makebox(.25,.25){\footnotesize 3}}
\put(3.275,1.05){\makebox(.25,.25){\footnotesize 5}}
\put(3.775,1.05){\makebox(.25,.25){\footnotesize 4}}
\put(4.275,1.05){\makebox(.25,.25){\footnotesize 3}}
\put(4.775,1.05){\makebox(.25,.25){\footnotesize 2}}
\end{picture}

\bigskip

\bigskip

\begin{picture}(2,1)(1.9,.5)
\thicklines
\put(1.9,1){\circle*{.075}}
\put(1.9,1){\line(1,0){.5}}
\put(2.4,1){\circle*{.075}}
\put(2.4,1){\line(1,0){.4625}}
\put(2.9,1){\circle{.075}}
\put(2.9,.9625){\line(0,-1){.4625}}
\put(2.9,.5){\circle*{.075}}
\put(2.9375,1){\line(1,0){.4625}}
\put(3.4,1){\circle*{.075}}
\put(3.4,1){\line(1,0){.5}}
\put(3.9,1){\circle*{.075}}
\put(1.775,1.05){\makebox(.25,.25){\footnotesize 1}}
\put(2.275,1.05){\makebox(.25,.25){\footnotesize 2}}
\put(2.775,1.05){\makebox(.25,.25){\footnotesize 3}}
\put(2.925,.375){\makebox(.25,.25){\footnotesize 2}}
\put(3.275,1.05){\makebox(.25,.25){\footnotesize 2}}
\put(3.775,1.05){\makebox(.25,.25){\footnotesize 1}}
\end{picture}
\hspace*{\fill}
\begin{picture}(3,1)(1.9,.5)
\thicklines
\put(1.9,1){\circle*{.075}}
\put(1.9,1){\line(1,0){.5}}
\put(2.4,1){\circle*{.075}}
\put(2.4,1){\line(1,0){.4625}}
\put(2.9,1){\circle{.075}}
\put(2.9,.9625){\line(0,-1){.4625}}
\put(2.9,.5){\circle*{.075}}
\put(2.9375,1){\line(1,0){.4625}}
\put(3.4,1){\circle*{.075}}
\put(3.4,1){\line(1,0){.5}}
\put(3.9,1){\circle*{.075}}
\put(3.9,1){\line(1,0){.5}}
\put(4.4,1){\circle*{.075}}
\put(4.4,1){\line(1,0){.5}}
\put(4.9,1){\circle*{.075}}
\put(1.775,1.05){\makebox(.25,.25){\footnotesize 2}}
\put(2.275,1.05){\makebox(.25,.25){\footnotesize 4}}
\put(2.775,1.05){\makebox(.25,.25){\footnotesize 6}}
\put(2.925,.375){\makebox(.25,.25){\footnotesize 3}}
\put(3.275,1.05){\makebox(.25,.25){\footnotesize 5}}
\put(3.775,1.05){\makebox(.25,.25){\footnotesize 4}}
\put(4.275,1.05){\makebox(.25,.25){\footnotesize 3}}
\put(4.775,1.05){\makebox(.25,.25){\footnotesize 2}}
\end{picture}

\caption{}
\label{figure1}
\end{figure}

Note that the primitive partial resolution graphs shown
are uniquely determined by the multiplicity of
$C$ in the maximal ideal cycle.  It
follows from our proof of the main theorem that every partial resolution graph
admits a deformation to a primitive graph through graphs for which $C$ has
constant multiplicity in the maximal ideal cycle.

\bigskip

We can now state our main theorem.

\begin{maintheorem}
The generic hyperplane section of
an isolated Gorenstein threefold singularity
which has an irreducible small resolution
defines one of the
 primitive partial resolution graphs in figure~\ref{figure1}.
Conversely, given any
such primitive
partial resolution graph, there exists an
irreducible small resolution $Y\to X$
whose general hyperplane
section is described by that partial resolution graph.
\end{maintheorem}

It follows that there are exactly six basic types of simple flops on
smooth threefolds.

\begin{corollary}
The general hyperplane section of $X$ is uniquely determined
by the length of the singular point $P$.
\end{corollary}

Two special cases of this theorem were known previously.  If $X$ has
some hyperplane section of type $A_n$, then the theorem asserts that
the general hyperplane section has type $A_1$.  This was known to
several people (including Mori, Pinkham,
Shepherd-Barron and the second author) about
10 years ago, but has apparently never been published.  And in \cite{[K]},
the first author
did the case in which $X$ has some hyperplane section of
type $D_4$.

\bigskip

The converse statement in the theorem follows from the discussion above:
Pinkham showed that examples exist for any partial resolution graph,
and as we have
observed, a hyperplane section leading to a primitive partial resolution
graph must be the
general hyperplane section.  The first statement in the main theorem
will be proved in section 8.

\section{Root systems, Weyl groups,
and distinguished polynomials.}

In this section, we establish some notation for certain
root systems and their Weyl groups.
As is customary in algebraic
geometry, we take root systems to have
negative definite inner product;
aside from this, we follow the notation
of Bourbaki \cite{[Bo]} fairly closely for $A_{n-1}$
and $D_ n$, making some departures in
the case $E_ n$.

We begin with an inner product space
$\H ^ {n+1}$ ($n \ge 1$)
over a field $k$ of characteristic 0,
with orthogonal basis
$e_ 0 , e _ 1, \dots ,
e _ {n}$ such that
$\ip{e_ 0}{e_ 0} = 1$ and
$\ip{e _ i}{
e _ i} = -1$ for $i \ge 1$.
We let $e_0^*,e_1^*,\ldots,e_n^*$ be the dual basis of the dual space
$(\H ^ {n+1})^*$.

$\H ^ {n+1}$ contains the
lattice
\[\H ^ {n+1}_ \Z  \coloneq
\{ x = \sum \xi_i e_i \in \H ^ {n+1} \ | \
\xi_i \in \Z  \text{ and } \sum \xi_i \in 2 \Z  \}.\]

We define three subspaces of
$\H ^ {n+1}$.  $V_ {E_ n}$ is
the orthogonal complement of
the special vector $k \coloneq -3
e_ 0  + \sum _ {i=1}^ n e _ i$,
$V_ {D_ n}$ is the orthogonal
complement of the first basis vector
$e_ 0 $, and $V_ {A_ {n-1}}
 = V_ {E_ n} \cap V_ {D_ n}$
is the orthogonal complement of the
subspace spanned by $k$
 and $e_ 0 $.
The lattice  $\H ^ {n+1}_ \Z $
induces (by intersection with
these subspaces) lattices $L_ {E_ n}$,
$L_ {D_ n}$, and $L_ {A_ {n-1}}$.
A {\em root\/} in one of these lattices
is an element of norm $-2$.

The set of roots in the lattice $L_ {E_ n}$ ($n \ge 3$),
$L_ {D_ n}$ ($n \ge 2$), or $L_ {A_ {n-1}}$ ($n \ge 1$)
is called a {\em root system\/}
of type ${E_ n}$, ${D_ n}$, or ${A_ {n-1}}$, respectively.
We denote this set of roots by $R_ {E_ n}$ (resp. $R_ {D_ n}$ or
$R_ {A_ {n-1}}$).
The lattice itself is called the {\em root lattice},
and the vector space $V_ {E_ n}$ (resp. $V_ {D_ n}$ or
$V_ {A_ {n-1}}$) is called the {\em root space over $k$}.
It is customary in Lie theory to take $k = \R $; here, we are
primarily interested in the case $k = \C$, and in that case
we refer to $V$ as the {\em complex root space}.

We have included a degenerate case $A_0$, as well as
two cases in which the root system is reducible:
$R_ {D_ 2} \isom R_ {A_ 1} \union R_ {A_ 1}$, and
$R_ {E_ 3} \isom R_ {A_ 2} \union R_ {A_ 1}$.
In addition, there are some isomorphisms among the irreducible ones:
$R_ {D_ 3} \isom R_ {A_ 3}$,
$R_ {E_ 4} \isom R_ {A_ 4}$, and
$R_ {E_ 5} \isom R_ {D_ 5}$.

When $n \ge 1$, the lattice $L_ {A_ {n-1}}$ can be
generated by the
{\em root basis\/} $v_ 1 , \dots ,
v_ {n-1}$, where
$v _ i \coloneq e _ i
- e _ {i+1}$.  This root
basis has as its Dynkin diagram $\Gamma_{A_{n-1}}$:

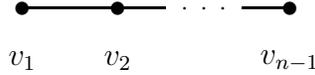
\begin{center}

\begin{picture}(2.4,.5)(.5,.6)
\thicklines
\put(1,1){\line(1,0){.5}}
\put(1,1){\circle*{.075}}
\put(1.5,1){\circle*{.075}}
\put(1.5,1){\line(1,0){.25}}
\put(1.85,1){\circle*{.02}}
\put(1.95,1){\circle*{.02}}
\put(2.05,1){\circle*{.02}}
\put(2.15,1){\line(1,0){.25}}
\put(2.4,1){\circle*{.075}}
\put(.875,.6){\makebox(.25,.25){$v _ 1$}}
\put(1.375,.6){\makebox(.25,.25){$v _ 2$}}
\put(2.275,.6){\makebox(.25,.25){$v _ {n-1}$}}
\end{picture}

\end{center}

\noindent
Notice that in the degenerate case $A_0$, we have the empty root basis,
which forms
a basis for the zero vector space $V_{A_0}$.

When $n \ge 2$, adding the root $v _ n\coloneq
e _ {n-1} + e _ n$ to the set $v_1,\ldots,v_{n-1}$
produces a root basis
$v_ 1 , \dots , v_ {n-1}, v _ n$
of $L_ {D_ n}$, which has Dynkin diagram $\Gamma_{D_n}$:

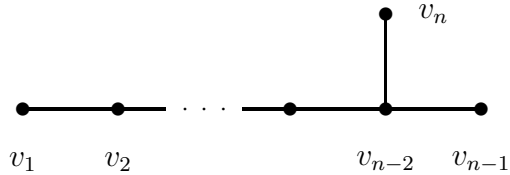
\begin{center}

\begin{picture}(3.4,1)(.5,.6)
\thicklines
\put(1,1){\line(1,0){.5}}
\put(1,1){\circle*{.075}}
\put(1.5,1){\circle*{.075}}
\put(1.5,1){\line(1,0){.25}}
\put(1.85,1){\circle*{.02}}
\put(1.95,1){\circle*{.02}}
\put(2.05,1){\circle*{.02}}
\put(2.15,1){\line(1,0){.25}}
\put(2.4,1){\circle*{.075}}
\put(2.4,1){\line(1,0){.5}}
\put(2.9,1){\circle*{.075}}
\put(2.9,1){\line(0,1){.5}}
\put(2.9,1.5){\circle*{.075}}
\put(2.9,1){\line(1,0){.5}}
\put(3.4,1){\circle*{.075}}
\put(.875,.6){\makebox(.25,.25){$v _ 1$}}
\put(1.375,.6){\makebox(.25,.25){$v _ 2$}}
\put(2.775,.6){\makebox(.25,.25){$v _ {n-2}$}}
\put(3.275,.6){\makebox(.25,.25){$v _ {n-1}$}}
\put(3.025,1.375){\makebox(.25,.25){$v _ n$}}
\end{picture}

\end{center}

\noindent
Only the two end vertices $v_{n-1}$ and $v_n$ appear in the reducible
case $D_2$.

Finally, if $3 \le n \le 8$,
adding the root
$v _ 0 \coloneq e_ 0  -
e _ 1 - e _ 2 -
e _ 3$ to the set $v_1,\ldots,v_{n-1}$
produces a root basis $v _ 0,
v_ 1 , \dots , v_ {n-1}$
of $L_ {E_ n}$.  This has Dynkin
diagram $\Gamma_{E_n}$:

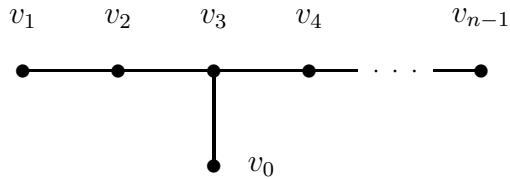
\begin{center}

\begin{picture}(3.4,1)(1.4,.4)
\thicklines
\put(1.9,1){\circle*{.075}}
\put(1.9,1){\line(1,0){.5}}
\put(2.4,1){\circle*{.075}}
\put(2.4,1){\line(1,0){.5}}
\put(2.9,1){\circle*{.075}}
\put(2.9,1){\line(0,-1){.5}}
\put(2.9,.5){\circle*{.075}}
\put(2.9,1){\line(1,0){.5}}
\put(3.4,1){\circle*{.075}}
\put(3.4,1){\line(1,0){.25}}
\put(3.75,1){\circle*{.02}}
\put(3.85,1){\circle*{.02}}
\put(3.95,1){\circle*{.02}}
\put(4.05,1){\line(1,0){.25}}
\put(4.3,1){\circle*{.075}}
\put(1.775,1.15){\makebox(.25,.25){$v _ 1$}}
\put(2.275,1.15){\makebox(.25,.25){$v _ 2$}}
\put(2.775,1.15){\makebox(.25,.25){$v _ 3$}}
\put(3.025,.375){\makebox(.25,.25){$v _ 0$}}
\put(3.275,1.15){\makebox(.25,.25){$v _ 4$}}
\put(4.175,1.15){\makebox(.25,.25){$v _ {n-1}$}}
\end{picture}

\end{center}

\noindent
In the reducible case $E_3$, this diagram consists of $v_1$ joined
to $v_2$ on the left, and $v_0$ below.

Many of our constructions for a root system $R$ will depend on its
{\em type\/} $S$, which is one of $A_{n-1}$, $D_n$, or $E_n$ in all
cases we consider.  Constructions which depend on the type may fail to
be invariant under the isomorphisms of root systems
$R_ {D_ 3} \isom R_ {A_ 3}$,
$R_ {E_ 4} \isom R_ {A_ 4}$, and
$R_ {E_ 5} \isom R_ {D_ 5}$.

We  single out certain linear functions on our root spaces.
Let $V$ be one of our root spaces $V_{A_{n-1}}$, $V_{D_n}$, or
$V_{E_n}$.  We define the {\em distinguished functionals\/} on the
root space $V$ to be the $n$ functions $t_1,\ldots,t_n$ given by
\[t_i \coloneq (\frac13 e_0^* + e_i^*)|_V.\]

It is not difficult to express
the distinguished functionals
$t_1,\ldots,t_n$  in terms of the dual basis $v_{\alpha}^*, v_{\alpha + 1}^*,
 \ldots, v_{\beta}^*$ of the root basis
$v_{\alpha}, v_{\alpha + 1}, \ldots, v_{\beta}$; the result is shown
in the third column of table~\ref{table12}.  (In order to simplify
the notation, we have introduced a few extra $v_i^*$'s into the formulas,
which should be set equal to zero, as indicated in the second column of
the table.)  In the fourth column of  table~\ref{table12}, we have
solved for  the dual
basis in terms of the distinguished functionals.  A key step in doing
this is to let $s_1\coloneq t_1+\cdots+t_n$  be the sum of the distinguished
functionals, and to note that $s_1=0$ in the case of $A_{n-1}$.

\begin{table}[t]
\begin{center}
\begin{tabular}{|c|c|l|l|}  \hline
 & & &  \\
$A_{n-1}$ & $v_0^*=0$ &
   $t_i = - v_{i-1}^* + v_i^*$  &
   $v_i^*  = t_1 + \cdots + t_i$ \\
 & $v_n^*=0$ & \quad ($1 \le i \le n$)  &  \quad ($1 \le i \le n-1$) \\
 & &  &  \\ \hline
 & &  &  \\
 & & $t_i      =  - v_{i-1}^* + v_i^*$
 &  $v_i^*      =  s_1 - t_{i+1} - \cdots - t_n $ \\
 & & \quad ($1 \le i \le n-2$)  &  \quad ($1 \le i \le n-2$) \\
 & &  &  \\
$D_n$ & $v_0^* = 0$ & $t_{n-1}  =  - v_{n-2}^* + v_{n-1}^* + v_n^*$
 &  $v_{n-1}^*  =  \frac12 s_1 - t_n$ \\
 & &  &  \\
 & & $t_{n}    =  - v_{n-1}^* + v_n^*$
 &  $v_n^*      =  \frac12 s_1$ \\
 & & &  \\ \hline
 & &  &  \\
 & & $t_1  =  - \frac23 v_0^* + v_1^* $
 &  $v_0^* = \frac{3}{n-9} s_1 $ \\
 & &  &  \\
 & & $t_2  =  - \frac23 v_0^* - v_1^* + v_2^* $
 &  $v_1^*  =  \frac{2}{n-9} s_1 + t_1 $ \\
$E_n$ & $v_n^*=0$ &  &  \\
 & & $t_3  =  - \frac23 v_0^* - v_2^* + v_3^* $
 &  $v_2^*  =  \frac{4}{n-9} s_1 + t_1 + t_2 $ \\
 & &  &  \\
 & & $t_i  =  \frac13 v_0^* - v_{i-1}^* + v_i^* $
 &  $v_i^*  =  \frac{9-i}{n-9} s_1 + t_1 + \cdots + t_i $ \\
 & & \quad ($4 \le i \le n$)  &  \quad ($3 \le i \le n-1$) \\
 & &  &  \\  \hline
\end{tabular}
\end{center}

\medskip

\caption{}
\label{table12}
\end{table}

The formulas in table~\ref{table12} imply that
the ring $\C [V]$ of polynomial functions on $V$ is generated
by $t_1,\ldots,t_n$, and is in fact isomorphic to
$\C [t_1,\ldots,t_n]/I$, where $I = (t_1+\cdots+t_n)$ in case
$A_{n-1}$, and $I=(0)$ otherwise.

We collect the distinguished functionals into the {\em distinguished
polynomial\/}
\begin{equation}
\label{eqf1}
f_S(U;t)\coloneq\prod_{i=1}^n (U + t_i),
\end{equation}
where $S$ denotes the type of the root system
 (one of $A_{n-1}$, $D_n$, or $E_n$)
and $t$ denotes $(t_1,\ldots,t_n)$.
If we let the symmetric group $\Sym_n$ act by permuting $\{t_1,\ldots,t_n\}$
in the usual way, then
the distinguished polynomial $f_S(U;t)$ lies in the subring
$(\C [V]^{\sym_n})[U]$ of polynomials
whose coefficients are invariant under $\Sym_n$.  It can thus
be written in the form
\[
f_S(U;t)=U^n + \sum_{i=1}^n s_i U^{n-i},
\]
where
the coefficients $s_1,\ldots,s_n$ are the elementary symmetric
functions of the distinguished functionals $t_1,\ldots,t_n$.
(Note that  $s_1$ is the sum of the distinguished functionals,
agreeing with the definition given above.)
However, the product expansion  for $f_S(U;t)$ given in
equation (\ref{eqf1}) is only valid
in the larger ring $(\C [V])[U]$.

It is important to observe that the construction of both
the distinguished functionals and the
distinguished polynomial depend on the
type of the root system having been identified as one of $A_{n-1}$,
$D_n$, or $E_n$.
For example, although $E_4$ and $A_4$
are isomorphic as root
systems, they have {\em different\/} distinguished functionals
and polynomials.
(Even the degrees of the distinguished polynomials are different.)

Notice also that when our definitions are applied to the degenerate case
$A_0$, the single distinguished functional $t_1$ is identically zero,
 and the distinguished
polynomial is simply $f_{A_0}(U;t) = U \in \C [U]$.

We next describe the action of the Weyl group $\W$ of our
root system.
$\W$ is the subgroup of $\Aut (V^*)$ generated by the reflections
$r_v$ in the roots $v$ which belong to the root basis.  We are particularly
interested in the invariant polynomials for the action of $\W$ on $V^*$.
A theorem of Coxeter \cite{[Cx]} and Chevalley \cite{[Chv]} guarantees
that the ring $\C [V]^{\w}$ of invariant polynomials is a free
polynomial algebra over $\C$.  This implies that the quotient
space $V/\W$ is smooth.

The reflections $r_v \in \Aut (V^*)$ which generate $\W$ are
restrictions of reflections in $\Aut ((\H ^{n+1})^*)$, which
we also denote by $r_v$.  The action of each $r_v$ on $(\H ^{n+1})^*$
is easily described.  For example,
the reflections $r_{v_i}$ for $1 \le i \le n-1$
act on $(\H ^ {n+1})^*$
by fixing $e_0^*$ and mapping $e_j^*$
to $e_{\sigma (j)}^*$, where
$\sigma$ is the simple transposition $(i,i+1)$ of the set $\{1,\ldots,n\}$.
It follows that the action of these $r_{v_i}$ on $V^*$ maps
$t_j^*$ to $t_{\sigma(j)}^*$ in each of the three cases
$V=V_{A_{n-1}}$, $V=V_{D_n}$, and $V=V_{E_n}$.

The Weyl group $\W_ {A_ {n-1}}$ is the
group generated by the reflections
 $r_ {v _ 1}, \dots , r_ {v _ {n-1}}$, and therefore coincides with
 the symmetric
group $\Sym_n$.
Since $t_1 + \cdots + t_n = 0$, the
invariant polynomials for the action of $\W_{A_{n-1}}$ are generated by
the nonzero coefficients
$s_2,\ldots,s_n$ of the distinguished polynomial $f_{A_{n-1}}(U;t)$.

In the case of the Weyl group $\W_ {D_ n}$, the reflections
$r_{v_i}$ for $1 \le i \le n-1$ again act as permutations of
$t_1,\ldots,t_n$, although this time the sum $s_1$ is not necessarily
zero.  $\W_ {D_ n}$ is generated by those reflections together
with $r_{v_n}$.  The action of $r_{v_n}$ on $(\H ^ {n+1})^*$
is not hard to compute:  we have
\[
r_ {v _ n}(e_i^*) =
\begin{cases}
\quad e_i^* & \text{if  $0 \le i \le n-2$}, \\
- e_ n^* & \text{if $i=n-1$}, \\
- e_ {n-1}^* & \text{if $i=n$}.
\end{cases}
\]
It follows that the action on $V_{D_n}^*$ is given by
\[
r_ {v _ n}(t_i) =
\begin{cases}
\quad t_i & \text{if $1 \le i \le n-2$}, \\
- t_ n & \text{if $i=n-1$}, \\
- t_ {n-1} & \text{if $i=n$}.
\end{cases}
\]
The invariant polynomials for this action are generated by the
product
$t_1 \cdots t_n$,
together with the elementary
symmetric functions of the squares $t_1^2,\ldots,t_n^2$.
These invariants are captured by the constant term $s_n = f_{D_n}(0;t)$
of the distinguished polynomial, together with the coefficients of
a polynomial $g_{D_n}(Z;t)$
whose defining property is
\begin{equation}
\label{eqg1}
g_{D_n}(Z;t) \coloneq \prod_{i=1}^n (Z + t_i^2).
\end{equation}
(Note that $g_{D_n}(0;t) = f_{D_n}(0;t)^2 = s_n^2$, so its constant term is not
needed to generate the $\W$-invariants.)
We call $g_{D_n}(Z;t)$ the {\em second distinguished polynomial\/}; it is
only defined for root systems of type $D_n$.

The second distinguished polynomial $g_{D_n}(Z;t)$ lies in the ring
$(\C [V]^{\w})[Z]$ of polynomials whose coefficients are $\W$-invariants,
although the defining property (\ref{eqg1}) only holds in the larger
ring $(\C [V])[Z]$.  The second distinguished polynomial could
also have been defined by the property
\begin{equation}
\label{eqg2}
g_{D_n}(- U^2;t) = f_{D_n}(U;t) \cdot f_{D_n}(-U;t)
\end{equation}
which holds in the auxiliary ring $(\C [V]^{\sym_n})[U]$.  Following
Tyurina~\cite{[T]}, we can express this property directly in the subring
$(\C [V]^{\sym_n})[Z]$ of $(\C [V])[Z]$ in the following way.
Collect terms of even and odd degree in the distinguished polynomial:
\begin{equation}
\label{eqg3a}
f_{D_n}(U;t) = U \cdot P_{D_n}(-U^2;t) + Q_{D_n}(-U^2;t).
\end{equation}
This defines two polynomials
\[P_{D_n}(Z;t), Q_{D_n}(Z;t) \in (\C [V]^{\sym_n})[Z].\]
Then equation (\ref{eqg2}) can be rewritten as
\begin{equation}
\label{eqg3b}
g_{D_n}(Z;t) = Z \cdot P_{D_n}(Z;t)^2 + Q_{D_n}(Z;t)^2.
\end{equation}
We also record for later use the existence of a polynomial
$G_{D_n}(Z,U;t) \in (\C [V]^{\sym_n})[Z,U]$ such that
\begin{equation}
\label{eqG1}
U \cdot P_{D_n}(Z;t) + Q_{D_n}(Z;t)
= (Z + U^2) \cdot G_{D_n}(Z,U;t) + f_{D_n}(U;t).
\end{equation}
Since $G_{D_n}(Z,U;t)$ has degree $n-2$ in $U$, if we define
\begin{equation}
\label{eqG2}
\widetilde{G}_{D_n}(Z,u,v;t) \coloneq v^{n-2} \cdot G_{D_n}(Z,u/v;t),
\end{equation}
then $\widetilde{G}_{D_n}(Z,u,v;t)$ is a polynomial which is
 homogeneous in $u, v$.

Finally, in the case of the Weyl group $\W_ {E_ n}$, the reflections
$r_{v_i}$ for $1 \le i \le n-1$ still act as permutations of
$t_1,\ldots,t_n$, and again the sum $s_1$ is not necessarily
zero.  $\W_ {E_ n}$ is generated by those reflections together
with $r_{v_0}$.  The action of $r_{v_0}$ on $(\H ^ {n+1})^*$
is not hard to compute:  we have
\[
r_ {v _ 0}(e_ i^*) =
\begin{cases}
e_ 0^* + (e_ 0^* + e_ 1^* + e_ 2^* + e_ 3^*) & \text{if $i=0$}, \\
e_ i^* - (e_ 0^* + e_ 1^* + e_ 2^* + e_ 3^*) & \text{if $1 \le i \le 3$}, \\
e_ i^* & \text{if $4 \le i \le n$}.
\end{cases}
\]
It follows that the action on $V_{E_n}^*$ is given by
\begin{equation}
\label{WE}
r_ {v _ 0}(t_ i) =
\begin{cases}
 t_ i - \frac 23
(t_ 1 + t_ 2 + t_ 3) & \text{if $1 \le i \le 3$}, \\
t_ i + \frac 13
(t_ 1 + t_ 2 + t_ 3) & \text{if $4 \le i \le n$}.
\end{cases}
\end{equation}

The invariant polynomials for the
action of $\W_ {E_ n}$ are
not as simple to describe as those for the actions of
$\W_{A_{n-1}}$ and $\W_{D_n}$.  In the reducible case $E_3$,
it is not hard to verify (using equation~(\ref{WE})) that the
polynomials
\begin{align}
\label{eqE3inv}
\begin{split}
 &\e_2^{(1)} \coloneq  s_1^2 \\
 &\e_2^{(2)} \coloneq  s_2 \\
 &\e_3       \coloneq  s_3 - \tfrac13 s_1 s_2 + \tfrac{2}{27} s_1^3
\end{split}
\end{align}
are invariant under $\W_{E_3}$, and generate the ring of invariants
$\C [V_{E_3}]^{\w_{E_3}}$.  In the cases of $E_4$ and $E_5$,
 bases for the rings of
invariants can be calculated by using the isomorphisms
$R_ {E_ 4} \isom R_ {A_ 4}$ and
$R_ {E_ 5} \isom R_ {D_ 5}$.
We carry this out in the corollary to lemma~\ref{lem71} in section 7.

The invariant polynomials are extremely complicated in the classical
cases $E_6$, $E_7$, and $E_8$.  Calculations of some of these were made
in the early 1950s by Coxeter \cite{[Cx]}, Frame \cite{[Fr]}, and
Racah \cite{[Rc]}, although published explicit formulas are limited
to the case of $E_6$.
One of the by-products
of the work discussed here is an
explicit description of a basis
for the invariant polynomials in these three cases.
(These invariants are of course polynomials in the symmetric functions
$s_1,\ldots,s_n$, since they are invariant under the subgroup
 $\W_ {A_ {n-1}}$.)  The basis we compute appears in Appendix 1 for
$E_6$ and in Appendix 2 for $E_7$.  The polynomials in the basis
are too large to be worth writing down
in the case of $E_8$, but an algorithm for computing with those polynomials
 is given in section 9.

\section{Rational double points and simultaneous resolution.}

Associated to each irreducible root system of type $S=A_{n-1}$, $D_n$,
 or $E_n$
is a rational double point with the same name.  The minimal
resolution of this complex
 surface singularity has a dual graph which
is isomorphic to
the corresponding Dynkin diagram, and the singularity is
determined up to isomorphism by the diagram.  For each singularity type,
we fix a representative for the isomorphism class:  the hypersurface
in $\C ^3$ whose defining polynomial is given in the middle column of
table~\ref{table2A}.  Notice
that as in section 2, our choices depend on the type $S$
 rather than just on the isomorphism class.

\begin{table}[t]
\begin{center}
\begin{tabular}{|c|c|rcl|} \hline
  &  & & & \\
$S$ & Defining Polynomial & \multicolumn{3}{c|}{Action of $\l  \in \C ^*$} \\
  &  & & & \\ \hline
  &  & & & \\
$A_{n-1}$ & $ - X Y + Z^n$ & $(X,Y,Z)$ & $\mapsto$ &
$(\l ^kX,\l ^{n-k}Y,\l Z)$ \\
$(n \ge 2)$  & & & & \\
  &  & & & \\
$D_n$  & $ - X^2 - Y^2 Z + Z^{n-1}$ & $(X,Y,Z)$ & $\mapsto$ &
$(\l ^{n-1}X,\l ^{n-2}Y,\l ^2Z)$\\
$(n \ge 3)$  &  & & & \\
  &  & & & \\
$E_4$   & $ - X Y + Z^5$ & $(X,Y,Z)$ & $\mapsto$ &
$(\l ^kX,\l ^{5-k}Y,\l Z)$\\
  &  & & & \\
$E_5$  & $ - X^2 - Y^2 Z + Z^4$ & $(X,Y,Z)$ & $\mapsto$ &
$(\l ^4X,\l ^3Y,\l ^2Z)$\\
  &  & & & \\
$E_6$  & $ - X^2 - X Z^2 + Y^3$ & $(X,Y,Z)$ & $\mapsto$ &
$(\l ^6X,\l ^4Y,\l ^3Z)$\\
  &  & & & \\
$E_7$  & $ - X^2 - Y^3 + 16 Y Z^3$ & $(X,Y,Z)$ & $\mapsto$ &
$(\l ^9X,\l ^6Y,\l ^4Z)$\\
  &  & & & \\
$E_8$  & $ - X^2 + Y^3 - Z^5$ & $(X,Y,Z)$ & $\mapsto$ &
$(\l ^{15}X,\l ^{10}Y,\l ^6Z)$\\
 & & & & \\ \hline
\end{tabular}
\end{center}

\medskip

\caption{}
\label{table2A}
\end{table}

Each of our chosen representatives admits a $\C ^*$-action, as
indicated in the last column of
table~\ref{table2A}:
if $\C ^3$ is given the specified $\C ^*$-action then the
defining polynomial is {\em weighted homogeneous}, i.e., is an
eigenfunction for the $\C ^*$-action.
In the cases of $A_{n-1}$ and $E_4$,
there are several possible $\C ^*$-actions,
any one of which will suit our purposes.

By a theorem of Pinkham \cite{C*},
a singularity with a $\C ^*$-action admits a $\C ^*$-semi-universal
deformation.  (Such a deformation is semi-universal for deformations
with a $\C ^*$-action, and is also semi-universal for arbitrary
deformations.)
For a surface in $\C ^3$ defined by a weighted homogeneous polynomial
$F$, this can be obtained as follows.  Choose
weighted homogeneous polynomials $G_1,\ldots,G_r$ which descend to
a basis
of the vector space
$\C [X,Y,Z]/(\frac{\partial F}{\partial X},
\frac{\partial F}{\partial Y}, \frac{\partial F}{\partial Z})$.
The polynomial
$\Phi \coloneq F + \mu_1 G_1 + \cdots + \mu_r G_r$
will then define a $\C ^*$-semi-universal deformation as a
hypersurface in $\C ^{3+r}$.  The variable $\mu_i$ is given
the difference of the weights of $F$ and $G_i$ as its weight.

In the case of rational double points, even after fixing the defining
polynomial as in table~\ref{table2A}, there are still some choices to
be made, for there may be more than one
choice of weighted homogeneous polynomials which give a basis of
$\C [X,Y,Z]/(\frac{\partial F}{\partial X},
\frac{\partial F}{\partial Y}, \frac{\partial F}{\partial Z})$.
We implicitly fix one such choice in
table~\ref{table2B}.

\begin{table}[t]
\begin{center}
$\begin{array}{|c|r@{\  + \ }l|} \hline
 & \multicolumn{2}{c|}{} \\
$S$ & \multicolumn{2}{c|}{\text{Preferred Versal Form}}  \\
 & \multicolumn{2}{c|}{} \\ \hline
 & \multicolumn{2}{c|}{} \\
A_{n-1}   &  - X Y + Z^n  &   \sum_{i=2}^{n}\a_i Z^{n-i} \\
(n \ge 2) & \multicolumn{2}{c|}{} \\
  & \multicolumn{2}{c|}{} \\
D_n  &  \multicolumn{1}{r@{\  - \ }}{ X^2 + Y^2 Z - Z^{n-1} } &
 \sum_{i=1}^{n-1}\d_{2i} Z^{n-i-1}  + 2 \c_n Y  \\
(n \ge 3)  & \multicolumn{2}{c|}{} \\
  & \multicolumn{2}{c|}{} \\
E_4   &  - X Y + Z^5  &   \e_2 Z^3 + \e_3 Z^2 + \e_4 Z + \e_5 \\
  & \multicolumn{2}{c|}{} \\
E_5  &  \multicolumn{1}{r@{\  - \ }}{  X^2 + Y^2 Z - Z^4 } &
 \e_2 Z^3 - \e_4 Z^2  + 2 \e_5 Y - \e_6 Z - \e_8  \\
  & \multicolumn{2}{c|}{} \\
E_6  &  - X^2 - X Z^2 + Y^3  &
      \e_2 Y Z^2 + \e_5 Y Z + \e_6 Z^2 + \e_8 Y \\
  &  &
      \e_9 Z + \e_{12} \\
  & \multicolumn{2}{c|}{} \\
E_7  &  - X^2 - Y^3 + 16 Y Z^3  &
    \e_2 Y^2 Z + \e_6 Y^2 + \e_8 Y Z + \e_{10} Z^2  \\
  &  &
       \e_{12} Y
     + \e_{14} Z + \e_{18} \\
  & \multicolumn{2}{c|}{} \\
E_8  &  - X^2 + Y^3 - Z^5 &
      \e_2 Y Z^3 + \e_8 Y Z^2 + \e_{12} Z^3 + \e_{14} Y Z \\
  &  &
      \e_{18} Z^2 + \e_{20} Y + \e_{24} Z
     + \e_{30} \\
 & \multicolumn{2}{c|}{} \\ \hline
\end{array}$
\end{center}

\medskip

\caption{}
\label{table2B}
\end{table}

\begin{definition}
The defining polynomial $\Phi_S$
of a $\C ^*$-semi-universal
deformation of a rational double point is said to be
in {\em preferred versal form\/} if it has
 the form given in table~\ref{table2B}.
\end{definition}

We have chosen the notation
 in table~\ref{table2B} to make the
$\C ^*$-action explicit:  the $\C ^*$-action on a
$\C ^*$-semi-universal deformation whose defining polynomial is in
preferred versal form is determined by giving
$X$, $Y$, $Z$ the same eigenvalues as in table~\ref{table2A},
and giving each coefficient
with subscript $i$ the eigenvalue $\l ^i$.

\begin{remark}
Tables~\ref{table2A} and \ref{table2B}
are the first of several places in this paper where the normalization
chosen in the case of $E_7$ looks a bit peculiar.  (The peculiarity
in this case is the coefficient ``$16$''.)  In all of these cases, our
notation is chosen to match that of Bramble \cite{[Bra]}, so that
we may directly compare our results with his.  The reader interested in
obtaining a defining polynomial of a $\C ^*$-semi-universal
deformation of  $E_7$ of the more natural form $- X^2 - Y^3 + Y Z^3
    +\e_2 Y^2 Z + \e_6 Y^2 + \e_8 Y Z + \e_{10} Z^2
     +  \e_{12} Y  + \e_{14} Z + \e_{18} $ may easily do so as follows. Start
with the defining polynomial of the $\C ^*$-semi-universal
deformation of  $E_7$ in preferred versal form as given by table~\ref{table2B}
and Appendix 2 (which will be calculated later).  Substituting $8X$ and $4Y$
for $X$ and $Y$ respectively, then dividing by 64 gives the desired defining
polynomial.

\end{remark}

Our first theorem is a slight refinement of the famous theorem on
simultaneous resolution of rational double points, due to
Brieskorn \cite{[Bri0]},
\cite{[Bri]},
\cite{[Bri-Nice]} and Tyurina \cite{[T]}.
(We follow Tyurina's approach, as amplified by Pinkham \cite{[P-RDP]}.)
Recall that a {\em simultaneous resolution\/} of a family
${\cal X} \to {\cal D}$ is a resolution of singularities of
${\cal X}$ which also resolves the singularities of each fiber of
the map ${\cal X} \to {\cal D}$.  In general,
one can only hope to find a simultaneous
resolution  after making a base change ${\cal R} \to {\cal D}$
and pulling back the original family to a family
${\cal X} \times_{{\cal D}} {\cal R}$.

Let $R$ be an irreducible root system of type $S=A_{n-1}$, $D_n$,
or $E_n$, in which a basis of simple roots $\{ v_i \}$
has been chosen. Let $V$ be the complex root space,  let
$\W$ be the Weyl group acting on $V$ by reflections,
and let $\rho\colon V \to V/\W$ be the quotient by $\W$.  Let $X_0$ be
the corresponding rational double point with $\C ^*$-action,
let $Z_0 \to X_0$ be the minimal resolution, and fix an identification
between the basis of simple roots $\{ v_i \}$
of $R$ and the irreducible exceptional
curves $\{ C_i \}$ on $Z_0$ which preserves the associated graphs.
For each positive root
$v = \sum a_i v_i$, the curve $C_v = \sum a_i C_i$ is an effective rational
curve of self-intersection $-2$ on $Z_0$.
The deformations of $Z_0$ to which $C_v$ lifts forms a natural subset
of all deformations of $Z_0$ (cf.\ \cite[\S 2]{[Wahl-disc]}).

We fix coordinates $X$, $Y$, $Z$ on $\C ^3$.

\begin{theorem}  
Let $R$ be an irreducible root system of type $S$, and let
$\rho \colon V \to V/\W$ and
$Z_0 \to X_0$ be as above.  Let $\C ^*$
act on the complex root space
$V$ via $v \mapsto \l  \cdot v$ for $\l  \in \C ^*$,
and give $V/\W$  the induced
$\C ^*$-action.  Then there is a $\C ^*$-semi-universal
deformation
${\cal X} \to V/\W$ of $X_0$ and a $\C ^*$-equivariant
simultaneous resolution ${\cal Z} \to {\cal X} \times_{\rho} V$
inducing $Z_0 \to X_0$ with the
following properties:
\begin{enumerate}
\item ${\cal X}$ can be embedded as
 a hypersurface in $\C ^3 \times V/\W$ whose
defining polynomial $\Phi_S$ with respect to the coordinates $X$, $Y$, $Z$
is in preferred versal form.

\item The coefficients of the defining polynomial $\Phi_S$ are explicitly
computable functions on $V/\W$, which give a basis of $\C [V]^{\w}$.

\item For each positive root $v \in R$, the set of deformations of $Z_0$
in the family ${\cal Z}$ to which the curve $C_v$ lifts is parametrized
by the hyperplane $v^{\perp}$ (the orthogonal complement of $v$ in $V$).
\end{enumerate}
\end{theorem}

The two novelties in this statement are the phrase ``explicitly
computable", and the third part of the theorem (although the latter is
implicit in Slodowy's exposition \cite{[Slod]}, and in some of Looijenga's
work \cite{[Lj]}).
We will give a proof of this theorem, based on two constructions of
Tyurina \cite{[T]},  in the next two sections.
We have explicitly computed  these coefficients as functions on $V/\W$,
and we give our results
 in equation (\ref{eqA}) for $A_{n-1}$, equation (\ref{eqD}) for
$D_n$, equations (\ref{eqE4}) and (\ref{eqE5}) for $E_4$ and $E_5$
respectively,
and Appendices 1 and 2 for $E_6$ and $E_7$ respectively.
(In the last two cases, we ``clear denominators'' before displaying
the result.)
We have not attempted to write down the formulas for $E_8$, although
we give an algorithmic method of computing with them in section 9.

\section{Simultaneous resolutions for $A_{n-1}$ and $D_n$.}

In this section, we will prove theorem 1 in the cases of $A_{n-1}$ and $D_n$.
It is easy to see that the theorem then follows in cases $E_4$ and $E_5$,
since the $\C ^*$-semi-universal deformation of $E_4$ (resp.~$E_5$)
in table~\ref{table2B} coincides with the one for $A_4$ (resp.~$D_5$).
We postpone the cases of $E_6$, $E_7$, and $E_8$ to the next section.

We will construct the required
deformations ${\cal X}$ of
rational double points and their simultaneous resolutions ${\cal Z}$,
 using the distinguished functionals
$t_1,\ldots,t_n$ (and the elementary symmetric functions
$s_1,\ldots,s_n$ in those distinguished functionals)
as the key ingredients in the construction.
The construction is due to Kas {\cite{[Kas]}} (for type $A_{n-1}$ only)
and Tyurina \cite{[T]}, both heavily influenced by work of Brieskorn
\cite{[Bri0]}.

\bigskip

We begin with the case of $A_{n-1}$.
Define functions $\a_i \in \C [V]^{\w}$ by
\begin{equation}
\label{eqA}
\a_i \coloneq
\begin{cases}
0 & \text{for $-1 \le i \le 1$}, \\
s_i & \text{for $\hphantom{-} 2 \le i \le n$}.
\end{cases}
\end{equation}
Since the Weyl group $\W$ coincides with the symmetric group
$\Sym_n$ and $s_1 = 0$ in $\C [V]$, it follows that
$\{ \a_i \}_{2 \le i \le n}$ is a basis of $\C [V]^{\w}$.

Define a hypersurface ${\cal X} \subset \C ^3 \times V/\W$ by means
of the polynomial
\[
\Phi_{A_{n-1}} \coloneq  - X Y + Z^n  + \sum_{i=2}^{n}\a_i Z^{n-i}
\]
which is in preferred versal form.
Notice that this
can also be written in the form
\[
\Phi_{A_{n-1}} = -XY + f_{A_{n-1}}(Z;t),
\]
and the functions $\a_i$ can be interpreted as the coefficients of
$f_{A_{n-1}}(Z;t)$.

Let $\rho\colon V \to V/\W$ be the quotient map.
To construct a simultaneous resolution of
${\cal X} \times_{\rho} V$, we recall that
in the ring $(\C [V])[U]$, the distinguished polynomial can be
written in the factored form given in equation (\ref{eqf1}).
Therefore,
\[\rho^*(\Phi_{A_{n-1}}) = -XY + \prod_{i=1}^n (Z + t_i).\]

We construct a simultaneous resolution
${\cal Z} \to {\cal X} \times_{\rho} V$ by
taking the closure of
the graph of the morphism
\begin{align*}
{\cal X} \times_{\rho} V &\to (\P
^ 1)^{n-1} \\
(X,Y,Z,t_ 1,\ldots ,t_ n) &\mapsto [X,\prod_{\nu=1}^
i(Z+t_\nu)]_ i
\end{align*}

Assertions 1 and 2 of theorem 1 are clear, so we need only
verify assertion 3.

Let $(u_ k,v_ k)$ be homogeneous coordinates on the
$k$th $\P ^ 1$ arising from the resolution described above.  Then
\begin{gather*}
-XY+\prod_{i=1}^ n(Z+t_ i)=0 \\
Xv_ j=u_ j\prod_{i=1}^ j(Z+t_ i)\quad  (1\le
j\le n-1), \\
\prod_{i=k+1}^ j(Z+t_ i)u_ jv_ k=u_ kv_ j\quad
(1\le k<j\le
n-1)
\end{gather*}
are equations for ${\cal Z}$.  Thinking of $V_{A_ {n-1}}$ as
a hyperplane in
$\C ^ n$, these equations show that the $\pi$-exceptional
part of
the fiber of ${\cal Z}$
over $(t_
1,\ldots ,t_ n)=(0,\ldots ,0)$ is given as the union of the
curves
$C_ i,\ 1\le i\le
n-1$, where $C_ i$ is defined by $X=Y=Z=0$, $u_ j=0$ for
$j<i$, $v_ k=0$
for $k>i$.

All positive roots are
of the form $v=v_ i+\cdots + v_{j-1}$.
{}From this, it is easy to see that the locus where $C_ v=C_
i+
\cdots +C_{j-1}$ deforms is given
by $t_ i=t_ j$,
which is indeed the orthogonal complement of $v$.  To see this,
note that if $t_ i=t_ j$, and are distinct from the other
$t_ k$, the
$\pi$-exceptional part of the fiber of ${\cal Z}$ over $(t_
1,\ldots ,
t_ n)$ is
given by $Z=-t_ i=-t_ j, u_ k=0\ (k<i),\ v_ k=0\ (k\ge
j),
u_ mv_ k\prod_{n=k+1}^ m(t_ n-t_ i)=u_ kv_ m,\
(i\le k<m<j)$.
For general parameter values, this is easily seen to be a $\P ^ 1$,
specializing to $C_ i+\cdots +C_{j-1}$ as the parameters
approach 0.  Since
the genus is constant, this is a flat family.  Since the general
fiber of this
family is irreducible, this is the only possible flat family over
the generic
point of $t_ i=t_ j$ which is contained in the exceptional
locus.

Since ${\cal Z}$
may be obtained from ${\cal X} \times_{\rho} V$ by successively blowing up
Weil divisors, the exceptional set lies entirely over the
discriminant
$\prod_{i<j}(t_ i-t_ j)$.  Hence the locus over which any
(union of)
exceptional curves  deform is contained in the locus where some
$t_ i$ equals
some $t_ j$.  This proves the assertion.

\bigskip

Next we turn to the case of $D_n$. Define functions
$\c_n, \d_{2i} \in \C [V]^{\w}$ by
\begin{align}
\label{eqD}
\begin{split}
\gamma_n &\coloneq t_1
\cdots
t_n  =  s_n \\
\delta_{2i} &\coloneq \text{the $i^{\text{th}}$ elementary symmetric}
 \\
 & \qquad \text{function of $ t_1^2, \ldots, t_n^2$}
\end{split}
\end{align}
(Note that $\d_{2n} = (\c_n)^2$.)  The functions
$\c_n, \d_2, \ldots, \d_{2n-2}$ generate $\C [V]^{\w}$.  Moreover,
comparing (\ref{eqD}) with the definition of the second distinguished
polynomial in (\ref{eqg1}), we see that
\begin{align*}
g_{D_n}(Z;t) &= Z^n + \sum_{i=1}^{n-1} \d_{2i} Z^{n-i} + (\c_n)^2 \\
f_{D_n}(0;t) &= \c_n.
\end{align*}

Define a hypersurface ${\cal X} \subset \C ^3 \times V/\W$ by means
of the polynomial
\[
\Phi_{D_{n}} \coloneq   X^2 + Y^2 Z   -
      Z^{n-1} - \sum_{i=1}^{n-1}\d_{2i} Z^{n-i-1}  + 2 \c_n Y
\]
which is in preferred versal form.
Notice that this
can also be written in the form
\[
\Phi_{D_{n}} =  X^2 + Y^2 Z - \frac{1}{Z}(g_{D_n}(Z;t)
- f_{D_n}(0;t)^2) + 2 f_{D_n}(0;t) Y.
\]

Let $\Sym_n \subset \W$ be the Weyl group of the subsystem of the
root system $R$ generated
by $v_1,\ldots,v_{n-1}$, which
acts on $V$ by permutations of the distinguished functionals
$t_1,\ldots,t_n$.  The ring of invariant functions
$\C [V]^{\sym_n}$ is generated by $s_1,\ldots,s_n$.

In the ring $(\C [V]^{\sym_n})[Z]$, we can use the properties of
$g_{D_n}(Z;t)$ given in equations (\ref{eqg3a}) and (\ref{eqg3b}),
which are expressed in terms of the polynomials $P_{D_n}(Z;t)$
and $Q_{D_n}(Z;t)$.  Now (\ref{eqg3a}) implies that $Q_{D_n}(0;t) = s_n$,
so we can also define
\begin{equation}
\label{eqS}
S_{D_n}(Z;t) \coloneq \frac{Q_{D_n}(Z;t) - s_n}{Z}.
\end{equation}

Let $\tau\colon V/\Sym_n \to V/\W$ be the natural map.
The first step in Tyurina's construction of the simultaneous resolution
is done most naturally on
${\cal X} \times_{\tau} V/\Sym_n$.
We abbreviate $P_{D_n}(Z;t)$, $Q_{D_n}(Z;t)$, and $S_{D_n}(Z;t)$ by
$P$, $Q$, and $S$, respectively.  Using (\ref{eqg3a}), (\ref{eqg3b}),
(\ref{eqS}) and some algebraic tricks as in Tyurina \cite{[T]},
we find that the defining polynomial of ${\cal X} \times_{\tau} V/\Sym_n$
can be written as
\begin{align*}
\tau^*(\Phi_{D_ n}) &=   X^2 + Y^2 Z
- \frac{1}{Z}(Z P^2 + (ZS + s_n)^2 - s_n^2)
+ 2 s_n Y \\
 &=  (X-P)(X+P)+(YZ+ZS+2s_n)(Y-S).
\end{align*}
Let ${\cal Y} \to {\cal X} \times_{\tau} V/\Sym_n$
be the blowup of the ideal
 $(X-P,Y-S)$.
It is easy to see that the entire singular locus of ${\cal Y}$
is contained in the coordinate chart ${\cal Y}^0$ defined by
the substitution $X - P = (Y - S) U$.  The defining polynomial in this chart
becomes
\begin{align*}
\Phi_{{\cal Y}^0} &=  U \cdot ( 2P+(Y-S)U ) + (YZ+ZS+2s_n) \\
 &=  (Y-S)(Z+U^2)+2UP+2Q \\
 &=  (Y-S+2G)(Z+U^2) + 2f_{D_n}(U;t)
\end{align*}
(using the polynomial $G = G_{D_n}(Z,U;t)$ defined in equation (\ref{eqG1})).

Now let $\sigma\colon V \to V/\Sym_n$ be the natural quotient map.
Since $f_{D_n}(U;t)$ factors in $(\C [V])[U]$, we have
\[
\sigma^*(\Phi_{{\cal Y}^0}) = (Y-S+2G)(Z+U^2) + 2 \prod_{i=1}^n (U+t_i),
\]
which after a change of coordinates is the same family used in the case
of $A_{n-1}$.
We can thus construct a resolution
${\cal Z} \to {\cal Y} \times_{\sigma} V$ (which also resolves
${\cal X} \times_{\rho} V$) by means of the construction used for $A_{n-1}$.

Again, we only need to verify assertion 3 from theorem 1.

Let $[u,v]$ be homogeneous coordinates on the $\P ^1$
of the first blow up for $D_n$ as described above, with $U=u/v$.
The homogeneous form of the equation
$\sigma^*(\Phi_{{\cal Y}^0})$ is
\[
(v^{n-2}(Y-S)+2{\widetilde
G})(v^2Z+u^2)+2\prod_{i=1}^n(u+t_iv)=0
\]
(using the homogeneous counterpart
${\widetilde G} = \widetilde{G}_{D_n}(Z,u,v;t)$
of $G$ which was defined in equation (\ref{eqG2})).
The $A_{n-1}$ singularity is then contained in the affine piece
$v=1$.
Let $[u_k,v_k]$ be homogeneous coordinates on the $k$th $\P ^1$ used in
resolving the $A_{n-1}$ singularity.  Then we get as equations
for ${\cal Z}_{D_n}$:
\begin{gather*}
(v^2Z+u^2)v_{n-1}=vu_{n-1}(u+t_nv) \\
v^{n-2-l}(v^2Z+u^2)v_l=u_l\prod_{i=l+1}^n(u+t_iv)\quad  (1\le
l\le n-2) \\
v^{k-j}v_ju_k=v_ku_j\prod_{i=j+1}^k(u+t_iv)\quad  (1\le j<k\le
n-1).
\end{gather*}
These equations show that the $\pi$-exceptional part of the fiber
of ${\cal Z}$
over $(t_1,\ldots ,t_n)=(0,\ldots ,0)$ is given as the union of
the curves
$C_i, 1\le i\le n,$ where for $i<n$, $C_i$ is defined by
$X=Y=Z=0$, $u=0$, $u_j=0$
for $j>i$, $v_j=0$ for $j<i$ and $C_n$ is defined by $X=Y=Z=0$,
$uv_{n-1}=vu_{n-1}$, $v^{n-l-2}v_l=u^{n-l-2}u_l$ for $1\le l\le
n-2$,
$v^{k-j}v_ju_k=v_ku_j\prod_{i=j+1}^k(u+t_iv)$  for $1\le j<k\le
n-1$.

The positive roots are of five types: $r=v_ i+\cdots +
v_{j-1}\quad (i\le j\le n)$, $s=v_n$, $t=v_j+\cdots +v_{n-2}+v_n
\quad (j\le n-2)$,
$ u=v_j+\cdots +v_n\quad (j\le n-2)$,
$ v=v_ j+\cdots + v_{k-1}+2v_ k+\cdots
+2v_{n-2}+v_{n-1}+v_
n\quad (j<k\le n-2)$.

Exactly as in the $A_{n-1}$ case, it is easy to see from the
above construction
that the locus where $C_ r$ deforms
is given by $t_ i=t_ j$.  For each of the remaining cases, we
exhibit the equations of a divisor in ${\cal Z}$ which defines a
flat family of generically irreducible curves over a hyperplane
in the
parameter space; the hyperplane will be given first.
\begin{gather*}
t_{n-1}+t_n=0 \tag*{{\bf s:}}\\
X=0,\ Y=-t_1\cdots t_{n-2},\ Z=-t_n^2 \\
(u-t_nv)v_{n-1}=vu_{n-1} \\
v_{n-2}=u_{n-2} \\
v^{n-l-2}v_l=u_l(u+t_{l+1}v)\cdots (u+t_{n-2}v)\quad (1\le l\le
n-2) \\
v^{k-j}v_ju_k=v_ku_j\prod_{i=j+1}^k(u+t_iv)\quad  (1\le j<k\le
n-1).
\end{gather*}
\begin{gather*}
t_j+t_n=0 \tag*{{\bf t:}}\\
X=0,\ Y=-t_1\cdots t_{j-1}t_{j+1}\cdots t_{n-1},\ Z=-t_n^2 \\
(u-t_nv)v_{n-1}=vu_{n-1} \\
(u-t_nv)v^{n-l-2}v_l=u_l(u+t_{l+1}v)\cdots (u+t_{n-1}v)\quad
(j+1\le l\le n-2) \\
v^{n-j-2}v_j=u_j(u+t_{j+1}v)\cdots (u+t_{n-1}v) \\
v^{n-l-2}v_l=u_l(u+t_{l+1}v)\cdots
(u+t_{j-1}v)(u+t_{j+1}v)\cdots
     (u+t_{n-1}v)\quad (1\le l\le j) \\
v^{k-j}v_ju_k=v_ku_j\prod_{i=j+1}^k(u+t_iv)\quad  (1\le j<k\le
n-1).
\end{gather*}
\begin{gather*}
t_j+t_{n-1}=0 \tag*{{\bf u:}} \\
X=0,\ Y=-t_1\cdots t_{j-1}t_{j+1}\cdots t_{n-2}t_n,\
Z=-t_{n-1}^2 \\
(u^2-t_{n-1}^2v^2)v_{n-1}=vu_{n-1}(u+t_nv) \\
(u-t_{n-1}v)v^{n-l-2}v_l=u_l(u+t_{l+1}v)\cdots
(u+t_{n-2}v)(u+t_nv)\quad
     (j+1\le l\le n-2) \\
v^{n-j-2}v_j=u_j(u+t_{j+1}v)\cdots (u+t_{n-2}v)(u+t_nv) \\
v^{n-l-2}v_l=u_l(u+t_{l+1}v)\cdots
(u+t_{j-1}v)(u+t_{j+1}v)\cdots
     (u+t_{n-2}v)(u+t_nv)\quad (1\le l\le j) \\
v^{k-j}v_ju_k=v_ku_j\prod_{i=j+1}^k(u+t_iv)\quad  (1\le j<k\le
n-1).
\end{gather*}
\begin{gather*}
t_j+t_k=0 \tag*{{\bf v:}}\\
X=0,\ Y=-t_1\cdots t_{j-1}t_{j+1}\cdots t_{k-1}t_{k+1}\cdots
t_n,\
Z=-t_k^2 \\
(u^2-t_k^2v^2)v_{n-1}=vu_{n-1}(u+t_nv) \\
(u^2-t_k^2v^2)v^{n-l-2}v_l=u_l(u+t_{l+1}v)\cdots (u+t_nv)\quad
     (k\le l\le n-2) \\
(u-t_kv)v^{n-j-1}v_{j-1}=u_{j-1}(u+t_{j+1}v)\cdots (u+t_nv) \\
(u-t_kv)v^{n-l-2}v_l=u_l(u+t_{l+1}v)\cdots
(u+t_{j-1}v)(u+t_{j+1}v)\cdots
     (u+t_nv)\quad (j\le l\le k-2) \\
v^{n-j-1}v_{j-1}=u_{j-1}(u+t_{j+1})\cdots
(u+t_{k-1})(u+t_{k+1})\cdots
     (u+t_n)
\end{gather*}
\begin{multline*}
v^{n-l-2}v_l=
u_l(u+t_{l+1}v)\cdots
(u+t_{j-1}v)(u+t_{j+1}v)\cdots \\
\cdots
     (u+t_{k-1}v)(u+t_{k+1}v)\cdots (u+t_nv)
\quad
(1\le l\le j-2)
\end{multline*}
\begin{gather*}
v^{k-j}v_ju_k=v_ku_j\prod_{i=j+1}^k(u+t_iv)\quad  (1\le j<k\le
n-1).
\end{gather*}

For each of the families described above in the four cases $w=s,t,u,v,$ the
following statement holds.  The general fiber is a $\P ^ 1$,
specializing to $C_w$ as the parameters approach 0.  Since
the genus is constant, this is a flat family.  Since the general
fiber of this
family is irreducible, this is the only possible flat family over
the generic
point of the corresponding hyperplane in the parameter space
which is contained
in the exceptional locus.

Since ${\cal Z}$
may be obtained from ${\cal X}_{D_n}$ by successively blowing up
Weil divisors, the exceptional set lies entirely over the
discriminant
$\prod_{i<j}(t_ i^2-t_ j^2)$.  Hence the locus over which
any (union of)
exceptional curves  deform is contained in the locus where some
$t_ i$ equals
plus or minus some $t_ j$.  In each case, this is the desired
orthogonal
complement.

\section{Simultaneous resolutions for $E_6$, $E_7$, and $E_8$.}

In this section, we complete the proof of theorem 1,
treating the cases of $E_6$, $E_7$, and $E_8$.  We use another
construction of Tyurina \cite{[T]}, one which had been anticipated by
Bramble \cite{[Bra]} in the case of $E_7$ in 1918.  This construction
is discussed in considerable detail by Pinkham \cite{[P-RDP]},
whose approach
we follow, and also by M\'erindol \cite{Mer}.  The strategy
in this case is to first build ${\cal Z}$, and then recover ${\cal X}$.
In fact, ${\cal Z}$ is constructed as an open subset of a relative
projective model $\bar{\cal Z} \to V$.

We begin with $\P ^2$
with homogeneous coordinates $[x, y, z]$, and let $C$ be the cuspidal
cubic with equation $x^3 = y z^2$.  The smooth points of this
 rational curve form an open set $C_0 \subset C$ isomorphic to the
affine line.  The map $\eta\colon \A ^1 \to \P ^2$ defined by
$\eta(U) \coloneq [U,U^3,1]$ gives such an isomorphism.

Let $V$ be the root space of $E_n$, and let the distinguished functionals
$t_1,\ldots,t_n$ serve as coordinates on $V$.
Let $\bar{\cal Z}_0 = \P ^2 \times V$,
and ${\cal C}_0 = C \times V$.
Define $n$ sections $\sigma_j\colon V \to \bar{\cal Z}_0$, $1 \le j \le n$
as follows:
\[\sigma_j(t_1,\ldots,t_n) \coloneq (\eta(t_j),(t_1,\ldots,t_n)).\]
Now for $j = 1,\ldots,n$, let $\bar{\cal Z}_j$ be the blowup of $\bar{\cal
Z}_{j-1}$
along the proper transform of the section $\sigma_j$, and let ${\cal C}_j$
be the proper transform of ${\cal C}_{j-1}$ on $\bar{\cal Z}_j$.  We use
$\bar{\cal Z}$ and ${\cal C}$ to denote $\bar{\cal Z}_n$ and ${\cal C}_n$,
respectively, and let $p\colon \bar{\cal Z} \to V$ be the natural map.  The
fibers
of $p$ are all smooth surfaces.

As Pinkham points out, each fiber $\bar{Z}_x \coloneq p^{-1}(x)$
of $p$ is the blowup
of $\P ^2$ in a collection of points $\eta(t_1),\ldots,\eta(t_n)$
in ``almost general position"
in the sense of Demazure \cite{[D-RDP]}.  In particular,
$\omega_{\bar{Z}_x}^{-1}$
is nef.  We use the notation $\omega_{\bar{\cal Z}/V}^{- k}$ to stand for
$(\omega_{\bar{\cal Z}/V}^{-1})^{\otimes k}$.  Let
\[\bar{\cal P} \coloneq \Proj _V \left( \bigoplus_{k \ge 0}
p_*(\omega_{\bar{\cal Z}/V}^{- k}) \right)\]
be the relative anti-canonical model.  The fibers of $\bar{\cal P} \to V$
are ``generalized del Pezzo surfaces", that is, del Pezzo surfaces with
rational double points allowed.

Let ${\cal Z} \coloneq \bar{\cal Z}  -  {\cal C}$ and
${\cal P} \coloneq \bar{\cal P} - {\cal C}$.
(We have abused notation, and denoted
the image of ${\cal C}$ in $\bar{\cal P}$ again by ${\cal C}$.)
As Pinkham shows, the Weyl group $\W = \W_{E_n}$ acts
on ${\cal P}$ by Cremona transformations (and permutations
of  $\{ \sigma_1,\ldots,\sigma_n \}$).  In Pinkham's version of the
construction,
the parameter space is $(C_0)^n$ rather than $V$; Pinkham computes the
induced action of $\W$ on $(C_0)^n$ obtaining the
formulas on p.~196 of \cite{[P-RDP]}.  Since those formulas
agree with our equation (\ref{WE}) which describes
 the action of $\W$ on $V$ when the distinguished
functionals are used as coordinates,\footnote{This is  why we
chose the distinguished functionals as we did.}
our identification of $V$ with
$(C_0)^n$ makes the map ${\cal P} \to V$ into a $\W$-equivariant map.
We define ${\cal X} \coloneq {\cal P}/\W.$

If we give $\P ^2$ and $V$ the linear $\C ^*$-actions defined by
\begin{align*}
[x,y,z] &\mapsto [\l x,\l^3y,z] \\
(t_1,\ldots,t_n) &\mapsto (\l t_1,\ldots,\l t_n)
\end{align*}
for $\l  \in \C ^*$,
then the entire construction becomes $\C ^*$-equivariant.
Pinkham shows that ${\cal X} \to V/\W$ is a $\C ^*$-semi-universal
deformation of the central fiber $X_0$, and that
${\cal Z} \to {\cal P} = {\cal X} \times_{V/\w} V$ is a
$\C ^*$-equivariant simultaneous resolution.

We need to verify the properties stated in theorem 1.  The third
property is the easiest this time, since it was essentially checked
by Pinkham.    For each $x \in V$, the singularities of the fiber
$\bar{P}_x$ correspond to the ``effective roots" of $\bar{Z}_x$.  These
can be seen from the geometry of the set of points blown up:  they
correspond to 2  points being infinitely near, 3 points being
collinear, 6 points being conconic, or (in the case of $E_8$) 8 points
lying on a nodal cubic with one of the points being the node.  For each
possible ``effective root", Pinkham computes the equation of the
locus in $V$ in which the root is effective, giving the results in
 a table on p.~193 of \cite{[P-RDP]}.  We display those equations in
table~\ref{table44}, using our identification of $(C_0)^n$ with $V$.  The other
column of the table gives the root $v$ in $V$ such that the equation
is proportional to the equation $v^{\perp} = 0$.

{\renewcommand{\arraystretch}{.6}

\begin{table}[t]
\begin{center}
\begin{tabular}{|c|c|} \hline
 &  \\
Equation & Root \\
 & \\ \hline
 & \\
$t_i - t_j = 0$ & $e_i - e_j$ \\
 & \\
$t_i + t_j + t_k = 0$ & $e_0 - e_i - e_j - e_k$ \\
 & \\
$\sum_{j=1}^6 t_{i_j} = 0$ & $2 e_0 - \sum_{j=1}^6 e_{i_j}$ \\
 & \\
$2 t_{i_1} + \sum_{j=2}^7 t_{i_j} = 0$ &
             $3 e_0 - 2 e_{i_1} - \sum_{j=2}^7 e_{i_j}$ \\
 & \\ \hline
\end{tabular}
\end{center}

\medskip

\caption{}
\label{table44}
\end{table}
}

To finish the proof of theorem 1, we must explain how to embed
${\cal X}$ into $\C ^3 \times V/\W$.  We will actually embed
${\cal P}$ into $\C ^3 \times V$, and then note the $\W$-invariance
of our construction.  And that embedding in turn will be a restriction
of a projective embedding of $\bar{\cal P}$.

To describe a projective embedding of $\bar{\cal P}$, we  extend
Demazure's analysis \cite{[D-RDP]}
of anti-pluricanonical mappings to the case of
families of surfaces. (A similar extension in another context has been
made by M\'erindol \cite{Mer}.) For a single generalized del Pezzo surface
$\bar{P}_x$, Demazure found the
following.\footnote{The weighted projective spaces occurring in the
description are only implicit in Demazure's paper.}
In the case of $E_6$, the anti-canonical map embeds  $\bar{P}_x$ into $\P ^3$,
and the image is a cubic surface with rational double points.
In the case of $E_7$, the anti-canonical map $\bar{P}_x \to \P ^2$ is
a finite map of degree 2, and the anti-bicanonical map embeds $\bar{P}_x$
into the weighted projective space $\P ^{(1,1,1,2)}$.
$\bar{P}_x$ can be described as a double cover of $\P ^2$ ramified along
a quartic curve.  In the case of $E_8$, the anti-canonical system is a
pencil with a base point, and the anti-bicanonical map is again of degree
2, this time mapping to the weighted projective space
$\P ^{(1,1,2)}$ (which can be embedded as a
quadric cone in $\P ^3$).  $\bar{P}_x$ is a double cover of
$\P ^{(1,1,2)}$, branched on
a curve of
graded\footnote{We use the term {\em graded degree\/} (rather than
{\em weighted degree\/}) for a polynomial
in a weighted projective space, to avoid confusion with the weights
under the background $\C ^*$-action.}
degree 6; the anti-tricanonical
map embeds $\bar{P}_x$ into $\P ^{(1,1,2,3)}$.

Our immediate goal is to show that these projective
embeddings can be described {\em globally\/} over $V$
as embeddings into $\P  \times V$ for the appropriate weighted
projective space $\P $.
Then we will show that
the image in the case of $E_6$, and branch loci in the cases of $E_7$ and
$E_8$, can be similarly globally defined over $V$ by a single polynomial.
We will explain how to explicitly compute these polynomials in section 9;
here we simply show that they exist.
A related approach to constructing
projective embeddings of $\bar{\cal P}$ appears in \cite{Mer}.

\begin{lemma} \label{lem55}
\quad
\begin{enumerate}
\item
The sheaf $p_*(\omega_{\bar{\cal Z}/V}^{- k})$ is locally
free of rank
 $h^0(\bar{Z}_x,\omega_{\bar{Z}_x}^{- k})
= 1 + \frac{(k^2+k)(9-n)}{2}$ in case $E_n$.

\item
In cases $E_7$ and $E_8$, the map
$\Symm ^2(p_*\omega_{\bar{\cal Z}/V}^{-  1})\to
p_*\omega_{\bar{\cal Z}/V}^{-  2}$ is injective as a morphism of vector
bundles, and its cokernel is locally free of rank 1.

\item
In the case of $E_8$, the natural map
$p_*\omega_{\bar{\cal Z}/V}^{-  1}\otimes
p_*\omega_{\bar{\cal Z}/V}^{-  2} \to
p_*\omega_{\bar{\cal Z}/V}^{-  3}$ has a cokernel which is
locally free of rank 1.
Its kernel coincides with
$\Ker (p_*\omega_{\bar{\cal Z}/V}^{-  1}\otimes
\Symm ^2(p_*\omega_{\bar{\cal Z}/V}^{-  1}) \to
\Symm ^3(p_*\omega_{\bar{\cal Z}/V}^{-  1}))$,
which is locally free of rank 2.

\end{enumerate}
\end{lemma}

\begin{pf}

(1)  Since
 $\omega_{\bar{Z}_x}^{-1}$
is nef,  $H^1(\bar{Z}_x,\omega_{\bar{Z}_x}^{- k}) = 0$ for all $k \ge 0$.
Thus, each of the sheaves $p_*(\omega_{\bar{\cal Z}/V}^{- k})$ is locally
free.  The rank  follows from Riemann-Roch, since $c_1^2 = 9-n$ for $E_n$.

(2)  According to Demazure \cite[V, Proposition 1b]{[D-RDP]}, the anticanonical
mapping in the case
of $E_7$ is 2-1.  In particular, its image is not contained in any quadric,
hence the desired injectivity statement.  In the case of $E_8$, the
anticanonical map maps to $\P ^1$, hence its image is not contained in a
quadric either.  The rank 1 assertion follows from the ranks listed in the
table of the lemma.

(3)    By Demazure [op.\ cit.], in the case
of $E_8$ the antibicanonical mapping is 2-1; in particular, its image is not
contained in a graded cubic hypersurface (thinking of the mapping as
factoring through the weighted projective space $\P ^{(1,1,2)}$).  This
gives the equality of the two kernels, where
$\hbox{Sym}^2p_*\omega_{\bar{\cal Z}/V}^{-  1}$ is identified with its
image in $p_*\omega_{\bar{\cal Z}/V}^{-  2}$ by part 2 already proven.
The second mentioned kernel is easily computed to have rank 2 (the mapping is
surjective); from which it follows that the cokernel of the first mapping has
rank 1.
\end{pf}

\begin{lemma} \label{lem56}
There exist $\C ^*$-invariant subspaces
\begin{align*}
&L_1 \subset H^0(V,p_*(\omega_{\bar{\cal Z}/V}^{-  1})) \text{ of dimension
$10-n$,
in case $E_n$}, \\
&L'_2 \subset H^0(V,p_*(\omega_{\bar{\cal Z}/V}^{-  2})) \text{ of dimension 1,
in cases $E_7$ and $E_8$, and} \\
&L'_3 \subset H^0(V,p_*(\omega_{\bar{\cal Z}/V}^{-  3})) \text{ of dimension 1
in case $E_8$,}
\end{align*}
such that the natural maps
\begin{align*}
L_1 \otimes {\cal O}_V &\to  p_*(\omega_{\bar{\cal Z}/V}^{-  1}) \\
L'_2 \otimes {\cal O}_V &\to
\Coker (\Symm ^2(p_*\omega_{\bar{\cal Z}/V}^{-  1})\to
p_*\omega_{\bar{\cal Z}/V}^{-  2}) \\
L'_3 \otimes {\cal O}_V &\to
\Coker (
p_*\omega_{\bar{\cal Z}/V}^{-  1}\otimes
p_*\omega_{\bar{\cal Z}/V}^{-  2} \to
p_*\omega_{\bar{\cal Z}/V}^{-  3})
\end{align*}
are $\C ^*$-equivariant isomorphisms of sheaves.  In other words,
the targets of these maps are {\em trivial\/} locally
free sheaves.
\end{lemma}

\begin{pf}
The key ingredient is Quillen's affirmative answer \cite{[Q]}
to Serre's
conjecture that a finitely generated projective module over a polynomial ring
must be free.  Since $p_*\omega_{\bar{\cal Z}/V}^{-  1}$ is coherent and $V$
is affine,
$M=H^0(V,p_*\omega_{\bar{\cal Z}/V}^{-  1})$ is a finitely generated module
over the polynomial algebra $\C [V]$.  Since $M$ is also locally free, it is
in fact projective.  (Cf.\  \cite[Chapter II, \S 5.2, Theorem 1]{[CA]}.)
Quillen's theorem
then implies that $M$ is free.  This proves the triviality of the target of the
first map.  The other two cases are similar; just observe that the sheaves on
the right hand side of the
arrows are coherent and locally free of rank 1
by lemma~\ref{lem55}.
Since all maps in question are $\C ^*$-equivariant, we may choose
$\C ^*$-eigensections as generators of these trivial bundles.
\end{pf}

We let $L = \bigoplus L_k$ be the graded
$\C [V]$-algebra generated by $L_1$ in the case of
$E_6$, by $L_1$ and $L'_2$ in the case of $E_7$, and by
$L_1$, $L'_2$, and $L'_3$
in the case of $E_8$.  (We have $L_2 = \Symm ^2L_1 \oplus L'_2$
in cases $E_7$ and $E_8$, and
$L_3 = \Symm ^3L_1 \oplus (L_1 \otimes L'_2) \oplus L'_3$ in case $E_8$.)
The algebra $L$ has 4 generators, and
 $\Proj _V(L)$
is a relative weighted projective space of dimension 3 over $V$,
 isomorphic to
$\P ^{(1,1,1,1)} \times V$, $\P ^{(1,1,1,2)} \times V$,
 or $\P ^{(1,1,2,3)} \times V$,
respectively.  Moreover, there is a natural embedding of $\bar{\cal P}$
into $\Proj _V(L)$, where it forms a hypersurface.

We now show that $\bar{\cal P}$ can be defined by a weighted homogeneous
polynomial globally
over $V$.

\begin{lemma} \label{lem57}
There is a weighted homogeneous polynomial $\Phi_{E_n} \in L_k$ which
generates the ideal of $\bar{\cal P}$ in $\Proj _V(L)$,
where $k=3,4,6$ in cases $E_6$, $E_7$, $E_8$, respectively.
In the $E_7$ case, by an appropriate choice of $L'_2$ and a generator
$X_2$ of $L_2'$,
the polynomial takes the form $X_2^2-f_4$, for some
$f_4\in\hbox{Sym}^4(L_1)\subset
L_4$.
In the $E_8$ case, by an appropriate choice of $L'_3$ and a generator
$X_3$ of $L_3'$,
the polynomial takes the form
 $X_3^2-f_6$, for some $f_6\in\hbox{Sym}^6(L_1)\oplus (\hbox{Sym}^4(L_1)\otimes
L_2')\oplus
(\hbox{Sym}^2(L_1)\otimes\hbox{Sym}^2(L_2'))\oplus\hbox{Sym}^3(L_2')\subset
L_6$.
\end{lemma}

\begin{pf}
By \cite[V.3]{[D-RDP]}, it follows in the case of
$E_6$ that each
$\bar{P}_x$ is a cubic.  Hence the sheaf Ker$(\hbox{Sym}^3
p_*(\omega_{\bar{\cal Z}/V}^{-  1})\to p_*(\omega_{\bar{\cal Z}/V}^
{-  3}))$ is invertible.  Since it is also coherent, being the kernel of
a morphism of coherent sheaves, its triviality follows from Quillen's theorem.
This yields a  cubic  with coefficients in $\C [V]$
which serves as a global defining polynomial for $\bar{\cal P}$.

In the case of $E_7$, pick a generator $X_2$ for $L'_2$, and define a
map
\[p_*(\hbox{Sym}^2\omega_{\bar{\cal Z}/V}^{-  1}) \oplus
p_*(\hbox{Sym}^4\omega_{\bar{\cal Z}/V}^{-  1}) \to
p_*(\omega_{\bar{\cal Z}/V}^{-  4})\]
by $(a,b) \mapsto aX_2+b$.
It follows from \cite[V.4]{[D-RDP]} that
this map is surjective on each fiber; by the triviality of the bundles,
it must be surjective on global sections as well.  Thus, since
$(X_2)^2 \in H^0(V,p_*(\omega_{\bar{\cal Z}/V}^{-  4}))$,
there exist global sections
$a\in H^0(V,p_*(\hbox{Sym}^2\omega_{\bar{\cal Z}/V}^{-  1})),$
$b\in H^0(V,p_*(\hbox{Sym}^4\omega_{\bar{\cal Z}/V}^{-  1}))$
such that
$(X_2)^2=aX_2+b$.  Since $p_*(\omega_{\bar{\cal Z}/V}^{-  1})$ is trivial,
it follows that $H^0(V,p_*(\hbox{Sym}^k\omega_{\bar{\cal Z}/V}^{-  1}))=
\hbox{Sym}^kH^0(V,p_*(\omega_{\bar{\cal Z}/V}^{-  1}))$, hence
$a$ and $b$ can be described as linear combinations of polynomials on $V$ times
monomials in a basis for $L_1$.  We thus get
$\bar{\cal P}$ as  defined by the graded quartic polynomial
$\Phi_{E_7} \coloneq -(X_2)^2+aX_2+b$
 with coefficients in
$\C [V]$.  To put the defining polynomial into the form claimed,
we need only complete the square as in
\cite[V.4]{[D-RDP]}.  This is tantamount to making a different choice for
$L_2'$.

The case of $E_8$ is similar.

Since  $\bar{\cal P}$ is
 invariant under the $\C ^*$-action, the defining polynomial $\Phi_{E_n}$ must
be weighted homogeneous.
\end{pf}

In order to calculate explicit generators for the algebra $L$,
we interpret global sections of
$p_*(\omega_{\bar{\cal Z}/V}^{-  k})$ as being the defining polynomials
of curves of degree $3k$ in
$\P ^2 \times V$ with base conditions.  The base conditions
state that the curve should pass through the zero-cycle
$\eta(t_1)+\cdots+\eta(t_n)$ with multiplicity
$k$.  Choosing a basis for the space of such polynomials determines
a rational map $\pi\colon \P ^2 \times V \to \P ^N$, which can be
interpreted as the blowup of $\sigma_1,\ldots,\sigma_n$ followed
by the anti-$k$-canonical map of ${\cal Z}$.

To guarantee the $\W$-invariance of our defining polynomial $\Phi_{E_n}$
for $\bar{\cal P}$
(and hence obtain a defining polynomial for ${\cal X}$) we must choose our
generators of $L$ quite carefully.

\begin{definition}
Let $\mm  = (t_1,\ldots,t_n) \subset \C [V]$ be the maximal
ideal of the origin $0 \in V$.
We say that the weighted homogeneous polynomials
$X, Y, Z, W \in \C [V][x,y,z]$ form a {\em good generating set\/}
for $L$ if they generate $L$ and satisfy the defining conditions given in
table~\ref{table45}.
\end{definition}

The defining conditions given in table~\ref{table45} were obtained as follows.
Each element of $L_k$ when restricted to the central fiber
$\bar{Z}_0$ becomes an
element of $H^0(\bar{Z}_0,\omega_{\bar{Z}_0}^{-  k})$.
In the case of $E_n$, the polynomials of
degree $3k$ which belong to that space are exactly the ones which have a zero
of order at least $nk$ at $U=0$ when they are pulled back via
$\eta\colon C_0 \to \P ^2$,
and whose partial derivatives with respect to $x$, $y$, and $z$
up through order $k-1$ have a zero of order at least $k$ when
{\em they\/} are pulled back via $\eta$.  It is easy to see that these
conditions on the partial derivatives are superfluous when applied
to a monomial, since $n \ge 6$, which is more than the largest order of
vanishing at $U=0$ among $x$, $y$, and $z$ after pulling back via $\eta$.
A set of generators for the anti-pluricanonical ring of $\bar{Z}_0$ is
easy to find using this description.  We have implicitly listed one such set
in
table~\ref{table45}, as the right-hand sides of  congruences for $X$, $Y$,
$Z$, and $W$.
(The last column of the table includes congruences for $X$ derived from
the defining conditions in the cases of $E_7$ and $E_8$.)
The $\C ^*$-action preserves the central fiber $\bar{Z}_0$, so we
chose a generating set of polynomials which are weighted homogeneous on the
central fiber.\footnote{The coefficients in the case of $E_7$ were
chosen to match the notation of Bramble \cite{[Bra]}.}

We can see the weights of the generators from the information
given in table~\ref{table45}:  they are $(9, 7, 6, 3)$ for
 $({X},{Y},{Z},{W})$ in the case of $E_6$,
$(15, 9, 7, 3)$ for $({X},{Y},{Z},{W})$ in the case of $E_7$,
and $(24, 16, 9, 3)$ for $({X},{Y},{Z},{W})$
in the case of $E_8$.
Moreover, the last column of the table shows the leading order terms
of the polynomial $\Phi_{E_n}$ (determined by elimination from the defining
conditions), and
determines its weight as 21, 30, or 48, respectively.

{\renewcommand{\arraystretch}{.6}

\begin{table}[t]
\begin{center}
\begin{tabular}{|c|rcl|c|} \hline
 & & & & \\
\multicolumn{1}{|c|}{} & \multicolumn{3}{c|}{Defining Conditions} &
\multicolumn{1}{c|}{Other Properties} \\
 & & & & \\ \hline
 & & & & \\
 & ${W}$ & $=$ & $x^3 - y z^2$ & \\
 & & & & \\
 & ${Z}$ & $\equiv$ & $y^2 z\mod \mm $  & \\
$E_6$ & & & &
${Y}^3 - {X} {Z}^2 - {X}^2 {W} \equiv 0\mod \mm $  \\
 & ${Y}$ & $\equiv$ & $x y^2\mod \mm $  & \\
 & & & & \\
 & ${X}$ & $\equiv$ & $y^3\mod \mm $  & \\
 & & & & \\ \hline

 & & & & \\
 & ${W}$ & $=$ & $x^3 - y z^2$ &  \\
 & & & & \\
 & ${Z}$ & $\equiv$ & $x y^2\mod \mm $  &
${X} \equiv 8 y^5 z\mod \mm $  \\
$E_7$ & & & & \\
 & ${Y}$ & $\equiv$ & $4 y^3\mod \mm $  &
$16 {Y} {Z}^3 - {X}^2 - {Y}^3 {W} \equiv 0\mod \mm $  \\
 & & & & \\
 & ${X}$ & $=$ &  $\frac13 \
    \frac{\partial({Y},{Z},{W})}{\partial(x,y,z)}$ & \\
 & & & & \\ \hline

 & & & & \\
 & ${W}$ & $=$ & $x^3 - y z^2$ & \\
 & & & & \\
 & ${Z}$ & $\equiv$ & $y^3\mod \mm $  &
${X} \equiv y^8 z\mod \mm $  \\
$E_8$ & & & & \\
 & ${Y}$ & $\equiv$ & $x y^5\mod \mm $  &
${Y}^3 - {X}^2 - {Z}^5 {W} \equiv 0\mod \mm $  \\
 & & & & \\
 & ${X}$ & $=$ &  $- \frac16 \
    \frac{\partial({Y},{Z},{W})}{\partial(x,y,z)}$ & \\
 & & & & \\ \hline
\end{tabular}
\end{center}

\medskip

\caption{}
\label{table45}
\end{table}
}

\begin{proposition} \label{prop51}
There exists a good generating set $X, Y, Z, W$ for $L$ such that when the
defining polynomial $\Phi_{E_n} \in \C [V][X,Y,Z,W]$
of $\bar{\cal P} \subset \Proj _V(L)$ is restricted to the affine
chart $W=1$,
it gives
a $\C ^*$-semi-universal deformation of the rational double point
which is in preferred versal form.
\end{proposition}

\begin{pf}
We first claim that there exist good generating sets for
$L$.  We can start with the polynomial $x^3-yz^2$ as one of the generators,
since it belongs to
$L$ in all three cases.  On the central fiber $\bar{Z}_0$, the generator
$x^3-yz^2$ can be
extended to a generating set for the entire algebra
 $\bigoplus H^0(\bar{Z}_0,\omega_{\bar{Z}_0}^{-  k})$
as indicated in table~\ref{table45}:  the additional generators are
 $(y^2z, xy^2, y^3)$ in
case $E_6$, $(xy^2, 4y^3, 8y^5z)$ in case $E_7$, and
$(y^3,xy^5,y^8z)$ in case $E_8$.

Since each of the bundles involved is trivial, these generators on the
central fiber can be lifted
to generators of $L$ itself.  The only thing left to show is that in the
cases of $E_7$ and $E_8$, we may use the Jacobian determinant as
the generator of top degree.  To see this, we only need to note that
this Jacobian determinant does indeed belong to the algebra (since
it satisfies the base conditions), and its restriction mod $\mm $
is $8y^5z$, resp.\ $y^8z$, as required.

Now let $\bar{X}, \bar{Y}, \bar{Z}, \bar{W}$ be a good generating set for
$L$, and let $\bar{\Phi}_{E_n}$ be the defining polynomial of $\bar{\cal P}$
with respect to this generating set.  In cases $E_7$ and $E_8$, the
map determined by $\bar{Y}, \bar{Z}, \bar{W}$ expresses $\bar{\cal P}$
as a double cover of a weighted projective space.
Since $\bar{X}$ is the Jacobian determinant of
this map (up to a constant), it vanishes exactly on the ramification
locus of this double cover.  It follows that $\bar{\Phi}_{E_n}$ takes the
form $-\bar{X}^2 + \widetilde{\Phi}_{E_n}(\bar{Y},\bar{Z},\bar{W})$ in these
two cases.

In all three cases, we know the weight $w$ of $\bar{\Phi}_{E_n}$,
and also its graded degree $d$ in the  algebra $L$.  Any monomial
$m$ in $\bar{X}, \bar{Y}, \bar{Z}, \bar{W}$ which appears in $\bar{\Phi}_{E_n}$
must
have graded degree $d$, and weight $w_m \le w$:  the weight of the
coefficient of the monomial will then be $w-w_m$.  Using
table~\ref{table45} to determine the leading order terms,
we can then write   $\bar{\Phi}_{E_n}$ as follows, with undetermined
coefficients.\footnote{The two strange terms in the expression for
$\bar{\Phi}_{E_7}$ are yet another artifact of making our notation
match that of Bramble \cite{[Bra]}.}
(We adopt the convention that equation numbers which are followed by
$a$, $b$, or $c$ refer to the cases of $E_6$, $E_7$, or $E_8$, respectively.)
The notation is chosen so that the subscript on a coefficient shows
its weight.

\refstepcounter{equation}\label{eq11}
\begin{align*}
\begin{split}
\bar{\Phi}_{E_6} &=
- \bar{X}^2 \bar{W}
- \bar{X} \bar{Z}^2
+ \bar{Y}^3
+ \bar{\f}_1 \bar{Y}^2 \bar{Z}
+ \bar{\f}_2 \bar{X} \bar{Y} \bar{W}
+ \bar{\e}_2 \bar{Y} \bar{Z}^2
\\ &\qquad
+ \bar{\f}'_3 \bar{X} \bar{Z} \bar{W}
+ \bar{\f}''_3 \bar{Z}^3
+ \bar{\f}_4 \bar{Y}^2 \bar{W}
+ \bar{\e}_5 \bar{Y} \bar{Z} \bar{W}
+ \bar{\f}_6 \bar{X} \bar{W}^2
\\ &\qquad
+ \bar{\e}_6 \bar{Z}^2 \bar{W}
+ \bar{\e}_8 \bar{Y} \bar{W}^2
+ \bar{\e}_9 \bar{Z} \bar{W}^2
+ \bar{\e}_{12} \bar{W}^3
\end{split}\displaybreak[0]\tag{\ref{eq11}a}\\[1.5ex]
\begin{split}
\bar{\Phi}_{E_7} &=
- \bar{X}^2 - \bar{Y}^3 \bar{W}
+ 16 \bar{Y} \bar{Z}^3
+ \bar{\e}_2 \bar{Y}^2 \bar{Z} \bar{W}
+ \bar{\f}_2 (16 \bar{Z}^4 -
\makebox[0pt][l]{$
\bar{Y}^2 \bar{Z} \bar{W})
$}
\\ &\qquad
+ \bar{\f}_4 \bar{Y} \bar{Z}^2 \bar{W}
+ \bar{\e}_6 \bar{Y}^2 \bar{W}^2
+ \bar{\f}_6 (16 \bar{Z}^3 \bar{W} - \bar{Y}^2 \bar{W}^2)
\\ &\qquad
+ \bar{\e}_8 \bar{Y} \bar{Z} \bar{W}^2
+ \bar{\e}_{10} \bar{Z}^2 \bar{W}^2
+ \bar{\e}_{12} \bar{Y} \bar{W}^3
+ \bar{\e}_{14} \bar{Z} \bar{W}^3
+ \bar{\e}_{18} \bar{W}^4
\end{split}\displaybreak[0]\tag{\ref{eq11}b}\\[1.5ex]
\begin{split}
\bar{\Phi}_{E_8} &=
- \bar{X}^2 + \bar{Y}^3
- \bar{Z}^5 \bar{W}
+ \bar{\e}_2 \bar{Y} \bar{Z}^3 \bar{W}
+ \bar{\f}_4 \bar{Y}^2 \bar{Z} \bar{W}
+ \bar{\f}_6 \bar{Z}^4 \bar{W}^2
\\ &\qquad
+ \bar{\e}_8 \bar{Y} \bar{Z}^2 \bar{W}^2
+ \bar{\f}_{10} \bar{Y}^2 \bar{W}^2
+ \bar{\e}_{12} \bar{Z}^3 \bar{W}^3
+ \bar{\e}_{14} \bar{Y} \bar{Z} \bar{W}^3
\\ &\qquad
+ \bar{\e}_{18} \bar{Z}^2 \bar{W}^4
+ \bar{\e}_{20} \bar{Y} \bar{W}^4
+ \bar{\e}_{24} \bar{Z} \bar{W}^5
+ \bar{\e}_{30} \bar{W}^6
\end{split}\tag{\ref{eq11}c}
\end{align*}

Notice that this polynomial restricts to one in preferred versal form
in the affine
$W=1$  exactly when all of the $\bar{\f}_i$
terms vanish.  We therefore wish to make a change of generating set
to ensure that this occurs.

Suppose that $X, Y, Z, W$ is another good generating set of
$L$.  Since the mod $\mm $ restriction of a good generating set is
determined by table~\ref{table45},
the change of generators must restrict to the identity mod $\mm$.
Moreover (as follows from properties of the Jacobian determinant),
 in cases $E_7$ and $E_8$ we have $\bar{X}=X$.  Considering as
before the graded degrees of elements of the  algebra $L$,
and the fact that each term in the equation expressing the
change of generators
must have a coefficient with nonnegative weight, we find that
the change of generators must take the form shown below, with
undetermined coefficients $\psi_i$.

\refstepcounter{equation}\label{eq23}
\begin{align*}
\begin{split}
\begin{tabular}{ccrcrcrcr}
$\bar{X}$ & $=$ & $X$ & + & $\psi_2 Y$ & + & $\psi'_3 Z$ & + & $\psi_6 W$ \\
$\bar{Y}$ & $=$ &   &  & $Y$ & + & $\psi_1 Z$ & + & $\psi_4 W$ \\
$\bar{Z}$ & $=$ &   &  &   & & $Z$ & + & $\psi''_3 W$ \\
$\bar{W}$ & $=$ &   &  &   & &   &  & $W$ \\
\end{tabular}
\end{split}\displaybreak[0]\tag{\ref{eq23}a}\\[1.5ex]
\begin{split}
\begin{tabular}{ccrcrcrcr}
$\bar{X}$ & $=$ & $X$ &  &  &  &  &  &  \\
$\bar{Y}$ & $=$ &   &  & $Y$ & + & $\psi_2 Z$ & + & $\psi_6 W$ \\
$\bar{Z}$ & $=$ &   &  &   & & $Z$ & + & $\psi_4 W$ \\
$\bar{W}$ & $=$ &   &  &   & &   &  & $W$ \\
\end{tabular}
\end{split}\displaybreak[0]\tag{\ref{eq23}b}\\[1.5ex]
\begin{split}
\begin{tabular}{ccrcrcrcr}
$\bar{X}$ & $=$ & $X$ &  &  &  &  &  &  \\
$\bar{Y}$ & $=$ &   &  & $Y$ & + & $\psi_4 Z W$ & + & $\psi_{10} W^2$ \\
$\bar{Z}$ & $=$ &   &  &   & & $Z$ & + & $\psi_6 W$ \\
$\bar{W}$ & $=$ &   &  &   & &   &  & $W$ \\
\end{tabular}
\end{split}\tag{\ref{eq23}c}
\end{align*}

To finish the proof,  we substitute
equation (\ref{eq23})
into  equation (\ref{eq11})
and collect the coefficients of the monomials
$\{Y^2 Z$, $XYW$, $ Z^3$, $XZW$, $Y^2 W$, $X W^2\}$ in case $E_6$,
$\{Z^4$, $Y Z^2 W$, $Z^3 W\}$ in case $E_7$, and
$\{Y^2 Z W$, $Z^4 W^2$, $Y^2 W^2\}$ in case $E_8$.
(These are the monomials which we cannot allow if we wish to achieve
preferred versal form.)  Equating all such coefficients to zero gives
the following set of equations.

\refstepcounter{equation}\label{eq51}
\begin{align*}
\begin{split}
\begin{aligned}
0 &= { { 3 {\psi_{1}}}+ {\bar{\f}_{1}}}  \\
0 &= {{- 2 {\psi_{2}}} +{\bar{\f}_{2}} }  \\
0 &= { {- {\psi'_{3}}} + {\bar{\f}''_{3}} +  { {\bar{\e}_{2}} {\psi_{1}}}
+ { {\psi_{1}}^{3}}  +  { {\bar{\f}_{1}} { {\psi_{1}}^{2}}}}  \\
0 &= { {- 2 {\psi'_{3}}} {- 2 {\psi''_{3}}} +
{\bar{\f}'_{3}} +  { {\bar{\f}_{2}} {\psi_{1}}} }  \\
0 &= { 3 {\psi_{4}}}+ { { {\bar{\f}_{1}} {\psi''_{3}}} { -  { {\psi_{2}}^{2}}}
+  { {\bar{\f}_{2}}
{\psi_{2}}}+ {\bar{\f}_{4}}}  \\
0 &= {{- 2 {\psi_{6}}} + { {\bar{\f}_{2}} {\psi_{4}}}
+  { {\bar{\f}'_{3}} {\psi''_{3
}}} { -  { {\psi''_{3}}^{2}}}+ {\bar{\f}_{6}}}
\end{aligned}
\end{split}\displaybreak[0]\tag{\ref{eq51}a}\\[3ex]
\begin{split}
\begin{aligned}
0 &= { { 16 {\psi_{2}}}+ { 16 {\bar{\f}_{2}}}}  \\
0 &= { { 48 {\psi_{4}}} { -  3 { {\psi_{2}}^{2}}}
+ {\bar{\f}_{4}} +  { 2 {\bar{\e}
_{2}} {\psi_{2}}} { -  2 {\bar{\f}_{2}} {\psi_{2}}}}  \\
0 &= { { 16 {\psi_{6}}} + { 16 {\bar{\f}_{6}}} { -  { {\psi_{2}}^{3}}}
+  { 48 {
\psi_{2}} {\psi_{4}}} +  { {\bar{\e}_{2}} { {\psi_{2}}^{2}}} }\\
 &\qquad {+ { 64 {\bar{\f}
_{2}} {\psi_{4}}} { -  {\bar{\f}_{2}} { {\psi_{2}}^{2}}}
+  { {\bar{\f}_{4}} {
\psi_{2}}}}
\end{aligned}
\end{split}\displaybreak[0]\tag{\ref{eq51}b}\\[3ex]
\begin{split}
\begin{aligned}
0 &= { { 3 {\psi_{4}}}+ {\bar{\f}_{4}}}  \\
0 &= { {- 5 {\psi_{6}}} + {\bar{\f}_{6}}
+  { {\bar{\e}_{2}} {\psi_{4}}}}  \\
0 &= { { 3 {\psi_{10}}}+ {\bar{\f}_{10}}
+  { {\bar{\f}_{4}} {\psi_{6}}}}
\end{aligned}
\end{split}\tag{\ref{eq51}c}
\end{align*}

These equations are  in a kind of triangular form, and  it is clear
they can be solved for the unknown coefficients
 $\psi_i$ in terms of the coefficients $\bar{\e}_i$ and $\bar{\f}_i$.
Making the corresponding change of generator produces the desired
generators for $L$.
\end{pf}

\begin{remark}
For the purposes of practical computation of these coefficients $\psi_i$
(which we carry out in section 9) it
 is important to observe that these equations  only depend
on a subset of the $\bar{\e}_i$ and $\bar{\f}_i$. Explicitly,  in case
$E_6$, they depend on
$\{ \bar{\f}_1, \bar{\f}_2, \bar{\e}_2, \bar{\f}'_3, \bar{\f}''_3,
\bar{\f}_4, \bar{\f}_6 \}$;
in case $E_7$, on
$\{ \bar{\e}_2, \bar{\f}_2, \bar{\f}_4, \bar{\f}_6 \}$;
and in case $E_8$, on
$\{ \bar{\e}_2, \bar{\f}_4, \bar{\f}_6, \bar{\f}_{10} \}$.
\end{remark}

To complete the proof of theorem 1, we first observe that the equation
$W=0$ defines ${\cal C} \subset \bar{\cal P}$ in the generators we
are considering.  What we still need to check is that the defining
polynomial $\Phi_{E_n}$ which we have found is invariant under the Weyl
group, and so can be used to define ${\cal X}$ as well as $\bar{\cal P}$.
This essentially follows from an argument of
Pinkham \cite[pp.\ 196-198]{[P-RDP]}, in the following way.

Let $\widetilde{\cal X} \to \Def (X_0)$ be the
$\C ^*$-semi-universal deformation
given by the hypersurface in preferred versal form.
(The defining polynomial of that hypersurface is the polynomial
$\Phi_{E_n}$, with $W$ set equal to 1.)
Since ${\cal P}$ is a deformation of $X_0$ with $\C ^*$-action,
there is a natural $\C ^*$-equivariant map
$\rho\colon V \to \Def (X_0)$ such that
${\cal P} \isom  \widetilde{\cal X} \times_{\rho} V$.
(The effect of proposition~\ref{prop51}
is to compute this map explicitly.)  The isomorphism
${\cal P} \isom  \widetilde{\cal X} \times_{\rho} V$
is in fact $\W$-equivariant, where $\W$ acts by Cremona transformations
on the left, and as the Weyl group on the right.  Thus, $\rho$ factors
through a map $V/\W \to \Def (X_0)$.  But the weights of the
$\C ^*$-actions are the same for these two spaces, so the map
is an isomorphism.  It follows that ${\cal P}/\W$ gives a semi-universal
deformation ${\cal X}$, and that $\Phi_{E_n}$ is $\W$-invariant.

This argument also implies that the $\W$-invariant functions $\e_i$
serve as generators of the ring $\C [V]^{\w}$ of  Weyl group
invariants.  Those generators (multiplied by appropriate integers)
are shown explicitly in cases $E_6$ and
$E_7$ in Appendices 1 and 2, respectively.

\section{Singularities of the  simultaneous partial  resolutions.}

In this section, we show that
the semi-universal deformations constructed in theorem 1 are essentially
unique (up to the
$\C ^*$-action).  We then use this uniqueness result to study
the singularities of simultaneous partial resolutions.

\begin{lemma} \label{lem61}
Let $V$ be the root space of an irreducible root system $R$.  Give $V$
 the linear $\C ^*$-action $x \mapsto \l  \cdot x$ for all $x \in V$,
$\l  \in \C ^*$.
Let $\gamma\colon V \to V$
be a $\C ^*$-equivariant map such that $\gamma (v^{\perp}) \subset
 v^{\perp}$
for all $v \in R$.  Then $\gamma = \lambda \cdot 1_V$ for some constant
$\lambda \in \C $.
\end{lemma}

\begin{pf}
We first note that $\gamma$ must be a linear map.  For if we expand
$\gamma$ in a Taylor series at $0$ and compare coefficients in the
equation $\gamma(\l x) = \l \gamma (x)$, we see that all higher
order terms must vanish.

Let $\{ v_i \}$ be a basis of $V$ consisting of roots $v_i \in R$.
Define vectors $ v_j^{\vee} \in V$ by the property
$\ip{v_j^{\vee}}{v_i} = \delta_{ij}$.  (The vectors $ v_j^{\vee}$
correspond to the elements of the dual basis $v_j^*$ of $V^*$
under the isomorphism $V \isom V^*$ determined by the inner product.)

Now $\C  \cdot v_j^{\vee} = \bigcap_{i \ne j} v_i^{\perp}$.  Thus
$v_j^{\vee} \in \bigcap_{i \ne j} v_i^{\perp}$, which implies that
$\gamma(v_j^{\vee}) \in \bigcap_{i \ne j} v_i^{\perp}$, and hence
$\gamma(v_j^{\vee}) = \lambda_j v_j^{\vee}$ for some constant $\lambda_j$.

If the vertex $v_j$ is connected to the vertex $v_k$ in the Dynkin
diagram, then $v_j + v_k$ must be a root as well.  Let
$w_{jk} = v_j^{\vee} - v_k^{\vee}$.  Then $w_{jk}$ lies in the space
$(v_j + v_k)^{\perp} \cap \bigcap_{i \ne j, k} v_i^{\perp}$,
and in fact it spans that space.  Thus, since that space is preserved
by $\gamma$, we must have $\gamma(w_{jk}) = \lambda_{jk} w_{jk}$ for
some constant $\lambda_{jk}$.  But then
\[\lambda_{jk} (v_j^{\vee} - v_k^{\vee}) = \gamma ( w_{jk})
= \gamma(v_j^{\vee}) - \gamma(v_k^{\vee})
= \lambda_j v_j^{\vee} - \lambda_k v_k^{\vee}.\]
It follows that $\lambda_j = \lambda_k = \lambda_{jk}$.

Thus, since  $R$ has a connected Dynkin diagram, all the $\lambda_j$
must be equal to the same constant $\lambda$.
\end{pf}

Let  $Z_0 \to X_0$ be the minimal
resolution of a rational double point of type $S$, and let ${\cal X} \to V/\W$
and
${\cal Z} \to {\cal X} \times_{V/\w} V$ be the deformation and
simultaneous resolution constructed in theorem 1 (which we fix
once and for all).
We call  ${\cal X} \to V/\W$ the {\em standard  deformation\/}
and ${\cal Z} \to {\cal X} \times_{V/\w} V$ the {\em standard
simultaneous resolution\/} of type $S$.
The coefficients
of the defining polynomial
$\Phi_S$ of ${\cal X}$ are specific functions
on the deformation space:
we call them the {\em standard coordinate functions\/} on $\Def (S)=V/\W$.
We call the standard coordinate function
 of highest weight the {\em ``constant term"},
since it occurs as the constant term in the defining polynomial
 of the hypersurface.
Notice that all of these definitions depend on identifying the type
as $S$.

\begin{theorem}  
Let ${\cal X'} \to {\cal D}$ be a nontrivial
$\C ^*$-equivariant deformation of $X_0$, and let
${\cal R}$
be a vector space with a linear $\C ^*$-action of pure weight one.
  Suppose that
${\cal R} \to {\cal D}$ is a $\C ^*$-equivariant map, and
that there is a simultaneous resolution
${\cal Z'} \to {\cal X'} \times_{{\cal D}} {\cal R}$
inducing $Z_0 \to X_0$.
For each root $v \in R$, let ${\cal H}_v \subset {\cal R}$
be the locus in ${\cal R}$ parametrizing deformations of $Z_0$
to which $C_v$
lifts.

Suppose that $\alpha\colon  {\cal R} \to V$ is a $\C ^*$-equivariant
surjective map such that $\alpha({\cal H}_v) \subset v^{\perp}$
for all $v \in R$.  Then $\alpha$ descends to a map
 $\bar{\alpha}\colon  {\cal D} \to V/\W$, and
 ${\cal X'}$ is isomorphic to ${\cal X} \times_{\bar{\alpha}} {\cal D}$,
the pullback of the standard deformation via $\bar{\alpha}$.
\end{theorem}

\begin{pf}
Let ${\cal X} \to V/\W$ be the standard deformation, which is
$\C ^*$-semi-universal by construction.
  By a theorem of Pinkham \cite{C*}, since
${\cal X'} \to {\cal D}$ is a
$\C ^*$-equivariant deformation of $X_0$,
there is a $\C ^*$-equivariant
map $\bar{\beta}\colon  {\cal D} \to V/\W$ such that
${\cal X'}$ is isomorphic to the pullback
${\cal X} \times_{\bar{\beta}} {\cal D}$.
Since ${\cal X'}$ admits a simultaneous resolution after base
change to
${\cal R}$, there is a map $\beta\colon {\cal R} \to V$ which induces
$\bar{\beta}$.
Since the locus in $V$ parametrizing deformations of $Z_0$ to which
$C_v$ lifts is $v^{\perp}$ by construction, the functoriality of
that property of deformations implies that
$\beta({\cal H}_v) \subset v^{\perp}$
for all $v \in R$.

It will suffice to show that $\beta = \lambda \cdot \alpha$ for some
$\lambda \in \C $.  For if that is the case, then $\lambda$ cannot
be zero since ${\cal X'} \to {\cal D}$ is a nontrivial deformation.
The desired map $\bar{\alpha}$ will then
simply be $\lambda^{-1} \cdot \bar{\beta}$,
and by using the action of $\lambda \in \C ^*$ on ${\cal X}$,
${\cal X} \times_{\bar{\alpha}} {\cal D}$ will be isomorphic to
${\cal X} \times_{\bar{\beta}} {\cal D}$, and therefore to ${\cal X'}$.

To show that $\beta = \lambda \cdot \alpha$, let $\{v_1,\ldots,v_r\}$
be a root basis of $V$, and define
\[{\cal K} = \bigcap_{i=1}^r {\cal H}_{v_i}.\]
Since $\bigcap_{i=1}^r {v_i}^{\perp} = (0)$, we
have ${\cal K} \subset \Ker  \alpha$ and
${\cal K} \subset \Ker  \beta$.  There are thus
induced maps $\widetilde{\alpha}\colon  {\cal R}/{\cal K} \to V$ and
$\widetilde{\beta}\colon  {\cal R}/{\cal K} \to V$, and it will suffice to show
that $\widetilde{\beta} = \lambda \cdot \widetilde{\alpha}$.

Since $\alpha$ is
surjective, if we define $m \coloneq \dim  {\cal R}$ then we have
\[m - r \le \dim  {\cal K} \le \dim  \Ker  \alpha = m - r,\]
which implies that $\Ker  \alpha = {\cal K}$.  Thus,
$\widetilde{\alpha}$ is an isomorphism.  Let
$\gamma \coloneq \widetilde{\beta} \circ \widetilde{\alpha}^{-1} \in \Aut  V$.
$\gamma$ is clearly $\C ^*$-equivariant.  And for
 each $v \in R$ we have
\[\gamma(v^{\perp}) = \widetilde{\beta}({\cal H}_v \bmod{\cal K})
\subset v^{\perp}.\]
Therefore, $\gamma$ satisfies the hypotheses of the lemma, so
$\gamma = \lambda \cdot 1_V$ for some $\lambda$.  It follows that
$\widetilde{\beta} = \lambda \cdot \widetilde{\alpha}$, as required.
\end{pf}

\begin{corollary}
If ${\cal X'} \to V/\W$ is  any $\C ^*$-semi-universal deformation
parametrized by $V/\W$
 such that the spaces ${\cal H}_v$ coincide with the
linear subspaces $v^{\perp} \subset V$, then ${\cal X'}$ is isomorphic
to the standard deformation ${\cal X}$.  In particular, ${\cal X'}$
can be embedded
as a hypersurface in $\C ^3 \times V/\W$ with
defining polynomial given in table~\ref{table2B},
where the coefficients in the defining polynomial
are the standard coordinate functions on $V/\W$.
\end{corollary}

We define {\em standard coordinate functions\/} on $V/\W$ in two other
cases.  The formulas for the standard coordinate functions on
$V_{D_n}/\W_{D_n}$
given in equation (\ref{eqD})
 make sense even
when $n=2$, and generate $\C [V_{D_2}]^{\w_{D_2}}$; we call them
the {\em standard coordinate functions\/} on $V_{D_2}/\W_{D_2}$.
And in the case of $E_3$, we define the {\em standard coordinate functions\/}
on $V_{E_3}/\W_{E_3}$
to be the generators $\e_2^{(1)}$, $\e_2^{(2)}$, and $\e_3$ of
$\C [V_{E_3}]^{\w_{E_3}}$ which are given in equation~(\ref{eqE3inv}).

In these cases, the functions do not directly have an interpretation
as coefficients of a  semi-universal deformation.  We relate the
standard coordinate functions on $V_{D_2}/\W_{D_2}$ to the deformation
theory in the next lemma.
A similar calculation could be done for $E_3$, but we omit it since
we do not need the result.

\begin{lemma} \label{lem72}
Let $f_{D_2}(U;t)$ be the distinguished polynomial
for $D_2$, let $\c_2 = f_{D_2}(0;t)$,
and let $g_{D_2}(Z;t) = Z^2 + \d_2 Z + (\c_2)^2$ be the
second distinguished polynomial,
so that $\c_2$ and $\d_2$ are the standard coordinate functions on
$V_{D_2}/\W_{D_2}$.
Then the ``constant terms" for
the two root systems of type $A_1$ which are the irreducible constituents
of $R_{D_2}$
 are $-\frac14(\d_2-2\c_2)$ and $-\frac14(\d_2+2\c_2)$
respectively.
\end{lemma}

\begin{pf}
The irreducible constituents of $R_{D_2}$
are spanned by the root vectors $v_1$
and $v_2$, respectively.
Let $t'_1, t'_2$ be the distinguished functionals for the $A_1$ spanned
by $v_1$, and $t''_1, t''_2$ be those for $v_2$.
The ``constant terms" are then $t'_1 t'_2$ and $t''_1 t''_2$, respectively.
Moreover, the definition of distinguished functionals implies
\begin{alignat*}{2}
t'_1 &= v_1^*,  & \qquad t'_2 &= -v_1^* \\
t''_1 &= v_2^*, & \qquad  t''_2 &= -v_2^*.
\end{alignat*}

Now according to table~\ref{table12},
\begin{gather*}
v_1^* =\frac12s_1 - t_2 = \frac12(t_1-t_2)\\
v_2^* =\frac12s_1 = \frac12(t_1+t_2),
\end{gather*}
where $t_1$ and $t_2$ are the distinguished functionals for $D_2$.
Thus,
\begin{gather*}
t'_1 t'_2 = -(v_1^*)^2 = -\frac14(t_1-t_2)^2 = -\frac14(s_1^2-4s_2)\\
t''_1 t''_2 = -(v_2^*)^2 = -\frac14(t_1+t_2)^2 = -\frac14s_1^2.
\end{gather*}
If we
write $f_{D_2}(U;t)  = U^2 + s_1 U + s_2$,  then
$\c_2 = s_2$.  Furthermore,  $g_{D_2}(-U^2;t) = U^4 + (2s_2-s_1^2)U^2 + s_2^2,$
which implies that $\d_2=s_1^2-2s_2$.
It follows that
$t'_1 t'_2 = - \frac14(\d_2-2\c_2)$ and
$t''_1 t''_2 = - \frac14(\d_2+2\c_2)$.
\end{pf}

We now turn to the study of simultaneous partial resolutions.
Let $Y_0 \to X_0$ be a partial resolution of a singularity of type $S$,
and let $Z_0 \to Y_0$
be the minimal resolution.  The $\C ^*$-action on $X_0$ lifts
to a $\C ^*$-action on $Y_0$.
As in section 1, there is an associated
partial resolution graph $\Gamma_0 \subset \Gamma$.  We write
$\Gamma - \Gamma_0 = \bigcup \Gamma_i$ as a union of its connected
components, and let $\W_0 = \prod \W_i$ be the subgroup of $\W$ generated
by reflections corresponding to vertices of $\Gamma - \Gamma_0$.
(Such vertices are illustrated with closed circles ($\bullet $)
in figures~\ref{figure1} and \ref{figure2}.)
The components $\Gamma_i$ correspond to the singular
points $Q_i$ of $Y_0$.

Let ${\cal Z} \to {\cal X}$ be the standard simultaneous resolution of
type $S$.  By the techniques of
\cite{[Wahl]}, deformations of $Z_0$ can be partially
blown down to deformations of $Y_0$, and $\C ^*$-actions are preserved
when this is done.  Doing this universally gives a family
$\widehat{\cal Z} \to V$.  Moreover, as Pinkham argues in \cite{[P]}, the map
$\widehat{\cal Z} \to V$ is $\W_i$-equivariant for each $i$, since
it provides a model for simultaneous resolution of $Q_i$.
(The argument uses a result of Burns and Wahl \cite{[BW]}, as
refined by Pinkham \cite{[P-Leop]}.)
Thus, if we let ${\cal Y} = \widehat{\cal Z}/\W_0$, then
${\cal Y} \to V/\W_0$  is a
$\C ^*$-semi-universal deformation of $Y_0$.
We call this the {\em standard simultaneous partial resolution\/} of type
$(S,\Gamma_0)$.

Let $Y_0^{(i)}$ resp.\ ${\cal Y}^{(i)}$ be the union of all $\C ^*$-orbits
in $Y_0$ resp.\ ${\cal Y}$ whose closure contains $Q_i$.  Then
$Y_0^{(i)}$ is a $\C ^*$-neighborhood of $Q_i \in Y_0$ and
${\cal Y}^{(i)} \to V/\W_0$ is a deformation of $Y_0^{(i)}$.

Now $Q_i \in Y_0^{(i)}$ is itself a rational double point, whose associated
root system $R_i$ is isomorphic to the subsystem of $R$ spanned by the
roots from $\Gamma_i$.  If we specify the type of this rational double point
as one of $A_{n-1}$, $D_n$, or $E_n$, and identify the complex root space
of the  root system $R_i$ with the subspace $V_i \subset V$ spanned by
the roots from $\Gamma_i$, then there is an associated standard deformation
${\cal X}_i \to V_i/\W_i$ of type $R_i$.

\begin{theorem}  
Let $\pr _i\colon V/\W_0 \to V_i/\W_i$ be the map induced by the orthogonal
projection $V \to V_i$.  Then ${\cal Y}^{(i)}$ is isomorphic to
${\cal X}_{i} \times_{\pr _i} V/\W_0$.  In other words,
there is a neighborhood of $Q_i \in {\cal Y}$ which can be embedded
in $\C ^3 \times V/\W_0$ in such a way that the coefficients of
the defining polynomial are pullbacks via $\pr _i$ of the
standard coordinate functions on $V_i/\W_i$.
\end{theorem}

\begin{pf}
Let ${\cal Z}\to {\cal X}\times_\rho V$ be the standard
simultaneous resolution.  Denoting the quotient map $V\to V/\W_0$ by
$\sigma$, there is a simultaneous resolution ${\cal Z}^{(i)}\to {\cal
Y}^{(i)}\times_{\sigma} V$, where ${\cal Z}^{(i)}$ is the union of all ${\bf
C}^*$-orbits in ${\cal Z}$ whose closure intersects some exceptional curve
lying
over $Q_i$.  For each root $v\in
R_i$, the locus in $V$ parametrizing deformations of $Z_0$ to which $C_v$
lifts is precisely the orthogonal complement of $v$ in $V$.  This space is
clearly mapped into the orthogonal complement of $v$ in $V_i$ by the orthogonal
projection $V\to V_i$.  The theorem follows by applying theorem 2 to this
situation, using the root system $R_i$ in place of $R$, and the orthogonal
projection in place of $\alpha$.
\end{pf}

Consider now the case of an irreducible small resolution, so
that $\Gamma_0$ consists of a single vertex $v$.  If
we write $\Gamma - \{v\} = \bigcup \Gamma_i$
as a union of its irreducible components, then
theorem 3 provides us with a natural isomorphism
$\PRes (S,v)  \isom  V/\W_0$, compatible with the
 maps
$\pr _i\colon V/\W_0 \to V_i/\W_i$
which are induced by the orthogonal projections.
In fact, $\W_0$ fixes the
hyperplane $\Ker (v^*) \subset V$, and the projections $\pr _i$
induce an isomorphism $\Ker (v^*)/\W_0 \isom \bigoplus V_i/\W_i$.
Here, $v^*\in V^*$ is the element dual to $v$ in the basis dual to the root
basis.

We define standard coordinate functions on the partial resolution
spaces $\PRes (S,v)$
in the following way.
The definition depends on choosing a decomposition
$\Gamma - \{ v \} = \bigcup \Gamma^{(j)}$,
where this time each $\Gamma^{(j)}$ is a union of connected components
of $\Gamma - \{ v \}$.  We assume that each corresponding root system
$R^{(j)}$ is either irreducible or of type $D_2$ or $E_3$; the definition
also depends on identifying the type of each $R^{(j)}$
as one of $A_{k-1}$, $D_k$, or $E_k$.
We define the {\em standard coordinate functions on
$\PRes (S,v) \isom V/\W_0$\/}
to be the linear functional $v^*$ together with the pullbacks of
the standard coordinate functions on the spaces $V^{(j)}/\W^{(j)}$ via the
 mappings $V/\W_0 \to V^{(j)}/\W^{(j)}$ which are induced by
orthogonal projection.
These will be the coefficients of the defining polynomials of various
open sets on the partial resolution, except in the cases
where one of the constituents is a reducible root system of type
$D_2$ or $E_3$.

\section{Relations among distinguished polynomials.}

Let $R$ be a root system of type $S=A_{n-1}$, $D_n$, or $E_n$.  The
complex root space $V_R$
 has a root basis which can be written in the form
$v_{\alpha}, v_{\alpha + 1}, \ldots, v_{\beta}$, where $\alpha,
\alpha + 1, \ldots, \beta$ is a sequence of consecutive integers.
There is a natural dual basis $v_{\alpha}^*, v_{\alpha + 1}^*,
 \ldots, v_{\beta}^*$ of the dual space $V^*$.

The basis element $v_k$ can be regarded as a vertex of the Dynkin
diagram $\Gamma_R$.
We fix  $v_k$, and regard
the complement $\Gamma_R - \{v_k\}$ as forming two root systems, either of
which may
be empty or reducible:  $R'$ is the part spanned by vertices
to the left of $v_k$ in $\Gamma_R - \{v_k\}$, and $R''$ is the part
spanned by vertices to the right.
(We use the orientation of the Dynkin diagrams displayed in section 2.)
This is admittedly ambiguous in a few cases, so to make everything
completely clear, we specify that when $(S,v_k) = (D_n,v_n)$ or $(E_n,v_0)$,
we take $R'$ to be spanned by $\Gamma_R - \{v_k\}$, and $R''$ to be empty.

We fix types $S'$ and $S''$ of the root systems $R'$ and $R''$,
as indicated in table~\ref{table4}.  (We have omitted the case
$(S,v_k)=(D_n,v_{n-1})$, since we will not need it later.)
Having identified the type $S'$, there
are distinguished functionals  $t_1',\ldots,t_{n'}'$, and
a distinguished polynomial  $f_{S'}(U;t')$.  The coefficients
of this polynomial are denoted by $s_i'$, and the standard coordinate
functions on $V'/\W'$ are denoted by $\a_i'$, $\c_i'$, $\d_i'$, or $\e_i'$, as
appropriate.  We use analogous notation for $R''$, replacing `` $'$ ''
by `` $''$ '' throughout.

\begin{table}[t]
\begin{center}
$\begin{array}{|c|c|c|c|l|}  \hline
 & & & & \\
S & k & S' & S'' & \multicolumn{1}{c|}{\widetilde{v}_k} \\
 & & & & \\ \hline
 & & & & \\
A_{n-1} & \text{any} & A_{k-1} & A_{n-k-1} & v_k +
    \sum_{i=1}^{k-1} \frac{i}{k} v'_i + \sum_{i=1}^{n-k-1}
    \frac{n-k-i}{n-k} v''_i \\
 & & & & \\
D_n & \le n-2 & A_{k-1} & D_{n-k} & v_k +
    \sum_{i=1}^{k-1} \frac{i}{k} v'_i + \sum_{i=1}^{n-k-2} v''_i
                   + \frac12 v''_{n-k-1} + \frac12 v''_{n-k} \\
 & & & & \\
D_n & n & A_{n-1} & & v_n + \sum_{i=1}^{n-2} \frac{2i}{n} v'_i
                 + \frac{n-2}{n} v'_{n-1}   \\
 & & & & \\
E_n & 0 & A_{n-1} & &  v_0 + \frac{n-3}{n} v'_1 + \frac{2n-6}{n} v'_2 +
               \sum_{i=3}^{n-1} \frac{3n-3i}{n} v'_i \\
 & & & & \\
E_n & 1 & D_{n-1} & & v_1 + \sum_{i=1}^{n-3} \frac{i}{2} v'_i +
               \frac{n-1}{4} v'_{n-2} + \frac{n-3}{4} v'_{n-1} \\
 & & & & \\
E_n & 2 & A_1 & A_{n-2} & v_2 + \frac12 v'_1 + \frac{n-3}{n-1} v''_1 +
               \sum_{i=2}^{n-2} \frac{2n-2i-2}{n-1} v''_i \\
 & & & & \\
E_n & \ge 3 & E_k & A_{n-k-1} & v_k
   + \frac{3 v'_0 + 2 v'_1  + 4 v'_2 }{9-k}
   + \sum_{i=3}^{k-1} \frac{9-i}{9-k} v'_i
   + \sum_{i=1}^{n-k-1} \frac{n-k-i}{n-k} v''_i \\
 & & & & \\ \hline
\end{array}$
\end{center}

\medskip

\caption{}
\label{table4}
\end{table}

By composing with the orthogonal projection maps $V_R \to V_{R'}$ and
$V_R \to V_{R''}$, we can regard the distinguished functionals on
$V_{R'}$ and $V_{R''}$ as linear functionals on $V_R$.  We denote
them again by $t'_i$ and $t''_i$, suppressing
 mention of the orthogonal projection maps for simplicity
of notation.

\begin{proposition} \label{prop71}
Let $v_k$ be a vertex of the Dynkin diagram $\Gamma_R$, and let
 $\widetilde{v}_k \in V_R$   be the
 vector specified in table~\ref{table4}.  Then there is an
orthogonal direct sum decomposition
\[V_R = (\C \cdot\widetilde{v}_k)  \oplus V_{R'} \oplus V_{R''},\]
and the dual space $V_R^*$ can be generated by the distinguished functionals
of $R'$ and $R''$ together with the functional $\m_1\coloneq v_k^*$.

The distinguished functionals of $R$ are therefore linear
combinations of these generators, and those linear combinations can
be expressed by means of a relation among distinguished polynomials.
The relation involves $f_S(U;t)$, $f_{S'}(U;t')$,  $f_{S''}(U;t'')$,
and $\m_1$,
with some linear changes of variable and possible extra linear factors,
and is given explicitly
 in table~\ref{tableAA}.
\end{proposition}

\begin{table}[b]
\begin{center}
\begin{tabular}{|c|c|rcl|} \hline
 & & & & \\
 $S$ & $k$ & & & \\
 & & & & \\ \hline
 & & & & \\
 $A_{n-1}$ &  any    & $f_{A_{n-1}}(U;t)$ & = &
   $f_{A_{k-1}}(U + \frac{1}{k} \m_1;t') \cdot
    f_{A_{n-k-1}}(U - \frac{1}{n-k} \m_1;t'')$  \\
 & & & & \\
 $D_n$ & $ \le n-2$ & $f_{D_n}(U;t)$ & = &
   $f_{A_{k-1}}(U + \frac{1}{k} \m_1;t') \cdot
    f_{D_{n-k}}(U;t'')$  \\
 & & & & \\
 $D_n$ & $n$       & $f_{D_n}(U;t)$ & = &
   $f_{A_{n-1}}(U + \frac{2}{k} \m_1;t')$  \\
 & & & & \\
 $E_n$ & $0$       & $f_{E_n}(U;t)$ & = &
   $f_{A_{n-1}}(U - \frac{9-n}{3n} \m_1;t')$  \\
 & & & & \\
 $E_n$ & $1$       & $(-1)^n \cdot f_{E_n}(U;t)$ & = &
   $(- U + \frac13 \rho_1 - \frac{9-n}{6} \m_1) \cdot
    f_{D_{n-1}}(- U - \frac16 \rho_1 + \frac{9-n}{12} \m_1;t')$  \\
 & & & & \\
 $E_n$ & $2$       & $f_{E_n}(U;t)$ & = &
   $f_{A_{1}}(U + \frac23 \sigma_1 + \frac{9-n}{6n-6} \m_1;t') \cdot
    \frac{f_{A_{n-2}}(U - \frac13 \sigma_1 - \frac{9-n}{3n-3} \m_1;t'')}
    {(U - \frac43 \sigma_1 - \frac{9-n}{3n-3} \m_1)}$  \\
 & & & & \\
 $E_n$ & $ \ge 3$   & $f_{E_n}(U;t)$ & = &
   $f_{E_k}(U;t') \cdot f_{A_{n-k-1}}(U
         - \frac{1}{9-k} \tau_1 - \frac{9-n}{(9-k)(n-k)} \m_1;t'')$  \\
 & & & & \\ \hline
\end{tabular}

\bigskip

{\it Notation:}

\medskip

\begin{tabular}{l}
 $\m_1$ denotes the coordinate function $v_k^*$ \\
 $\rho_1$ denotes the coefficient of $U^{n-2}$ in $f_{D_{n-1}}(U;t')$ \\
 $\sigma_1$ denotes a root of $f_{A_{n-2}}(U;t'')$ \\
 $\tau_1$ denotes the coefficient of $U^{k-1}$ in $f_{E_k}(U;t')$
\end{tabular}
\end{center}

\bigskip

\caption{}
\label{tableAA}
\end{table}

\begin{pf}
We identify the root bases of $V_{R'}$ and $V_{R''}$
with subsets of the root basis of $V_{R}$ as follows.  If $S \ne E_n$
or $k \ne 1, 2$
 we identify $v'_i$ with $v_i$ for all
basis vectors of $V_{R'}$, and $v''_i$ with $v_{k+i}$ for all basis vectors
of $V_{R''}$.
If $S=E_n$
and $k=1$ or $2$, we use the identifications indicated in
figure~\ref{figure2}.
With this notation established, it is easy to
check that the vector $\widetilde{v}_k$ defined in table~\ref{table4}
is orthogonal to both $V_{R'}$ and $V_{R''}$.
Moreover, since $\Gamma_{R'}$ is disjoint from $\Gamma_{R''}$, the
spaces $V_{R'}$ and $V_{R''}$ are themselves
 mutually perpendicular.  The claimed orthogonal direct
sum decomposition follows.

The orthogonal direct sum decomposition can be regarded as a change
of basis from $\{v_i\}$ to
$\{v'_i\} \union \{v''_i\} \union \{\widetilde{v}_k\}$
in the space $V$.  If we use maps $\s'$ and $\s''$ to describe
the identification of root bases, then this change of basis
can be written
\begin{align*}
v_i' & =  v_{\s'(i)} \\
v_i'' & =  v_{\s''(i)} \\
\widetilde{v}_k & =  {\textstyle\sum a'_i v_i' + \sum a''_i v_i''},
\end{align*}
where the coefficients $a'_i$ and $a''_i$ are found in table~\ref{table4}.
(In all but two cases, $\s'(i)=i$ and $\s''(i)=k+i$.)
The corresponding change of dual basis
 takes the form
\begin{align*}
v_{\s'(i)}^* & = a'_i \ \widetilde{v}_k^* +  {v_i'}^* \\
v_{\s''(i)}^* & =  a''_i \ \widetilde{v}_k^* + {v_i''}^* \\
v_{k}^* & =   \widetilde{v}_k^*.
\end{align*}
It follows that $\m_1 = v_{k}^*  =   \widetilde{v}_k^*$ can be used
along with the distinguished functionals on $V_{R'}$ and $V_{R''}$
to generate $V_R^*$.

To finish the proof, we must carry out the calculation which leads
to table~\ref{tableAA}.  We will do this in a few cases, and leave
the remaining ones to the reader.

We first treat an easy case:  the case $S = A_{n-1}$.
Using the fourth column of
table~\ref{table12} applied to $R'$ and $R''$,
we can write the change of basis as
\begin{alignat*}5
v_i^* & = &  \tfrac{i}{k} &\widetilde{v}_k^* + {v_i'}^*
&& = & \tfrac{i}{k} &\m_1 + t_{1}' + \cdots + t_i'
, \quad  && 1 \le i \le k-1 \\
v_k^* & = & &\widetilde{v}_k^*
&& = & &\m_1
  \\
v_{k+i}^* & = & \tfrac{n-k-i}{n-k} &\widetilde{v}_k^* +  {v_i''}^*
&& = & \tfrac{n-k-i}{n-k} &\m_1 + t_{1}'' + \cdots + t_{i}''
, \quad  && 1 \le i \le n-k-1
\end{alignat*}
Then using the third column of
 table~\ref{table12} applied to $R$, we get
\begin{alignat*}2
t_i & =  \tfrac{1}{k} \m_1 + t_i'
   , \quad &&  1 \le i \le k
\\
t_{k+i} & =  - \tfrac{1}{n-k} \m_1 + t_i''
    , \quad &&   k+1 \le k+i \le  n.
\end{alignat*}
It follows that
\[f_{A_{n-1}}(U;t) = f_{A_{k-1}}(U + \tfrac{1}{k} \m_1;t')
     \cdot f_{A_{n-k-1}}(U - \tfrac{1}{n-k} \m_1;t'').\]

\begin{figure}[t]
\begin{picture}(2.6,1)(1.9,.5)
\thicklines
\put(1.9,1){\circle{.075}}
\put(1.9375,1){\line(1,0){.4625}}
\put(2.4,1){\circle*{.075}}
\put(2.4,1){\line(1,0){.5}}
\put(2.9,1){\circle*{.075}}
\put(2.9,1){\line(0,-1){.5}}
\put(2.9,.5){\circle*{.075}}
\put(2.9,1){\line(1,0){.5}}
\put(3.4,1){\circle*{.075}}
\put(3.4,1){\line(1,0){.25}}
\put(3.75,1){\circle*{.02}}
\put(3.85,1){\circle*{.02}}
\put(3.95,1){\circle*{.02}}
\put(4.05,1){\line(1,0){.25}}
\put(4.3,1){\circle*{.075}}
\put(2.275,1.15){\makebox(.25,.25){$v'_{n-2}$}}
\put(2.775,1.15){\makebox(.25,.25){$v'_{n-3}$}}
\put(3.025,.375){\makebox(.25,.25){$v'_{n-1}$}}
\put(3.275,1.15){\makebox(.25,.25){$v'_{n-4}$}}
\put(4.175,1.15){\makebox(.25,.25){$v'_{1}$}}
\end{picture}
\hspace*{\fill}
\begin{picture}(2.6,1)(1.9,.5)
\thicklines
\put(1.9,1){\circle*{.075}}
\put(1.9,1){\line(1,0){.4625}}
\put(2.4,1){\circle{.075}}
\put(2.4375,1){\line(1,0){.4625}}
\put(2.9,1){\circle*{.075}}
\put(2.9,1){\line(0,-1){.5}}
\put(2.9,.5){\circle*{.075}}
\put(2.9,1){\line(1,0){.5}}
\put(3.4,1){\circle*{.075}}
\put(3.4,1){\line(1,0){.25}}
\put(3.75,1){\circle*{.02}}
\put(3.85,1){\circle*{.02}}
\put(3.95,1){\circle*{.02}}
\put(4.05,1){\line(1,0){.25}}
\put(4.3,1){\circle*{.075}}
\put(1.775,1.15){\makebox(.25,.25){$v'_1$}}
\put(2.775,1.15){\makebox(.25,.25){$v''_2$}}
\put(3.025,.375){\makebox(.25,.25){$v''_1$}}
\put(3.275,1.15){\makebox(.25,.25){$v''_3$}}
\put(4.175,1.15){\makebox(.25,.25){$v''_{n-2}$}}
\end{picture}

\caption{}
\label{figure2}
\end{figure}

We next treat the case of $R = E_n$, with $k=1$, which is displayed in the
left half of figure~\ref{figure2}.  In this case,
\begin{alignat*}5
v_0^* & = & \tfrac{n-3}{4} &\widetilde{v}_1^* + {v_{n-1}'}^*
&& = & \tfrac{n-3}{4} &\m_1 + \tfrac12 s_1'
 \\
v_1^* & = & &\widetilde{v}_1^*
&& = & &\m_1
 \\
v_2^* & = & \tfrac{n-1}{4} &\widetilde{v}_1^* + {v_{n-2}'}^*
&& = & \tfrac{n-1}{4} &\m_1 + \tfrac12 s_1' - t_{n-1}'
 \\
v_{n-i}^* & = & \tfrac{i}{2} &\widetilde{v}_1^* + {v_i'}^*
&& = & \tfrac{i}{2} &\m_1 + s_1' - t_{i+1}' - \cdots - t_{n-1}'
, \quad &&    1 \le i \le n-3
\end{alignat*}
which implies
\begin{alignat*}2
t_1 & =  \tfrac{9-n}{6} \m_1 - \tfrac13 s_1'
\\
t_{i+1} & =  \tfrac{n-9}{12} \m_1 + \tfrac16 s_1' - t_{n-i}'
    , \quad &&   2 \le i+1 \le n.
\end{alignat*}
It follows that
\[f_{E_n}(U;t) =
(-1)^{n-1} \cdot (U + \tfrac{9-n}{6} \m_1 - \tfrac13 s_1') \cdot
     f_{D_{n-1}}(- U - \tfrac{n-9}{12} \m_1 - \tfrac16 s_1';t'),\]
since the right-hand side is equal to
\begin{multline*}
 (-1)^{n-1} \cdot (U + \tfrac{9-n}{6} \m_1 - \tfrac13 s_1')
\cdot \prod (- U - \tfrac{n-9}{12} \m_1 - \tfrac16 s_1' + t_j') \\
\begin{aligned}
 & =  (U + \tfrac{9-n}{6} \m_1 - \tfrac13 s_1') \cdot
\prod (U + \tfrac{n-9}{12} \m_1 + \tfrac16 s_1' - t_j') \\
& =  \prod (U + t_i).
\end{aligned}
\end{multline*}

Finally, we treat the case $R = E_n$, $k = 2$, which is displayed in the
right half of figure~\ref{figure2}.  In this case,
\begin{alignat*}5
v_0^* & = & \tfrac{n-3}{n-1} &\widetilde{v}_2^* + {v_1''}^*
&& = & \tfrac{n-3}{n-1} &\m_1 + t_1''
 \\
v_1^* & = & \tfrac12 &\widetilde{v}_2^* + {v_1'}^*
&& = & \tfrac12 &\m_1 + t_1'
 \\
v_2^* & = & &\widetilde{v}_2^*
&& = & &\m_1
 \\
v_{i+1}^* & = & \tfrac{2n-2i-2}{n-1} &\widetilde{v}_2^* + {v_i''}^*
&& = & \tfrac{2n-2i-2}{n-1} &\m_1 + t_1'' + \cdots + t_i''
, \quad &&   2 \le i \le n-2
\end{alignat*}
which implies
\begin{alignat*}2
t_i & =  \tfrac{9-n}{6n-6} \m_1 - \tfrac23 t_1'' + t_i'
        ,  \quad &&   1 \le i \le 2
\\
t_{i+1} & =  \tfrac{n-9}{3n-3} \m_1 + \tfrac13 t_1'' + t_{i}''
        ,  \quad &&   2 \le {i} \le n-1.
\end{alignat*}
(Notice that the  functional
$\tfrac{n-9}{3n-3} \m_1 + \tfrac13 t_1'' + t_1''$
is ``missing" here.)
It follows that
\[f_{E_n}(U;t) =
f_{A_1}(U + \tfrac{9-n}{6n-6} \m_1 - \tfrac23 t_1'';t') \cdot
   \frac{f_{A_{n-2}}(U + \tfrac{n-9}{3n-3} \m_1 + \tfrac13 t_1'';t'')}
       {(U + \tfrac{n-9}{3n-3} \m_1 + \tfrac43 t_1'')}.\]
Since $-t''_1$ is a root of $f_{A_{n-1}}(U;t'')$, if we define
$\s_1\coloneq -t''_1$, the formula in the table follows.

All remaining cases are left to the reader.
\end{pf}

Proposition~\ref{prop71} provides a method for explicitly calculating
the map $\PRes (S,v_k) \to \Def (S)$ (which can also be written as
$V/\W_0 \to V/\W$).  What we wish to calculate explicitly is the
map on coordinate rings $\C [V]^{\w} \subset \C [V]^{\w_0}$.
In other words, we want to express the standard coordinate functions
on $V/\W$ as polynomials in the standard coordinate functions on
$V/\W_0$.

Now each standard coordinate function $\f_j$ on $V/\W$ is a function
of $s_1,\ldots,s_n$ (the coefficients of the distinguished polynomial
$f_S(U;t)$).  Using proposition~\ref{prop71},
these in turn are expressed as functions of the coefficients $s'_i$
and $s''_i$ of the distinguished polynomials of $f_{S'}(U;t')$
and $f_{S''}(U;t'')$, together
with $\m_1$ (and possibly an auxiliary variable $\rho_1$, $\sigma_1$,
or $\tau_1$ which will eliminate itself in the end).  The expression
for $\f_j$ in terms of these variables
is invariant under $\W_0$, and so can be expressed as a polynomial
in $\m_1$ together with the pullbacks of the
standard coordinate functions on
$V'/\W'$ and $V''/\W''$.  One approach to finding this polynomial expression
is the method of undetermined coefficients.

We carry this out  for the cases of $A_{n-1}$ and $D_n$,
collecting the information we require in the form of certain
congruences.  The part of the computation which we need in the $E_n$
cases will be stated in table~\ref{table-key1} in section 8, and verified
in section 10.

\begin{proposition} \label{prop72}
\quad
\begin{enumerate}
\item  If $R = A_{n-1}$, then
$\a_{n-1} \equiv
\a'_{k-1} \a''_{n-k} + \a'_k \a''_{n-k-1}\mod{\m_1}$, and
$\a_{n} \equiv \a'_k  \a''_{n-k}\mod{\m_1}$.

\item If $R = D_n$ and $k=1$, then
$\d_{2n-4} \equiv \d''_{2n-4}\mod{\m_1}$.

\item If $R = D_n$ and $k \le n-2$, let $J_4$ denote the ideal generated by all
monomials of degree 4 in the standard coordinate functions on
$\PRes (D_n,v_k)$.
Then
$\c_n \equiv \a'_k  \c''_{n-k}\mod{\m_1}$, and
$\d_{2n-2} \equiv  (\a'_k)^2 \d''_{2n-2k-2}\mod{J_4}$.

\item If $R = D_n$ and $k = n$, then
$\c_n \equiv \a'_{n}\mod{\m_1}$.

\end{enumerate}
\end{proposition}

\begin{pf}
We will prove the third statement, and leave the others
(which are easier) to the reader.
When $R = D_n$ and $k \le n-2$ we have
\begin{align*}
f_{D_{n}}(U;t) & =
   f_{A_{k-1}}(U + \tfrac{1}{k} \m_1;t') \cdot
    f_{D_{n-k}}(U;t'') \\
 & \equiv  f_{A_{k-1}}(U;t') \cdot
    f_{D_{n-k}}(U;t'')\mod{\m_1} \\
 & \equiv  \a'_k
    \c''_{n-k}\mod{(\m_1,U)}
\end{align*}
It follows that $\c_n \equiv \a'_k
    \c''_{n-k}\mod{\m_1}$.  Moreover,
if we define $\widetilde{g}_{A_{k-1}}(Z;t')$ by
\[
\widetilde{g}_{A_{k-1}}(- U^2;t')
= f_{A_{k-1}}(U+\tfrac{1}{k}\m_1;t') \cdot f_{A_{k-1}}(-U+\tfrac{1}{k}\m_1;t')
\]
then
\begin{multline*}
\widetilde{g}_{A_{k-1}}(Z;t')
 \equiv
({\a_{k-1}'}^2-2\a_{k-2}'\a_k'+
\tfrac{2}{k}\m_1\a_{k-1}'\a_{k-2}')Z  \\
+
({\a_k'}^2+\tfrac{2}{k}\m_1\a_{k-1}'\a_k')\mod{(J_4,Z^2)}
\end{multline*}
while
\[
g_{D_{n-k}}(Z;t'') \equiv
\d''_{2n-2k-2} Z +(\c''_{n-k})^2\mod{Z^2}.
\]
Hence, if we multiply these congruences and retain only terms of degree
at most three in the standard coordinate functions
on $\PRes (D_n,v_k)$, we get
\[
g_{D_n}(Z;t) \equiv (\a'_{k})^2
 \d''_{2n-2k-2} Z\mod{(J_4,Z^2)}.
\]
The  congruence for $\d_{2n-2}$ follows.
\end{pf}

In order to effectively apply the last line of table~\ref{tableAA}
when $k=4$ or $5$,
we need to
 compute the distinguished polynomials $f_{E_4}(U;t)$ and $f_{E_5}(U;t)$.

\begin{lemma} \label{lem71}
Let $t_1,\cdots,t_n$ be the distinguished functionals for $E_4$, resp.\
$E_5$, and let $\widetilde{t}_1,\cdots,\widetilde{t}_5$ be the distinguished
functionals for $A_4$, resp.\ $D_5$.  If we identify these root systems
in such a way that $\widetilde{v}_i = v_i$ for $1 \le i \le n-1$ and
$\widetilde{v}_n = v_0$ then
\[f_{A_4}(U;\widetilde{t}) = (U + \frac35 s_1 ) \cdot
         f_{E_4}(U - \frac25 s_1; t)\]
and
\[f_{D_5}(U;\widetilde{t}) = f_{E_5}(U - \frac12 s_1; t).\]
\end{lemma}

\begin{pf}
In the case of $E_4$, we calculate with table~\ref{table12}
as follows.
\begin{alignat*}6
\widetilde{t}_1 & = & &\widetilde{v}_1^* && = & &v_1^* && = &
   - \tfrac25 &s_1 + t_1 \\
\widetilde{t}_2 & = & - \widetilde{v}_1^* &+ \widetilde{v}_2^* && = &
   - v_1^* &+ v_2^* && = & - \tfrac25 &s_1 + t_2 \\
\widetilde{t}_3 & = & - \widetilde{v}_2^* &+ \widetilde{v}_3^* && = &
   - v_2^* &+ v_3^* && = & - \tfrac25 &s_1 + t_3 \\
\widetilde{t}_4 & = & - \widetilde{v}_3^* &+ \widetilde{v}_4^* && = &
   - v_3^* &+ v_0^* && = & \tfrac35 &s_1 - t_1 - t_2 - t_3 \\
\widetilde{t}_5 & = & - &\widetilde{v}_4^* && = &
   - &v_0^* && = & \tfrac35 &s_1
\end{alignat*}
Now since $\frac35 s_1 - t_1 - t_2 - t_3 = - \frac25 s_1 + t_4$
we get
\[(U+\widetilde{t}_1)\cdots(U+\widetilde{t}_5)
=(U-\frac25s_1+t_1)\cdots(U-\frac25s_1+t_4)
\cdot (U+\frac35s_1),\]
and the first equation follows.

The case of $E_5$ is similar (and easier), and will be left to the
reader.
\end{pf}

\begin{corollary}
The standard coordinate functions on $\Def (E_4)$ and $\Def (E_5)$,
expressed in terms of
the elementary symmetric functions $s_i$ of their respective distinguished
functionals $t_i$, are given by
\begin{align}
\begin{split}
\e_2 &= s_2 - \tfrac35 s_1^2\\
\e_3 &= s_3 - \tfrac15 s_2s_1 + \tfrac{2}{25}s_1^3\\
\e_4 &= s_4 + \tfrac15 s_3s_1 - \tfrac{8}{25}s_2s_1^2 + \tfrac{12}{125}s_1^4\\
\e_5 &= \tfrac35 s_1s_4 - \tfrac{6}{25}s_3s_1^2 + \tfrac{12}{125}s_2s_1^3 -
\tfrac{72}{3125}s_1^5
\end{split}\label{eqE4}
\end{align}
in the case of $E_4$, and by
\begin{align}
\begin{split}
\e_2 &=  -2s_2 + \tfrac54 s_1^2\\
\e_4 &= s_2^2 - 2s_2s_1^2 + \tfrac58 s_1^4 + s_1s_3 + 2s_4\\
\e_5 &=  -\tfrac18 s_2s_1^3 + \tfrac{1}{32}s_1^5 + \tfrac14 s_1^2s_3 -
\tfrac12s_1s_4 +s_5\\
\e_6 &=   \tfrac34 s_2^2s_1^2 - \tfrac34 s_2s_1^4 - s_2s_1
s_3 - 2s_2s_4 + \tfrac5{32} s_1^6 + \tfrac34 s_1^3s_3 + \tfrac12 s_1^2
s_4
\\ & \qquad
-
 3s_1s_5 - s_3^2\\
\e_8 &=  \tfrac3{16}s_2^2s_1^4 - \tfrac18 s_2s_1^6 - \tfrac12 s_2s_1^3s_3 +
 3s_2s_1s_5 +
\tfrac5{256} s_1^8 + \tfrac3{16} s_1^5s_3
\\  & \qquad
-
 \tfrac18 s_1^4s_4 - \tfrac12 s_1^3s_5 + \tfrac12 s_1^2s_3^2 - s_1s_4s_3 -
2s_5s_3 + s_4^2\\
\end{split}\label{eqE5}
\end{align}
in the case of $E_5$.
\end{corollary}

\begin{pf}
Write $f_{E_n}(U;t)=U^n+\sum_{i=0}^{n-1}s_iU^{n-i}$ for $n=4\hbox{ or }5$.  Use
the  formulas of the  lemma to calculate
$f_{A_4}(U;\tilde{t})$ and $f_{D_5}(U;\tilde{t})$.  In the case of $E_4$,
the
standard coordinate functions
$\e_i$ can be read off as the coefficients of $U^{5-i}$ in
$f_{A_4}$.  In the case of $E_5$, form the second distinguished polynomial
$g_{D_5}$ via $g_{D_5}(-U^2;\tilde{t})=
f_{D_5}(U;\tilde{t})f_{D_5}(-U;\tilde{t})$.  The standard coordinate
function $\e_{2i}$ can be read off as
the coefficient of $U^{5-i}$ in $g_{D_5}(U;\tilde{t})$, while $\e_5$ is simply
the
coefficient of $U^0$ in $f_{D_5}(U;\tilde{t})$.
\end{pf}

\section{Proof of the main theorem.}

In this section, we prove the main theorem, assuming the validity of
some results to be stated in table~\ref{table-key1}.  The proof will be
complete
when we verify that table in section 10.

The partial resolution graphs shown in figure~\ref{figure1}
determine a singularity
type $S_{\ell}$ for each length $\ell$ between $1$ and $6$, which
we call the {\em associated type\/} of the length.  Explicitly, this type is
\begin{center}
\begin{tabular}{|c|cccccc|} \hline
$\ell$ & 1 & 2 & 3 & 4 & 5 & 6 \\ \hline
$S_{\ell}$ & $A_1$ & $D_4$ & $E_6$ & $E_7$ & $E_8$ & $E_8$ \\ \hline
\end{tabular}.
\end{center}
Figure~\ref{figure1}
(in section 1) illustrates the partial resolution graphs
$\{ v \} \subset \Gamma_{S_{\ell}}$, where $v$ is the vertex corresponding
to the unique component in the maximal ideal cycle of length $\ell$.

Our aim is to show that for $\pi \colon Y \to X$ of length $\ell$, the
singularity type of the general hyperplane section is $S_{\ell}$.
We say that the singularity type of a rational double point is {\em
at worst\/} $S$ if its dual graph is isomorphic to
 a (proper or improper) subgraph
 of $\Gamma_S$.

\begin{lemma} \label{lem82}
\quad
\begin{enumerate}
\item
The partial resolution graphs shown in figure~\ref{figure1} are primitive.

\item
Let $\pi\colon Y \to X$ be an irreducible small resolution
of an isolated Gorenstein threefold singularity, and let $\ell$ be
the length.  Suppose that $X$ has a hyperplane section
whose singularity type is at worst $S_{\ell}$.
Then the generic hyperplane section defines the primitive partial
resolution graph corresponding to $S_{\ell}$ given in figure~\ref{figure1}.

\end{enumerate}
\end{lemma}

\begin{pf}
\quad

(1) For each $\ell$ between 1 and 6, let
$n(\ell)$
be the minimum $n$ such that there is a rational double point whose
dual graph has $n$ vertices, and at least one component has multiplicity
exactly $\ell$ in the maximal
ideal cycle.
Examining the maximal ideal cycles of the rational double points, it is
easy to see that for each  $\ell$, there is a  unique such rational
double point with $n(\ell)$ vertices, namely the one shown in
figure~\ref{figure1}.
Now for any nontrivial 1-parameter deformation of a rational double
point, the dual graph of the minimal resolution of the general fiber is
isomorphic to a proper subgraph of the dual graph of the special fiber.
It follows that each graph shown in figure~\ref{figure1}
is primitive:  any proper
subgraph
will have fewer vertices, and so cannot have a component in its
maximal ideal cycle of multiplicity $\ell$.

(2)
Fix $\ell$.  For the singularity of type $S_{\ell}$,
there is a unique component of multiplicity exactly $\ell$ in the
maximal ideal cycle.  Moreover, no proper subgraph
of $\Gamma_{S_{\ell}}$ has any component
with multiplicity exactly $\ell$ in {\em its\/} maximal ideal cycle.
Thus, since the length is $\ell$, the partial
resolution graph determined by the given hyperplane section must
be of the type shown in figure~\ref{figure1}
(which indicates the unique component
of multiplicity $\ell$).  On the other hand, since the graphs in
figure~\ref{figure1} are primitive, it follows that this is also the type of
the generic hyperplane section.
\end{pf}

Consider  now $\pi\colon Y \to X$,
 an irreducible small resolution
of an isolated Gorenstein threefold singularity $P \in X$,  and
 a hyperplane section $\{f=0\}$ which has a rational double point.
This determines a partial resolution graph $\{ v \} \subset \Gamma_S$, and
the length $\ell$ coincides with the multiplicity of $v$ in the maximal
ideal cycle.
There is a natural classifying map $\mu_f\colon \Delta \to \Def (S) = V/\W$
which allows us to recover a neighborhood of $P \in X$ as the pullback
of the standard deformation ${\cal X} \to V/\W$.
The map $\mu_f$ determines a discrete valuation
$\nu_f\colon \C [V]^{\w} \to \Z $
which is defined by
\[\nu_f(\f) = \text{order of vanishing at $0$ of $\mu^*_f(\f)$}.\]

Thanks to
lemma~\ref{lem82}, in order to prove the main theorem it suffices
 to show that $X$ has some hyperplane section whose
singularity type is at worst $S_{\ell}$.
  The following proposition shows how to use the
discrete valuation $\nu_f$ applied to the
standard coordinate functions  on $\Def (S)$ (or in the $E_7$
case, to certain simple polynomial expressions in these functions)
 to bound the
singularity type of the general hyperplane section.

\begin{proposition} \label{prop81}
Let $X$ be  a threefold with an isolated rational Gorenstein
 singular point $P$, and let $\{f=0\}$ be a hyperplane section
through $P$ with a rational double point of type $S$.  Let
$\mu_f\colon \Delta \to \Def (S) = V/\W$ be the classifying map, and let
$\nu_f\colon \C [V]^{\w} \to \Z $ be the
associated discrete valuation.
\begin{enumerate}
\item
Suppose that $S=A_{n-1}$, and let $\{ \a_i \}$ be the standard coordinate
functions on $\Def (S)$. If $\nu_f(\a_{n-1})=1$ or $\nu_f(\a_n)=2$
then the general hyperplane section of $X$ has  singularity type at worst
$A_1$.

\item
Suppose that $S=D_{n}$, and let $\{ \c_n,\ \d_i \}$ be the standard coordinate
functions on $\Def (S)$. If $\nu_f(\c_{n})=1$ or $\nu_f(\d_{2n-4})=1$
then the general hyperplane section of $X$ has  singularity type at worst
$A_1$, while if $\nu_f(\c_{n})=2$ or $\nu_f(\d_{2n-2})=3$
then the general hyperplane section of $X$ has  singularity type at worst
$D_4$.

\item
Suppose that $S=E_6$, $E_7$, or $E_8$, and let $\{\e_i\}$ be the
standard coordinate functions on $\Def (S)$.
Define
$\widetilde{\e}_i = \e_i$ if $S \ne E_7$ or $i \ne 12, 18$, and in
the case of $E_7$, define
\[\widetilde{\e}_{12} = \e_{12} + \frac13 \e_6^2
\quad  \text{and} \quad
\widetilde{\e}_{18} = \e_{18} + \frac13\e_6\e_{12} + \frac{2}{27}\e_6^3.\]
Let $M_f$ be the set of monomials $T^dY^kZ^{\ell}$ such that
$\e_iY^kZ^{\ell}$ is one of the terms in the polynomial in preferred versal
form, and $\nu_f(\widetilde{e}_i)=d<\infty$.
If any of the monomials in $M_f$ are listed in the right half of
table~\ref{tableMONOrev}, then the general hyperplane
section of $X$ has singularity type at worst $S$, where $S$ is the
label on the leftmost column in the right half of the
table which contains some monomial from $M_f$.

{\parindent 1.5em
Moreover, if $X$ has an irreducible small resolution of
length $\ell$, then $M_f$ contains no monomials
to the left of the column whose label is the associated type $S_{\ell}$.
}

\end{enumerate}
\end{proposition}

The left half of table~\ref{tableMONOrev} has been included to make it
easier to find which monomials in $M_f$ come from which standard coordinate
functions $\e_i$.  It is not actually necessary for the description of
the link between $M_f$ and the singularity type of the general hyperplane
section.

{\renewcommand{\arraystretch}{1.4}

\begin{table}[b]
\begin{center}
\begin{tabular}{|c|c|c||c|c|c|c|c|c|c|} \hline
$E_6$     & $E_7$ & $E_8$ & $A_0$ & $A_1$ & $A_2$ & $D_4$ & $D_k$  & $E_6$ &
$E_7$  \\ \hline
          &           & $\e_{8}$  & &    & &         & &          & $T Y Z^2$
\\
          & $\e_{6}$  &           & &    & &         &  & $T Y^2$     & \\
          &           & $\e_{12}$ & &    & &         & & $T Z^3$      & \\
$\e_{5}$  & $\e_{8}$  & $\e_{14}$ & &    & & $T Y Z$ & &          & $T^2 Y Z$
\\
$\e_{6}$  & $\e_{10}$ & $\e_{18}$ & &    & & $T Z^2$ & & $T^2 Z^2$    & \\
$\e_{8}$  & $\widetilde{\e}_{12}$ & $\e_{20}$ & & $T Y$ & & & $T^2 Y$ &
& $T^3 Y$ \\
$\e_{9}$  & $\e_{14}$ & $\e_{24}$ & & $T Z$ & & $T^2 Z$ & & $T^3 Z$   & \\
$\e_{12}$ & $\widetilde{\e}_{18}$ & $\e_{30}$ & $T$ & & $T^2$ & & $T^3$ & $T^4$
& \\ \hline
\end{tabular}
\end{center}

\medskip

\caption{}
\label{tableMONOrev}
\end{table}
}

\begin{pf}
The proof is based on the classification of rational double points
by means of their Newton polygons.  A convenient reference for this is
\cite[(4.9)(3)]{[YPG]}.  In brief,  suppose that
$\{f=0\} \subset \C ^3$ has a
rational double
point at the origin.  Then the  type $S$ is determined by the defining
polynomial
 $f$ as follows.
\begin{enumerate}
\item
$S=A_0$ (i.e. the surface is smooth at the origin) if and only if
$f$ contains a linear term.  (Notice that $f$ contains no constant
term, since the origin lies on the surface.)

\item
$S=A_1$ if and only if the quadratic part $f_2$ of $f$ has rank 3.

\item
$S=A_{n-1}$, $n>2$ if and only if the quadratic part $f_2$ of $f$
has rank 2.  The value of $n$ is determined by the higher order terms.
In particular, if the cubic part
$f_3$ of $f$ is nonzero and involves none of the variables appearing in
$f_2$, then $S=A_2$.

\item
If $f_2$ has rank 1, choose coordinates so that $f = x^2 + g(y,z)$,
and $g$ has no quadratic part.  Note that $g_3$ is a homogeneous
cubic in 2 variables.
\begin{enumerate}
\item
$S=D_4$  if and only if the cubic part $g_3$ of $g$ has three distinct
linear factors.

\item
$S=D_n$, $n>4$  if and only if the cubic part $g_3$ of $g$ has two
distinct linear factors.  (The value of $n$ is determined by the higher
order terms.)

\item
If $g_3$ has a unique linear factor, write $g_3 = h^3$.
\begin{enumerate}
\item
$S=E_6$  if and only if $h$ does not divide the quartic part $g_4$
of $g$.

\item
$S=E_7$  if and only if the quartic part $g_4$ of $g$ is divisible by
$h$ but not by $h^2$.

\item
$S=E_8$ otherwise.

\end{enumerate}
\end{enumerate}
\end{enumerate}

The threefold $X$ has defining polynomial $\mu_f^*(\Phi_S)$,
where $\Phi_S$ is the polynomial in preferred versal form of type $S$.
If $\f_i Y^k Z^{\ell}$ is a term in $\Phi_S$ and if $\nu_f(\f_i)=d<\infty$,
then the monomial $T^d Y^k Z^{\ell}$ appears in the defining polynomial of $X$
with
a nonzero coefficient.  (Here, $T$ is the coordinate on the disk $\Delta$.)
In this way, we can analyze the low-degree terms
appearing in the defining polynomial of $X$ by using the set $M_f$.

We define the  {\em leading terms\/} of $\Phi_S$  to be those
which have constant coefficients; there are  2 or 3 such terms.
For all other coefficients $\f_i$, we have $\nu_f(\f_i) \ge 1$.
Table~\ref{tableMONOrev} has been constructed so that all potential
low-degree terms (other than leading terms) are shown there.

Suppose first that $S=A_{n-1}$.  The only possible monomial of degree 1
in the defining polynomial is $T$, and this occurs if and only if
$\nu_f(\a_n)=1$.
This is the condition for $X$ (and its general hyperplane section) to
be smooth at the origin; in this case, the singularity type is certainly
``at worst" $A_1$.

If $\nu_f(\a_n)>1$, we consider quadratic terms.  The leading term of
degree 2 is $-XY$, while other potential quadratic terms must be
chosen from $\{TZ,T^2\}$.  If at least one of those potential terms occurs
with a nonzero coefficient, then the rank of the quadratic part of
the defining polynomial is at least 3.  And this implies that the quadratic
part of
the general hyperplane section will have rank 3, and so will have a
 singularity of type $A_1$.  But to guarantee that at least one of the
potential
terms occurs, we simply need   $\nu_f(\a_{n-1})=1$ or $\nu_f(\a_n)=2$.

Suppose next that $S=D_n$.  Again, the only possible monomial of degree 1
in the defining polynomial is $T$, and this occurs if and only if
$\nu_f(\d_{2n-2})=1$.
This is the condition for $X$ (and its general hyperplane section) to
be smooth at the origin; in this case, the singularity type is certainly
``at worst" $A_1$ or even $D_4$.

If $\nu_f(\d_{2n-2})>1$, we consider quadratic terms.  The leading term of
degree 2 is $X^2$, while other potential quadratic terms must be
chosen from $\{TY,TZ,T^2\}$.  If at least one of the terms $TY$, $TZ$ occurs
with a nonzero coefficient, then the rank of the quadratic part of
the defining polynomial is at least 3.  And this implies that the quadratic
part of
the general hyperplane section will have rank 3, and so will have a
 singularity of type $A_1$.  But to guarantee that at least one of those
terms occurs, we simply need  $\nu_f(\c_n)=1$ or  $\nu_f(\d_{2n-4})=1$.

On the other hand, if neither of those terms occurs, yet $T^2$ occurs,
then the rank of the quadratic part is 2.  (This happens when $\nu_f(\c_n)>1$,
$\nu_f(\d_{2n-4})>1$, and $\nu_f(\d_{2n-2})=2$.)  Now the leading
term $Y^2Z$ also appears in our defining polynomial.  Since this is a
term of degree 3 which involves neither of the variables
$X$, $T$ which appear in the quadratic part, there will be hyperplane
sections of type $A_2$.  (For example, the hyperplane section defined
by $Y=Z$ will be of type $A_2$.)
It follows that the singularity type of the
general hyperplane section is at worst $A_2$,
and hence is certainly at worst $D_4$.

So we may assume that $\nu_f(\c_n)>1$,
$\nu_f(\d_{2n-4})>1$, and $\nu_f(\d_{2n-2})>2$.  The defining polynomial can
then
be written in the form $X^2 + G(Y,Z,T)$, and the cubic part of $G$ takes
the form $G_3(Y,Z,T) = Y^2Z + T \cdot H(Y,Z,T)$.  Moreover, the only
monomials which could appear in $T \cdot H(Y,Z,T)$ are $TZ^2$, $T^2Y$,
$T^2Z$, $T^3$.  In order for the general hyperplane section of $G_3$
to fail to have 3 distinct linear factors, $G_3$ must be nonreduced.
Since $G_3=0$ defines a plane cubic, this implies that $G_3$ itself factors
in the form $H_1^2 H_2$, where $H_1$, $H_2$ are two
(possibly equal) linear polynomials.

Considering this factorization mod $T$, we see that it must take the
form
\[G_3(Y,Z,T) = (Y + \alpha T)^2 \cdot (Z + \beta T).\]
  Moreover,
since $TY^2$ is not one of the monomials which can occur in $G_3$,
$\beta$ must in fact be 0.  Thus, if the general hyperplane section
fails to have type $D_4$, $Z$ must divide $G_3$.  The presence of
either of the monomials $T^2Y$ or $T^3$ will prevent this, and their
presence is guaranteed by the conditions
$\nu_f(\c_{n})=2$ and $\nu_f(\d_{2n-2})=3$, respectively.  So when
either of these conditions holds, the general hyperplane section must
be of type $D_4$.

Suppose finally that $S=E_n$.  We define $\widetilde{X}=X+\frac12Z^2$
in case $E_6$ and $\widetilde{X}=X$
in cases $E_7$ and $E_8$.  Then the leading terms take the form
 $- \widetilde{X}^2 +\frac14 Z^4 + Y^3$,
$-\widetilde{X}^2-Y^3+16YZ^3$, $-\widetilde{X}^2+Y^3-Z^5$, respectively.
The monomials in these leading terms together with the monomials in
$M_f$ will include all monomials of low degree in the defining polynomial
$\mu_f^*(\Phi_S)$.

The analysis of  the cases with a linear part or with a quadratic part of
rank bigger than 1 proceeds almost exactly as in the case of $D_n$, using
$-\widetilde{X}^2$ in place of $X^2$ for the leading term of degree 2, and
$\pm Y^3$ for the leading term of degree 3.  It yields the criteria
for having  hyperplane sections of types $A_0$, $A_1$, and $A_2$ which
are stated in  table~\ref{tableMONOrev}.

(The only remarks that need to be added to the argument given in the $D_n$
case concern the $E_7$ case, since we use two modified coefficients
$\widetilde{\e}_{12}$ and $\widetilde{\e}_{18}$ in that case.  The remarks
(which follow from the defining formulas for the modified coefficients)
are that $\nu_f(\widetilde{\e}_{12})=1$ if and only if $\nu_f(\e_{12})=1$,
and that when $\nu_f(\e_{12})>1$, we have
$\nu_f(\widetilde{\e}_{18})=2$ if and only if $\nu_f(\e_{18})=2$.
Thus, the orders of vanishing which predict the presence of the monomials
$TY$ and $T^2$ are being calculated properly.)

So we may assume that none of the monomials $T^2$, $TY$, $TZ$, or $T$
appear in our defining polynomial.  The defining polynomial can be written
in the form $-\widetilde{X}^2 + G(Y,Z,T)$, and this time
the cubic part of $G$ takes
the form $G_3(Y,Z,T) = \pm Y^3 + T \cdot H(Y,Z,T)$.  As before,
the general hyperplane section will be of type $D_4$ unless
$G_3$ can be factored
in the form $\pm H_1^2 H_2$, where $H_1$, $H_2$ are two
(possibly equal) linear polynomials.

Considering this factorization mod $T$, we see that
it takes the form
\[\pm G_3(Y,Z,T) = (Y + \alpha T)^2 \cdot (Y + \beta T).\]
  In particular,
for such a factorization to exist, $G_3$ must be a function of $Y$ and
$T$ alone.  The presence of any of the monomials $TYZ$, $TZ^2$, or $T^2Z$
in $M_f$ prevents this, and forces the general hyperplane section to have
type $D_4$.

If none of those monomials is present in $M_f$, then $G_3$ is a homogeneous
binary cubic, and the general hyperplane section
has type $D_k$  unless $G_3$ is the cube of
a linear polynomial.  (The type will be $D_4$ if there are three distinct
linear factors of  $G_3$, and will be $D_k$, $k>4$ if there are only two.)
Now for $S \ne E_7$, the monomial $TY^2$ cannot occur in $G_3$ (as is
clear from table~\ref{tableMONOrev}).  In this
case, if $G_3$ is a cube it must be $Y^3$, and the presence of either of
the monomials $T^2Y$ or $T^3$ in $M_f$ will prevent this from happening,
and lead to the general hyperplane section having type $D_k$.

The argument is more complicated in the case of
$E_7$.\footnote{The argument could have been simplified, eliminating the
use of $\widetilde{\e}_{12}$ and $\widetilde{\e}_{18}$, were it not
for our desire to match notation with Bramble \cite{[Bra]}.}
We define
$\widetilde{Y} = Y - \frac13 \e_6$, and note that
\[-\widetilde{Y}^3 + \widetilde{\e}_{12} \widetilde{Y} + \widetilde{\e}_{18}
= -Y^3 + \e_6 Y^2 + \e_{12} Y + \e_{18}.\]
In this case, if $G_3$ is a cube, then its cube root $H$ must be
the linear part
(with respect to $Y$, $T$)
of $- \widetilde{Y}$.  Thus, if either
$\nu_f(\widetilde{\e}_{12})=2$ or $\nu_f(\widetilde{\e}_{18})=3$,
then $G_3$ cannot be a cube, and the singularity type must be $D_k$.

We now assume that $G_3$ is in fact a cube, and let $H$ be its cube
root.  If $S=E_6$, then the general hyperplane section is at worst
 $E_6$ (which is certainly at worst $E_7$), and we are finished.

Suppose instead that $S=E_7$ and $\nu_f(\e_6)=1$.  If
$\e_6 \equiv \alpha T \mod{T^2}$, then $H = - Y + \frac13 \alpha T$.
The quartic part of $G$ includes the leading term $16YZ^3$.  But since
the monomial $TZ^3$ cannot occur in $G$ and $\alpha \ne 0$,
it follows that the quartic part of $G$ cannot be divisible by $H$.
Thus, the general hyperplane section has type $E_6$.

We may therefore assume that either $S = E_7$ and $\nu_f(\e_6)>1$, or that
$S=E_8$.  In
either case, the cube root $H$ is exactly $\pm Y$.  This divides the
quartic part of the leading term (which is $16YZ^3$ or 0, respectively).
We can therefore identify which cases have general hyperplane
section $E_6$ or $E_7$ by finding in $M_f$ a monomial not divisible by $Y$,
or one not divisible by $Y^2$, respectively.  This is exactly what
is done in the final two columns of table~\ref{tableMONOrev}.

To prove the last statement in the proposition,
note that all labels $L$ to the left of $S_{\ell}$
in table~\ref{tableMONOrev} correspond to singularities with the property that
the
maximum multiplicity
which occurs in the maximal ideal cycle for the singularity is strictly
less than $\ell$. In addition, any singularity whose type is at worst
 $L$ has this same property.  But if $X$ has an irreducible small
resolution of length $\ell$,
 no such singularity can be the general hyperplane
section. Thus, there can be no monomials in columns to the left of that
labeled by $S_{\ell}$.
\end{pf}

\begin{lemma} \label{lem83}
Let $\pi\colon Y \to X$ be an irreducible small resolution
of an isolated Gorenstein threefold singularity,  and let
$\{f=0\}$ define a hyperplane section with a rational double point.
Let $\mu_{f \circ \pi}\colon \Delta \to \PRes (S,v) = V/\W_0$
be the classifying map determined by $f$, and define a discrete valuation
$\nu_{f \circ \pi}\colon \C [V]^{\w_0} \to \Z $
by
\[\nu_{f \circ \pi}(\f) =
\text{order of vanishing at $0$ of $\mu^*_{f \circ \pi}(\f)$}.\]
If
$\f_i \in \C [V]^{\w_0}$ is any standard coordinate function on
$\PRes (S,v)$  , then
$\nu_{f \circ \pi}(\f_i) \ge 1$.  If $\f_i$ is in fact a ``constant
term'', then $\nu_{f \circ \pi}(\f_i) = 1$.
\end{lemma}

\begin{pf}
The first assertion holds since $\mu_{f\circ\pi}(0)=0$.  The
second holds since $Y$ is smooth at the singular point associated with $\f_i$.
\end{pf}

Consider again the inclusion of rings
$\C [V]^{\w} \subset \C [V]^{\w_0}$,
which corresponds to the natural projection
$\sigma\colon \PRes (S,v) \to \Def (S)$.
The larger
ring $\C [V]^{\w_0}$ is a free polynomial ring generated by the standard
coordinate functions on $\PRes (S,v)$.  Thus, each element of the smaller
ring $\C [V]^{\w}$ can be written as a polynomial in those standard
coordinate functions.  Each monomial in such an expression has a
{\em degree\/} (in the standard coordinate functions
which are generating the ring)
as well as a {\em weight\/} under the background $\C ^*$-action.
For a fixed weight $i$ and degree $d$, we denote by $P_{i,d}$
the subspace of polynomials in
$\C [V]^{\w_0}$ whose weight is $i$ and whose degree is less than
or equal to $d$.

\begin{lemma} \label{lem84}
Fix a partial resolution type $(S,v)$.
Let $\pi\colon Y \to X$ be an irreducible small resolution of an isolated
Gorenstein
threefold singular point, and let $\{f=0\}$ be a hyperplane section
with partial resolution type $(S,v)$.  Let
$\nu_f\colon \C [V]^{\w} \to \Z $ be the
associated discrete valuation.

Suppose that $\f_i \in \C [V]^{\w}$ is a function
of weight $i$, and
 $m_i \in P_{i,d}$ is a monomial of degree exactly $d$.
Suppose further that there is  an ideal
$I \subset \C [V]^{\w_0}$ whose intersection with $P_{i,d}$ is $\{0\}$
such that
\[\f_i \equiv c \cdot m_i \mod I\]
for some nonzero constant $c$.
If $\nu_{f\circ\pi}(m_i) = d$, then $\nu_f(\f_i) = d$.
\end{lemma}

\begin{pf}
Consider the classifying map
$\mu_{f}\colon \Delta \to \Def (X_0) = \Def (S) = V/\W$
determining our threefold $X$, and its associated discrete valuation
$\nu_{f }\colon \C [V]^{\w} \to \Z $, as well as the related
map and valuation $\mu_{f \circ \pi}$ and $\nu_{f \circ \pi}$.   Note that
$\nu_{f\circ\pi}$ extends $\nu_f$.
Suppose that $\f_i - c\cdot m_i$ has degree at most $d$.
Then it lies in $I \cap P_{i,d}$ and so must be $0$.  Thus, we may
assume that
$\f_i-c\cdot m_i$ has degree strictly greater than $d$.  It follows by
lemma~\ref{lem83} that
$\nu_{f \circ \pi}(\f_i - c\cdot m_i) > d$.    Since
$\nu_{f \circ \pi}(c\cdot m_i) = d$, it follows that $\nu_{f }(\f_i) = d$.
\end{pf}

\begin{pf*}{Proof of the Main Theorem}
If $\pi\colon Y\to X$ is as above, we set about
showing that the orders of the standard coordinate functions on
$\Def (S)$ (or the
expressions
$\widetilde{\e}_i$ in the case of $E_7$) satisfy the relevant hypothesis from
proposition~\ref{prop81}.
We show that a
particular monomial occurs in the defining polynomial of $X$ with nonzero
coefficient,
then we use lemmas~\ref{lem83} and \ref{lem84} together with
proposition~\ref{prop81}
to conclude that the general hyperplane section of $X$ is as
claimed.
We assume henceforth that $(S,v_k)$ is {\em not\/} one of the pairs illustrated
in figure~\ref{figure1}.
This puts certain restrictions on $n$ and $k$ which we will exploit
without comment.

The first case to consider is $S=A_{n-1}$.
The standard coordinate functions on $\PRes (A_{n-1},v_k)$ are
$\mu_1, {\a'_2}, \ldots, {\a'_k}, {\a''_2}, \ldots, {\a''_{n-k}}$.
By proposition~\ref{prop72},
$\a_{n-1} \equiv
\a'_{k-1} \a''_{n-k} + \a'_k \a''_{n-k-1}\mod{\m_1}$, and
$\a_{n} \equiv \a'_k  \a''_{n-k}\mod{\m_1}$.

Suppose that $k=1$ and $n \ge 2$.  The length is 1, and
there is only one ``constant term" in this case:
$\a''_{n-1}$.  (Note that there would be no constant terms whatsoever
in case $n=1$.)  Since $\a'_0=1$ and $\a'_1=0$ by definition, we have
\[\a_{n-1} \equiv {\a''_{n-1}}  \mod{(\mu_1)}.\]
Furthermore $P_{n-1,1} = \C  \cdot \a''_{n-1}$, so its intersection with
the ideal $I = (\m_1)$ is $\{0\}$.  (Again we have used $n \ge 2$.)
By lemma~\ref{lem83}, $\nu_f(\a''_{n-1}) = 1$.
By lemma~\ref{lem84}, we conclude that
$\nu_f(\a_{n-1}) = 1$.  By proposition~\ref{prop81}, the singularity type
of the general hyperplane
section is at worst $A_1$.

Suppose instead that $1 < k < n-1$.  (We can omit the case $k=n-1$ by
symmetry.)  In this case, the length is 1,
there are two ``constant terms"
$\a'_k$ and $\a''_{n-k}$, and we have
\[\a_{n} \equiv \a'_k  \a''_{n-k}\mod{\m_1}.\]
If $\f \in (\mu_1) \cap P_{n,2}$ then $\f = \mu_1 \psi$ for some
$\psi \in P_{n-1,1}$.  But
since the maximum weight among the standard coordinate functions
on $\PRes (A_{n-1},v_k)$ is
$\max \{1,k,n-k\} \le n-2$, there can be no affine linear function
in these variables of weight $n-1$.  Thus $\psi = 0$, so $\f = 0$ and
the intersection of $P_{n,2}$ with $(\m_1)$ must be $\{0\}$.

By lemma~\ref{lem83}, $\nu_f(\a'_k\a''_{n-k})=2$.
By lemma~\ref{lem84}, we conclude that
$\nu_f(\a_{n}) = 2$.  By proposition~\ref{prop81}, it follows that the
singularity type of the general
hyperplane section is at worst $A_1$.

The next case to consider is $S=D_{n}$, $k \le n-2$, with $n \ge 4$.
The standard coordinate functions on $\PRes (D_{n},v_k)$ are
$\mu_1, {\a'_2}, \ldots, {\a'_k},
{\d''_2}, \ldots, {\d''_{2n-2k-2}}, {\c''_{n-k}}$.
By proposition~\ref{prop72}, we have
$\c_n \equiv \a'_k  \c''_{n-k}\mod{\m_1}$, and
$\d_{2n-2} \equiv  (\a'_k)^2 \d''_{2n-2k-2}\mod{J_4}$,
where $J_4$ is the ideal generated by all
monomials of degree 4 in the standard coordinate functions on
$\PRes (D_n,v_k)$.

Suppose that $k=1$.  The length is 1 and
there is only one ``constant term" in this case:
$\d''_{2n-4}$.  Moreover, proposition~\ref{prop72} provides us with
an additional congruence in this case:
\[\d_{2n-4} \equiv \d''_{2n-4}\mod{\m_1}.\]
We have $P_{2n-4,1} = \C  \cdot  \d''_{2n-4}$, whose intersection with
the ideal $I=(\m_1)$ is $\{0\}$.
By lemma~\ref{lem83}, $\nu_f(\d''_{2n-4})=1$.
By lemma~\ref{lem84}, we conclude that
$\nu_f(\d_{2n-4}) = 1$.  By proposition~\ref{prop81}, the singularity
type of the general hyperplane
section is at worst $A_1$.

If $1< k< n-1$, the length is 2.  Moreover, proposition~\ref{prop72} gives
\[\d_{2n-2}\equiv (\a_k')^2\d_{2n-2k-2}''\mod{J_4}.\]
  It is clear from
considering degrees that $P_{2n-2,3}$ intersects $J_4$ trivially.

If $k < n-2$, then there are two constant terms $\a'_k$ and $\d''_{2n-2k-2}$.
By lemma~\ref{lem83}, $\nu_f((\a'_k)^2\d''_{2n-2k-2})=3$.
By lemma~\ref{lem84}, we conclude that
$\nu_f(\d_{2n-2}) = 3$.  By proposition~\ref{prop81}, it follows that the
singularity type of the
general hyperplane section is at worst $D_4$.

If $k=n-2$, then by lemma~\ref{lem72} the constant terms are
$\a'_{n-2}$,  $-\frac14(\d''_2-2\c''_2)$, and $-\frac14(\d''_2+2\c''_2)$.
According to lemma~\ref{lem83}, each of these has order 1 with
respect to the discrete valuation $\nu_f$; hence, at least one of
$\nu_f(\c''_2)$ and $\nu_f(\d''_2)$ is also equal to 1.

If $\nu_f(\d''_2)=1$, then $\nu_f((\a'_k)^2\d''_{2n-2k-2})=3$ and we
can use the same argument as in the case $k<n-2$ to conclude that
the general hyperplane section has type at worst $D_4$.  On the other hand,
if $\nu_f(\c''_2)=1$, we use the additional congruence
\[\c_n \equiv \a'_{n-2}  \c''_{2}\mod{\m_1}\]
provided by proposition~\ref{prop72}.
We have $P_{n,2} \cap (\m_1) = \{0\}$
since the maximum weight among standard coordinate functions
on $\PRes (D_n,v_{n-1})$ is $n-2$,
and we have $\nu_f(\a'_{n-2}  \c''_{2}) = 2$.  By lemma~\ref{lem84},
we conclude that
$\nu_f(\c_n) = 2$.  By proposition~\ref{prop81}, it follows that the
singularity type of the general
hyperplane section is at worst $D_4$.

The third case to consider is $S=D_{n}$, $k = n$, in which the length is 1.
  (By symmetry, we can omit
the case $k=n-1$.)
The standard coordinate functions on $\PRes (D_{n},v_n)$ are
$\mu_1, \a'_2, \ldots, \a'_n$, and $\a'_n$ is the unique ``constant term".
By proposition~\ref{prop72}, we have
\[\c_n \equiv \a'_{n}\mod{\m_1},\]
and lemma~\ref{lem83} implies that $\nu_f(\a'_n)=1$.
Moreover, since $P_{n,1} = \C \cdot\a'_n$, its intersection with $(\m_1)$
is $\{0\}$.
By lemma~\ref{lem84}, we conclude that
$\nu_f(\c_n) = 1$.  By proposition~\ref{prop81}, it follows that the
singularity type of the general
hyperplane section is at worst $A_1$.

{\renewcommand{\arraystretch}{2}

\begin{table}[p]
\begin{center}
\begin{tabular}{|c|c|c|c|c|c|c|} \hline
$(S,v_k)$  & length & $S'$ & $S''$ & congruence & monomial & type \\ \hline
\hline
$(E_6,v_0)$  & 2 & $A_5$ & -- &
     $\e_6 \equiv - {\a_6}' \mod I$
& $TZ^2$ & $D_4$ \\ \hline
$(E_6,v_4)$  & 2 & $E_4$ & $A_1$ &
     $\e_5 \equiv - {\e_5}' \mod I$
& $TYZ$ & $D_4$ \\ \hline
$(E_6,v_5)$  & 1 & $E_5$ & $A_0$ &
     $\e_8 \equiv - \frac14 {\e_8}' \mod I$
& $TY$ & $A_1$ \\ \hline
\hline
$(E_7,v_0)$  & 2 & $A_6$ & -- &
     $\e_{14} \equiv 64 ({\a_7}')^2 \mod I$
& $T^2Z$ & $D_4$ \\ \hline
$(E_7,v_1)$  & 2 & $D_6$ & -- &
     $\e_{10} \equiv 16 {\d_{10}}' \mod I$
& $TZ^2$ & $D_4$ \\ \hline
$(E_7,v_2)$  & 3 & $A_1$ & $A_5$ &
    $\e_{6} \equiv - 12 {\a_6}'' \mod I$
& $TY^2$ & $E_6$ \\ \hline
$(E_7,v_4)$  & 3 & $E_4$ & $A_2$ &
     $\e_{10} \equiv 16 ({\e_5}')^2 \mod I$
& $T^2Z^2$ & $E_6$ \\ \hline
$(E_7,v_5)$  & 2 & $E_5$ & $A_1$ &
     $\e_{8} \equiv - 4 {\e_8}' \mod I$
& $TYZ$ & $D_4$ \\ \hline
$(E_7,v_6)$  & 1 & $E_6$ & $A_0$ &
     $\e_{12} \equiv 16 {\e_{12}}' \mod I$
& $TY$ & $A_1$ \\ \hline
\hline
$(E_8,v_0)$  & 3 & $A_7$ & -- &
    $\e_{24} \equiv  ({\a_8}')^3 \mod I$
& $T^3Z$ & $E_6$ \\ \hline
$(E_8,v_1)$  & 2 & $D_7$ & -- &
    $\e_{24} \equiv
      - \frac{1}{16} ({\d_{12}}')^2 \mod I$
& $T^2Z$ & $D_4$ \\ \hline
$(E_8,v_2)$  & 4 & $A_1$ & $A_6$ &
    $\e_{14} \equiv  ({\a_7}'')^2 \mod I$
& $T^2YZ$ & $E_7$ \\ \hline
$(E_8,v_5)$  & 4 & $E_5$ & $A_2$ &
    $\e_{8} \equiv - \frac14 {\e_8}' \mod I$
& $TYZ^2$ & $E_7$ \\ \hline
$(E_8,v_6)$  & 3 & $E_6$ & $A_1$ &
    $\e_{12} \equiv  {\e_{12}}' \mod I$
& $TZ^3$ & $E_6$ \\ \hline
$(E_8,v_7)$  & 2 & $E_7$ & $A_0$ &
    $\e_{18} \equiv
          \frac{1}{64} {\e_{18}}' \mod I$
& $TZ^2$ & $D_4$ \\ \hline
\end{tabular}
\end{center}

\medskip

\caption{Key computations}
\label{table-key1}
\end{table}
}

Finally, we consider the cases with $S=E_n$.  Among the standard coordinate
functions
on  $\PRes (S,v_k)$, let $\widetilde{\f}_N$ be
the ``constant term" of highest weight, say weight $N$.
(This is unique, since we are
avoiding the case $(S,v_k)=(E_6,v_3)$.)  Let $I$ be the ideal
in $\C [V]^{\w}$ which is generated by all the standard coordinate
functions on $\PRes (S,v_k)$ other than $\widetilde{\f}_N$.  We select a
standard
coordinate function $\e_i$ as indicated in table~\ref{table-key1},
and calculate it mod $I$ using the relations given in
proposition~\ref{prop71}.  Table~\ref{table-key1} shows the results
of this calculation:  we will describe the calculation itself in
section~10.

The calculated result takes the form
\[\e_i \equiv c \cdot (\widetilde{\f}_N)^d \mod{I}\]
for some nonzero constant $c$.
(The key point of the calculation is showing that this constant is
not 0.)   Moreover, since $i=Nd$ and $N$ is the highest
weight among standard coordinate functions on $\PRes (S,v_k)$,
any monomial of weight $i$ other than $(\widetilde{\f}_N)^d$ must have degree
strictly greater than $d$.  It follows that $I \cap P_{i,d} = \{0\}$
and thus by lemma~\ref{lem84}, $\nu_f(\e_i)=d$.

Let $M_f$ be the set of monomials from proposition~\ref{prop81}.
Since $\nu_f(\e_i)=d$, we conclude that $M_f$ contains the monomial
shown in the next-to-last column of table~\ref{table-key1}.
The label $L$ of the column in table~\ref{tableMONOrev} in which
that monomial appears has been reproduced in the last column of
table~\ref{table-key1}.
In each case, the label $L$ coincides with the type $S_{\ell}$ which
is associated with the length $\ell$.  (The length itself is shown in the
second column.)  Thus, by proposition~\ref{prop81},
the set $M_f$ can
contain no monomials to the left of the column labeled by $L$ in
table~\ref{tableMONOrev}, and the singularity type of
the general hyperplane section is at worst $L=S_{\ell}$.
The main theorem then follows from lemma~\ref{lem82}.
\end{pf*}

\section{The computation of preferred versal form in the $E_n$ cases.}

In this section, we explain how to explicitly compute a defining polynomial
in  preferred versal form for the standard deformation
 in the $E_n$ cases.  The result will be a formula for
the standard coordinate functions $\e_i$ in terms of the elementary symmetric
functions
$s_1,\ldots,s_n$ of the distinguished functionals $t_1,\ldots,t_n$.  One's
natural inclination is to expand everything completely as polynomials
in $s_1,\ldots,s_n$, and simply work with things in expanded form.  For
this computation, however, that would not be a wise strategy:
the ``constant term" in the case of $E_8$ is a polynomial with
2462 terms.\footnote{It is possible that the coefficients of a few of
these terms may be zero.}

To keep things explicit, and yet in a compact format, we introduce two
notions.  By a {\em set of substitution rules}, we mean a set of expressions
of the form $v_i = f_i(x_1,\ldots,x_k)$ which express certain variables
$v_i$ as polynomial functions of other variables $x_j$.
(This notion is more flexible than the notion of a ring homomorphism,
since the rings to which $v_i$ and $x_j$ belong do not need to be
specified until the substitution rules are used.)
We say that a set of
substitution rules ${\cal R}$ is given in {\em solve-list format\/} when it is
specified by three objects ${\cal R}'$, $P$, and $L$ which satisfy a certain
condition, as follows.
${\cal R}'$ is another set of substitution
rules, $P$ is  a polynomial,
and $L$ is an ordered list, the {\em solve-list}, consisting of
pairs $(m_i,v_i)$ where $m_i$ is a monomial and $v_i$ is a variable.
(The substitution rules ${\cal R}'$ take the form
$w_{\alpha} = g_{\alpha}(y_1,\ldots,y_{\ell})$, where the $w_{\alpha}$ are
distinct from the $v_i$ but the $y_{\beta}$ may include some $v_i$'s.)
The condition which must be satisfied is this:  if $c_i$ is the
coefficient of the monomial $m_i$ in the expression obtained by
substituting the rules ${\cal R}'$ into the polynomial $P$, then the
variable $v_i$ appears linearly in $c_i$ with a nonzero constant
coefficient, and does not appear in any $c_j$ with $j<i$.
(It is allowed,
however, that $v_i$ appear in $c_j$ with $j>i$, and it may even appear
in a nonlinear fashion.)

The algorithm for producing the rules ${\cal R}$ from the triple ${\cal R}'$,
$P$, $L$
is simple:  compute the coefficients $c_i$, and for $i = 1, 2, \ldots ,$
successively solve the equations $c_i = 0$ for the variables $v_i$.
In solving the $i^{\text{th}}$ equation, one uses previously found values
of $v_j$ for $j<i$ to eliminate the variables $v_j$ from the expression
$c_i$.
We refer to this process as {\em expanding the solve-list}.

Specifying the rules ${\cal R}$ by means of ${\cal R}'$, $P$, and $L$ (without
actually
expanding the solve-list) can give a relatively compact representation of a
complicated set of
rules.  Moreover, our later application of these explicit calculations
will be of the following form:  calculate what happens when the rules
${\cal R}$ are restricted to a subspace which is parametrized in a simple
way.
In carrying out those applications, it will be much to our advantage to
begin by pulling back the ingredients of the solve-list format
(i.e. ${\cal R}'$, $P$, and $L$) to the parameter space,
and then computing the pullback of the rules ${\cal R}$ by directly
expanding the pulled-back solve-list.
  This part of the computation is explained in
section 10.

We have already encountered (during the proof of proposition~\ref{prop51})
 a set of substitution rules which is best
expressed in solve-list format.  Let ${\cal R}_{\mu}$ denote the set of
substitution rules
given in equation (\ref{eq23}),
which
describes a change of generators in the algebra $L$.  Then the substitution
${\cal R}_{\psi}$ (which describes the coefficients to use in ${\cal R}_{\mu}$
which will
produce the preferred versal form) can be given in solve-list format by
means of the substitution rules
${\cal R}_{\mu}$, the polynomial  $\bar{\Phi}_{E_n}$, and the solve-list
\refstepcounter{equation}\label{eq29}
\begin{gather*}
\begin{tabular}{|c|c|c|c|c|c|} \hline
$Y^2Z$   & $XYW$   & $Z^3$     & $XZW$      & $Y^2W$   & $XW^2$ \\ \hline
$\psi_1$ & $\psi_2$ & $\psi'_3$ & $\psi''_3$ & $\psi_4$ & $\psi_6$ \\ \hline
\end{tabular}
\tag{\ref{eq29}a}\\[1.5ex]
\begin{tabular}{|c|c|c|} \hline
$Z^4$     & $YZ^2W$   & $Z^3W$    \\ \hline
$\psi_2$  & $\psi_4$  & $\psi_6$  \\ \hline
\end{tabular}
\tag{\ref{eq29}b}\\[1.5ex]
\begin{tabular}{|c|c|c|} \hline
$Y^2ZW$   & $Z^4W^2$   & $Y^2W^2$    \\ \hline
$\psi_4$  & $\psi_6$  & $\psi_{10}$  \\ \hline
\end{tabular}
\tag{\ref{eq29}c}
\end{gather*}
(We continue to use the convention that equation numbers which are followed by
$a$, $b$, or $c$ refer to the cases of $E_6$, $E_7$, or $E_8$, respectively.)
We computed the corresponding expressions $c_i$ explicitly in
equation (\ref{eq51}),
and verified that $\{ c_i=0\}$ has the appropriate triangular
form.

\bigskip

Throughout this section, the subscript on a variable
(when present)
indicates
its weight with respect to the $\C ^*$-action.
We retain the notation introduced in section 5.

The first step in our computation is to express the anti-pluricanonical
mappings explicitly in coordinates, and thereby obtain a good generating set
$\bar{X}$, $\bar{Y}$, $\bar{Z}$, $\bar{W}$ for $L$.
We identify these generators with polynomials in $\C [V][x,y,z]$
which satisfy certain base conditions.
The anti-canonical mapping (which corresponds to $L_1$)
is given by the cubics passing through the zero-cycle
$\eta(t_1)+\cdots+\eta(t_n)$.  In other words, we want cubics
$F \in \C [V][x,y,z]$ such that $\Psi_n(U)$ divides $\eta^*(F)$,
where
$\Psi_n(U) \coloneq U^n - s_1 U^{n-1} + \dots + (-1)^n s_n$
is the monic polynomial of degree $n$
whose roots are $t_1,\ldots,t_n$.
It is not difficult
to find a basis for these in each case (by hand), and to make the basis
 match the
first part of the normalizations established in table~\ref{table45}.
In the case of $E_6$, we get
\refstepcounter{equation}\label{eq1}
\begin{align*}
\begin{split}
\bar{W} & \coloneq   x^3 - y z^2 \\
\bar{Z} & \coloneq  y^ 2z
- s_ 1x^ 2y
+ s_ 2xyz
- s_ 3x^ 3
+ s_ 4x^2z
- s_ 5xz^ 2
+ s_ 6z^ 3 \\
\bar{Y} & \coloneq  x y^ 2
- s_ 1y^ 2z
+ s_ 2x^ 2y
- s_ 3xyz
+ s_ 4x^3
- s_ 5x^ 2z
+ s_ 6xz^ 2 \\
\bar{X} & \coloneq  y^ 3
+( s_ 2- s_ 1^ 2)xy^ 2
-( s_ 3- s_ 1 s_ 2)y^ 2z
+( s_ 4- s_ 1 s_ 3)x^ 2y
\\ &\qquad
-( s_ 5- s_ 1 s_ 4)xyz
+( s_ 6- s_ 1 s_ 5)x^ 3
+ s_ 1 s_ 6x^ 2z
\end{split}
\tag{\ref{eq1}a}
\end{align*}
which gives a good generating set for $L$, since $L$ is generated by
$L_1$ in this case.

In the case of $E_7$, a basis for the cubics
is given by the first three lines of equation~(\ref{eq1}b)
below.\footnote{This basis is chosen to match the one used by Bramble
\cite{[Bra]}; a simpler choice would have used $\bar{Y} - s_1^2 \bar{Z}$.}
This is completed to a good generating set for $L$ by using $\frac13$ of
the Jacobian determinant as the fourth generator.
\begin{align*}
\begin{split}
\bar{W} & \coloneq   x^3
- y z^2 \\
\bar{Z} & \coloneq  x y^ 2
- s_ 1y^ 2z
+ s_ 2x^ 2y
- s_ 3xyz
+ s_ 4x^
3
- s_ 5x^ 2z
+ s_ 6xz^ 2
-
\makebox[0pt][l]{$
s_7 z^3
$}
\\
\bar{Y} & \coloneq  4 y^ 3
+( 4 s_ 2
- 4 s_ 1^ 2
+ s_1^2)xy^ 2
-( 4 s_ 3
-
4 s_ 1 s_ 2
+ s_1^3)y^ 2z
+( 4 s_ 4
\\ &\qquad
- 4 s_ 1 s_ 3
+ s_1^2 s_2)x^ 2y
- ( 4 s_ 5
- 4 s_ 1 s_ 4
+ s_1^2 s_3)xyz
+(4  s_ 6
\\ &\qquad
- 4 s_ 1 s_ 5
+ s_1^2 s_4)
x^ 3
- (4 s_7
- 4 s_ 1 s_ 6
+ s_1^2 s_5) x^ 2z
+ (- 4 s_1 s_7
\\ &\qquad
+ s_1^2 s_6) x z^2
-  s_1^2 s_7 z^3 \\
\bar{X} & \coloneq  \tfrac13 \
\tfrac{\partial(\bar{Y},\bar{Z},\bar{W})}{\partial(x,y,z)}.
\end{split}
\tag{\ref{eq1}b}
\end{align*}

In the case of $E_8$, a basis for the cubics in $L_1$ is given by
\begin{align*}
\begin{split}
\bar{W} & \coloneq   x^3
- y z^2 \\
\bar{Z} & \coloneq  y^ 3
+( s_ 2
- s_ 1^ 2)xy^ 2
-( s_ 3
-
 s_ 1 s_ 2)y^ 2z
+( s_ 4
- s_ 1 s_ 3)x^ 2y
\\ &\qquad
- ( s_ 5
- s_ 1 s_ 4)xyz
+ ( s_ 6
- s_ 1 s_ 5)
x^ 3
- (s_7
- s_ 1 s_ 6) x^ 2z
\\ &\qquad
+ (s_8
- s_1 s_7) x z^2
+ s_1 s_8 z^3.
\end{split}
\tag{\ref{eq1}c}
\end{align*}
A basis for the sextics which determine the anti-bicanonical map is then
given by quadratic expressions in these cubics (i.e. $\Symm ^2L_1$)
together with a new sextic $F$.  To match the normalizations from
table~\ref{table45}, we assume that $F$ has weight 16 and satisfies
$F \equiv xy^5 \mod \mm $.
We then need to ensure that $F$ has multiplicity 2 along the zero-cycle
$\eta(t_1)+\ldots+\eta(t_8)$.  We do this by imposing two conditions
on $F$:
(i) $\eta^*(F) = \Psi_8(U)^2$ and
(ii) $\Psi_8(U)$ divides $\eta^*({\partial F}/{\partial x})$.
The first condition guarantees that $F$ meets $C$ with multiplicity at least
2 at each point in the zero-cycle.  (Note that $F$ and $\Psi_8(U)^2$ are monic
of the same degree).  The second guarantees that {\em one\/} of the partial
derivatives of $F$ vanishes at those points.  But now by the chain rule,
\[\eta^*(\frac{\partial F}{\partial x}) + 3 U^2
\eta^*(\frac{\partial F}{\partial y})
  =  2 \Psi_8(U)' \Psi_8(U).\]
It follows that $\Psi_8(U)$ divides $\eta^*(\partial F/\partial y)$ as well,
which implies that the other partial derivative vanishes at those points,
and therefore that $F$ has multiplicity 2 along the zero-cycle
 $\eta(t_1)+\cdots+\eta(t_n)$.

Finding $F$ explicitly is now a matter of solving the equations in the
coefficients of $F$ implied by these conditions.  The answer is not unique,
but still depends on 2 free parameters.  We used {\sc maple} and {\sc reduce}
 to solve the
equations, and made a choice for the free parameters which makes the
coefficients of $x^3y^3$ and $x^6$ both 0.
We use this polynomial $F$ as the third element $\bar{Y}$ of a good
generating set for $L$, and complete the good generating set
by using $- \frac16$ of
the Jacobian determinant as $\bar{X}$, that is,
\[ \bar{X}  \coloneq  - \frac16 \
\frac{\partial(\bar{Y},\bar{Z},\bar{W})}{\partial(x,y,z)}.\]
This good generating set is shown explicitly in Appendix 0;
in particular, the polynomial used as $\bar{Y}$ is displayed there
explicitly.

\bigskip

The second step of our computation is to compute the defining polynomial
$\bar{\Phi}_{E_n}$ of $\bar{\cal P} \subset \Proj _V(L)$ with respect to
the generating set $\bar{X}$, $\bar{Y}$, $\bar{Z}$, $\bar{W}$.
Such a polynomial was shown to exist in lemma~\ref{lem57}, and its
 general form  was given
 in equation (\ref{eq11});
we must
compute the unknown coefficients which appear in that equation.
Thanks to the remark following proposition~\ref{prop51}, we only
need the coefficients of low weight at this stage of the computation.

Let $\pi\colon \P ^2 \times V \to \bar{\cal P} \subset \Proj _V(L)$
be the rational
map determined by the  generating set
$\bar{X}$, $\bar{Y}$, $\bar{Z}$, $\bar{W}$.
Now $\bar{\Phi}_{E_n}$ vanishes on the image of $\pi$.  Thus, if we write
$\bar{\Phi}_{E_n}$ with undetermined coefficients and
compute $\pi^*(\bar{\Phi}_{E_n})$,
every coefficient in the expression for $\pi^*(\bar{\Phi}_{E_n})$ must vanish.
This produces some equations for the undetermined coefficients.

\begin{proposition} \label{prop91}
Let $\bar{\cal R}_{\pi}$ be the set of substitution rules which describes
the map $\pi$ with respect to the coordinates
$\bar{X}$, $\bar{Y}$, $\bar{Z}$, $\bar{W}$, as given in equations
(\ref{eq1}a), (\ref{eq1}b)
and Appendix 0.
Then the coefficients
of low weight in the equation $\bar{\Phi}_{E_n}$ can be described by a set of
substitution rules ${\cal R}_{\nu}$, which are given in solve-list
format by means of the substitution rules
$\bar{\cal R}_{\pi}$, the polynomial $\bar{\Phi}_{E_n}$, and the solve-list
\refstepcounter{equation}\label{eq17}
\begin{gather*}
\begin{tabular}{|c|c|c|c|c|c|c|c|c|} \hline
$x^2y^6z$    & $x^4y^5$     & $xy^6z^2$    & $x^3y^5z$     &
   $y^6z^3$       & $x^5y^4$     & $x^4y^4z$    & $x^6y^3$     &
       $x^3y^4z^2$
           \\ \hline
$\bar{\f}_1$ & $\bar{\f}_2$ & $\bar{\e}_2$ & $\bar{\f}'_3$ &
   $\bar{\f}''_3$ & $\bar{\f}_4$ & $\bar{\e}_5$ & $\bar{\f}_6$ &
       $\bar{\e}_6$
            \\ \hline
\end{tabular}
\tag{\ref{eq17}a}\\[1.5ex]
\begin{tabular}{|c|c|c|c|c|} \hline
$x^4 y^8$    & $x y^9 z^2$  & $x^5 y^7$    & $x^6 y^6$
    & $x^3 y^7 z^2$        \\ \hline
$\bar{\e}_2$ & $\bar{\f}_2$ & $\bar{\f}_4$ & $\bar{\e}_6$
    & $\bar{\f}_6$   \\ \hline
\end{tabular}
\tag{\ref{eq17}b}\\[1.5ex]
\begin{tabular}{|c|c|c|c|c|} \hline
$x^{4} y^{14}$ & $x^{5} y^{13}$ & $x^{6} y^{12}$ & $x^{7} y^{11}$
    & $x^{8} y^{10}$          \\ \hline
$\bar{\e}_{2}$ & $\bar{\f}_{4}$ & $\bar{\f}_{6}$ & $\bar{\e}_{8}$
    & $\bar{\f}_{10}$  \\ \hline
\end{tabular}
\tag{\ref{eq17}c}
\end{gather*}
\end{proposition}

\begin{pf}
We  carry out the expansion specified by the solve-list format, but
work  mod $\mm $.  The congruences in
table~\ref{table45} give substitutions for
$\bar{X}$, $\bar{Y}$, $\bar{Z}$, $\bar{W}$ which are valid mod $\mm $.
If we make those substitutions into equation
(\ref{eq11})
and collect terms of low
degree in $z$ we obtain:
\begin{align*}
\begin{split}
\pi^*(\bar{\Phi}_{E_6}) &\equiv
  \bar{\f}_1 x^2 y^6 z
+ \bar{\f}_2 x^4 y^5
+ (\bar{\e}_2 - \bar{\f}_2) x y^6 z^2
+ \bar{\f}'_3 x^3 y^5 z
+ (\bar{\f}''_3 - \bar{\f}'_3) y^6 z^3
\\ &\qquad
+ \bar{\f}_4 (x^5 y^4 - x^2 y^5 z^2)
+ \bar{\e}_5 (x^4 y^4 z - x y^5 z^3)
+ \bar{\f}_6 x^6 y^3
+ (\bar{\e}_6 - 2 \bar{\f}_6) x^3 y^4 z^2
\\ &\qquad
+ \bar{\e}_8 (x^7 y^2 - 2 x^4 y^3 z^2)
+ \bar{\e}_9 (x^6 y^2 z - 2 x^3 y^3 z^3)
+ \bar{\e}_{12} (x^9 - 3 x^6 y z^2) \\
&\qquad \mod{(\mm ,z^4)}
\end{split}
\\[1.5ex] \begin{split}
\pi^*(\bar{\Phi}_{E_7}) &\equiv
16 \bar{\e}_2 x^4 y^8
+ (16 \bar{\f}_2 - 16 \bar{\e}_2) x y^9 z^2
+ \bar{\f}_4 (4 x^5 y^7 - 4 x^2 y^8 z^2)
+ 16 \bar{\e}_6 x^6 y^6
\\ &\qquad
+ (16 \bar{\f}_6 - 32 \bar{\e}_6) x^3 y^7 z^2
+ \bar{\e}_8 (4 x^7 y^5 - 8 x^4 y^6 z^2)
+ \bar{\e}_{10} (x^8 y^4 - 2 x^5 y^5 z^2)
\\ &\qquad
+ \bar{\e}_{12} (4 x^9 y^3 - 12 x^6 y^4 z^2)
+ \bar{\e}_{14} (x^{10} y^2 - 3 x^7 y^3 z^2)
+ \bar{\e}_{18} (x^{12} - 4 x^9 y z^2)
\\ &\qquad \mod{(\mm ,z^3)} \end{split}
\\[1.5ex] \begin{split}
\pi^*(\bar{\Phi}_{E_8}) &\equiv
  \bar{\e}_2 x^4 y^{14}
+ \bar{\f}_4 x^5 y^{13}
+ \bar{\f}_6 x^6 y^{12}
+ \bar{\e}_8 x^7 y^{11}
+ \bar{\f}_{10} x^8 y^{10}
+ \bar{\e}_{12} x^9 y^9
\\ &\qquad
+ \bar{\e}_{14} y^8 x^{10}
+ \bar{\e}_{18} x^{12} y^6
+ \bar{\e}_{20} x^{13} y^5
+ \bar{\e}_{24} x^{15} y^3
+ \bar{\e}_{30} x^{18}
\\ &\qquad \mod{(\mm ,z)}.
\end{split}
\end{align*}

Since each equation which is to be solved is a coefficient of a monomial
in $x, y, z$, it is homogeneous (with respect to the
background $\C ^*$-action).  For each such equation, the $\bar{\e}_i$'s
and $\bar{\f_i}$'s involved are either the ones given above, or ones of
strictly lower weight (since they are multiplied by nontrivial functions
of the $t_i$'s).  Thus, if we proceed from lower weight to higher weight
in solving these equations, and use the leading order terms above as
a guide to the order in which equations of the same weight should be
solved, we arrive at the  solve-lists stated in the proposition.
\end{pf}

The solve-lists given in proposition~\ref{prop91} can be extended to
determine the entire defining polynomial in these coordinates.  We did this,
and
expanded the extended solve-lists using {\sc maple} and {\sc reduce}
in the cases of $E_6$ and $E_7$.
In the case of $E_6$, we obtained the expanded defining polynomial
\begin{align*}
\begin{split}
\bar{\Phi}_{E_6} &=
- \bar{X}^2 \bar{W} - \bar{X} \bar{Z}^2 + \bar{Y}^3
+ {s_{1}} \bar{Y}^2 \bar{Z}
- { {s_{2}}} \bar{X} \bar{Y} \bar{W}
+ 0 \bar{Y} \bar{Z}^2
- { {s_{3}}} \bar{X} \bar{Z} \bar{W}
\\ &\qquad
+ 0 \bar{Z}^3
- { {s_{4
}}} \bar{Y}^2 \bar{W}
+ ({ {s_{5}} { -  {s_{1}} {s_{4}}}}) \bar{Y} \bar{Z} \bar{W}
+ ({ { 2 {s_{6}}} { -  {s_{1}} {s
_{5}}}}) \bar{X} \bar{W}^2
+ 0 \bar{Z}^2 \bar{W}
\\ &\qquad
+ ({ { {s_{2}} {s_{6}} +  { {
{s_{1}}^{2}} {s_{6}}} { -  {s_{2}} {s_{1}} {s_{5}}} +  { {s_{5}} {s_{3
}}}}}) \bar{Y} \bar{W}^2
+ ({ { {s_{2}} {s_{1
}} {s_{6}}} { -  {s_{3}} {s_{6}}}}) \bar{Z} \bar{W}^2
\\ &\qquad
+ ({ { {s_{1}} {s_{6
}} {s_{2}} {s_{3}}} {- { {s_{6}}^{2}}} +  { {s_{5}} {s_{1}} {s_{6}}} {
 -  { {s_{1}}^{2}} {s_{4}} {s_{6}}} { -  { {s_{3}}^{2}} {s_{6}}}})
 \bar{W}^3.
\end{split}
\end{align*}
In the case of $E_7$, the defining polynomial which we found
agrees with the one found by Bramble \cite[p.\ 357]{[Bra]}, and we have
not reproduced it here.  The defining polynomial for $E_8$ is very large, and
we have not
attempted to write it down.

\bigskip

We can now describe a set of substitution rules ${\cal R}_{\pi}$
which describes the map $\pi$ with respect to the coordinates
$X$, $Y$, $Z$, $W$ as being  a composition
\[{\cal R}_{\pi} = {\cal R}_{\mu}^{-1} \circ \bar{\cal R}_{\pi}
\circ {\cal R}_{\psi} \circ {\cal R}_{\nu},\]
where we have used the fact that the ${\cal R}_{\mu}$ as
given in equation (\ref{eq23})
can be solved for $X$, $Y$, $Z$, $W$ as functions of the other
variables, yielding ${\cal R}_{\mu}^{-1}$.

\begin{proposition} \label{prop92}
Let ${R}_{\pi}={\cal R}_{\mu}^{-1} \circ \bar{\cal R}_{\pi}
\circ {\cal R}_{\psi} \circ {\cal R}_{\nu}$
be the set of substitution rules which describes
the map $\pi$ with respect to the coordinates
${X}$, ${Y}$, ${Z}$, ${W}$.  Then the coefficients
in  the defining polynomial ${\Phi}_{E_n}$ can be described by a set of
substitution rules  which are given in solve-list
format by means of the polynomial ${\Phi}_{E_n}$ (with undetermined
coefficients),
 the substitution rules
${R}_{\pi}$, and the solve-list
\refstepcounter{equation}\label{eq41}
\begin{gather*}
\begin{tabular}{|c|c|c|c|c|c|} \hline
    $xy^6z^2$    &
       $x^4y^4z$    &
       $x^3y^4z^2$  & $x^7y^2$     & $x^6y^2z$    & $x^6yz^2$
           \\ \hline
    $\e_2$ &
       $\e_5$ &
       $\e_6$ & $\e_8$ & $\e_9$ & $\e_{12}$
            \\ \hline
\end{tabular}
\tag{\ref{eq41}a}\\[1.5ex]
\begin{tabular}{|c|c|c|c|c|c|c|} \hline
  $x^4 y^8$
    & $x^6 y^6$    & $x^7 y^5$    & $x^8 y^4$       & $x^9 y^3$
       & $x^{10} y^2$    & $x^{12}$        \\ \hline
  $\e_2$
    & $\e_6$ & $\e_8$ & $\e_{10}$ & $\e_{12}$
       & $\e_{14}$ & $\e_{18}$ \\ \hline
\end{tabular}
\tag{\ref{eq41}b}\\[1.5ex]
\begin{tabular}{|c|c|c|c|c|c|c|c|} \hline
$x^{4} y^{14}$ & $x^{7} y^{11}$
    & $x^{9} y^{9}$   & $x^{10} y^{8}$  & $x^{12} y^{6}$
       & $x^{13} y^{5}$  & $x^{15} y^{3}$  & $x^{18}$        \\ \hline
$\e_{2}$ & $\e_{8}$
    & $\e_{12}$ & $\e_{14}$ & $\e_{18}$
       & $\e_{20}$ & $\e_{24}$ & $\e_{30}$ \\ \hline
\end{tabular}
\tag{\ref{eq41}c}
\end{gather*}

\end{proposition}

\begin{pf}
The defining polynomial in preferred versal form
with undetermined coefficients is:
\refstepcounter{equation}\label{eq35}
\begin{align*}
\begin{split}
\Phi_{E_6} &=
- X^2 W - X Z^2 + Y^3
+ \e_2 Y Z^2 +  \e_5 Y Z W + \e_6 Z^2 W + \e_8 Y W^2
\\ &\qquad
+ \e_9 Z W^2 + \e_{12} W^3
\end{split}
\tag{\ref{eq35}a}
\\[1.5ex] \begin{split}
\Phi_{E_7} &=
- X^2 - Y^3 W + 16 Y Z^3
+ \e_2 Y^2 Z W + \e_6 Y^2 W^2
+ \e_8 Y Z W^2
\\ &\qquad
+ \e_{10} Z^2 W^2 + \e_{12} Y W^3
+  \e_{14} Z W^3
+ \e_{18} W^4
\end{split}
\tag{\ref{eq35}b}
\\[1.5ex] \begin{split}
\Phi_{E_8} &=
 - X^2 + Y^3 - Z^5 W
+ \e_2 Y Z^3 W + \e_8 Y Z^2 W^2 + \e_{12} Z^3 W^3
\\ &\qquad
+ \e_{14} Y Z W^3
+ \e_{18} Z^2 W^4
+  \e_{20} Y W^4
+ \e_{24} Z W^5
+ \e_{30} W^6.
\end{split}
\tag{\ref{eq35}c}
\end{align*}

As before, we  carry out the procedure specified by the solve-list format,
working  mod $\mm $.  The congruences in
table~\ref{table45} give substitutions for
${X}$, ${Y}$, ${Z}$, ${W}$ which are valid mod $\mm $.
If we make those substitutions into the equation for $\Phi_{E_n}$ above
and collect terms of low
degree in $z$ we obtain:
\begin{align*}
\begin{split}
\pi^*({\Phi}_{E_6}) &\equiv
  {\e}_2  x y^6 z^2
+ {\e}_5 x^4 y^4 z
+ {\e}_6  x^3 y^4 z^2
+ {\e}_8 (x^7 y^2 - 2 x^4 y^3 z^2)
\\ &\qquad
+ {\e}_9 x^6 y^2 z
+ {\e}_{12} (x^9 - 3 x^6 y z^2)
\\ &\qquad \mod{(\mm ,z^3)}
\end{split}
\\[1.5ex] \begin{split}
\pi^*({\Phi}_{E_7}) &\equiv
16 {\e}_2 x^4 y^8
+ 16 {\e}_6 x^6 y^6
+ {\e}_8 4 x^7 y^5
+ {\e}_{10} x^8 y^4
\\ &\qquad
+ {\e}_{12} 4 x^9 y^3
+ {\e}_{14} x^{10} y^2
+ {\e}_{18} x^{12}
\\ &\qquad \mod{(\mm ,z)} \end{split}
\\[1.5ex] \begin{split}
\pi^*({\Phi}_{E_8}) &\equiv
  {\e}_2 x^4 y^{14}
+ {\e}_8 x^7 y^{11}
+ {\e}_{12} x^9 y^9
+ {\e}_{14} y^8 x^{10}
+ {\e}_{18} x^{12} y^6
\\ &\qquad
+ {\e}_{20} x^{13} y^5
+ {\e}_{24} x^{15} y^3
+ {\e}_{30} x^{18}
\\ &\qquad \mod{(\mm ,z)}.
\end{split}
\end{align*}

The same argument used in the proof of proposition~\ref{prop91} shows that
if we proceed from lower weight to higher weight
in solving these equations,
 we arrive at the  solve-lists stated in the proposition.
\end{pf}

It is important to notice that the solve-lists for $E_7$ and $E_8$
given in  proposition~\ref{prop92} do
not involve $z$; it is therefore possible (and desirable) to set $z=0$
at the beginning of any computation involving these solve-lists.

We have expanded these solve-lists using {\sc maple} and {\sc reduce}
in the cases of $E_6$ and $E_7$;
the results of this computation
 are displayed in Appendices 1 and 2, respectively.
(The results for $E_7$ can be found, with some errors, in
Bramble \cite{[Bra]}; Appendix 2 gives the corrected results.)
The results for $E_8$ are too large to contemplate writing down.
We stress however that even in the cases of $E_6$ and $E_7$,
further calculations with these formulas are best done by leaving
them in solve-list format as long as possible, and only expanding
the solve-lists at the very end, after restricting to a suitable subspace.

\section{Restricted polynomials and
the main computation.}

In this section, we will justify the congruences stated in
table~\ref{table-key1}.
These congruences all take the form
$\e_i \equiv c \cdot (\widetilde{\f}_N)^d \mod I$.
We will in fact compute $\e_i$ modulo another
 ideal $J \subset I$ which is generated by a certain subset of the set of
standard
coordinate functions, and note that our desired congruence is an
immediate consequence of this computation.
In order to describe the ideal $J$ and our computational method,
 we must first describe in parametric form some subspaces
of the deformation spaces $\Def (S)$.

Let $R$ be a root system of type $S$, let $V_R$ be the complex root
space, and let $t_1,\ldots,t_n$ be the distinguished functionals,
with $s_1,\ldots,s_n$ their elementary symmetric functions.
Suppose that we are given a vector space $W$ and a monic
 polynomial $r_S(U)$ of degree $n$ in $U$, whose coefficients lie in
the ring
$\C [W]$ of polynomial functions on $W$.
This determines a map $\psi\colon W \to V_R/\Sym_n$ by means of the action on
polynomials
$\psi^*\colon\C [V_R]^{\sym_n} \to \C [W]$
defined by
\[\psi^*(s_i) = \text{the coefficient of } U^{n-i} \text{ in } r_S(U).\]
In particular, the pullback of the distinguished polynomial is
$\psi^*(f_S(U;t)) = r_S(U)$.  We call $r_S(U)$ the {\em restricted
polynomial\/} associated to $\psi$.

We wish to describe a particular case of this construction
 in which the image of $\psi$ is defined
by an ideal $J_S$ which is generated by a subset of the standard
coordinate functions on $\Def (S)$.  When we have done so, we will call the
generators of $J_S$ the {\em vanishing coordinates}, and refer to
coordinates on $W$ as {\em parameters}.  We want $J_S$ to be as large
as possible, yet not to contain the ``constant term".
Equivalently, $W$ parametrizes a subspace of $\Def (S)$ on which many
of the standard coordinate functions (but not the ``constant term") vanish.

If  $S=A_{n-1}$, let $J_{A_{n-1}}$ be the ideal
generated by  all the standard
coordinate functions on $\Def (A_{n-1})$ other than the ``constant term".
Since the standard coordinate functions are the elementary symmetric
functions $s_i$ themselves,
it is easy to construct a restricted polynomial
for the ideal $J_{A_{n-1}}$.  There are two natural choices:
we use either $r_{A_{n-1}}(U) = U^n + \l_n$, or
 $r_{A_{n-1}}(U) = U^n - \l_1^n$.  (We will
use the second form when we need an explicit
root $\l_1$ of the restricted polynomial.)

If $S=D_n$ and $n$ is even, we again let $J_{D_n}$ be the ideal
generated by  all the standard
coordinate functions on $\Def (D_{n})$ other than the ``constant term".
The construction of a restricted polynomial in this case
is based on a special
factorization property:  if we define
\[F(U) = U^n - \l_{n-1} U,\]
 and
\[G(-U^2) = F(U) \cdot F(-U)\]
then
\[G(Z) = Z^n + \l_{n-1}^2 Z.\]
Thus, if we let $r_{D_{n}}(U) = F(U)$, then
the pullbacks via
$\psi^*$ of the standard coordinate functions will be the coefficients
of $G(Z)$, together with one coefficient of $F(U)$.  It follows that
$\psi^*$ of all standard coordinate functions  except for the
``constant term" vanish.

In the remaining cases, we can find an ideal $J_S$ and a
 restricted polynomial as follows.
Begin with parameters $\l_1,\ldots,\l_n$,  the initial restricted
 polynomial
$U^n + \sum \l_i U^{n-i}$, and the map given by $\psi^*(s_i) = \l_i$.
The pullbacks of the standard coordinate functions $\f_j$ via $\psi^*$
 can be computed
as functions of the $\l_i$.  If the weight $j$ of a standard
coordinate function $\f_j$ is at most $n$,
then $\l_j$ appears in the formula for that
coordinate function with a nonzero constant coefficient.
(This must be checked case by case.)  It follows that
if we set all such standard coordinate functions
 equal to zero, we get a triangular
system of linear equations in a subset of the set $\{\l_i\}$ of parameters.
These equations can be derived from equations (\ref{eqD}),
 (\ref{eqE4}), and
(\ref{eqE5}) in the   cases of $D_7$,
 $E_4$, and $E_5$ respectively, and from Appendices~1 and 2 in the
case of $E_6$ and $E_7$.
We used {\sc maple} and {\sc reduce} to solve the equations in those
5 cases; the resulting
 restricted polynomials are shown in table~\ref{tableABC}.

If $\f_N$ is the ``constant term" (or for that matter any of the
standard coordinate functions on $\Def (S)$) it is a straightforward
matter to compute $\psi^*(\f_N)$, based on the description of the
mapping $\psi$ which is given by the coefficients of the restricted
polynomial as shown in table~\ref{tableABC}.  For this purpose, one
again uses equations (\ref{eqD}),  (\ref{eqE4}), and
(\ref{eqE5}) in the cases of $D_7$, $E_4$, and $E_5$, respectively.
We have carried out this computation using {\sc maple} and
{\sc reduce}, and displayed the answers in the last column of
table~\ref{tableABC}.

In the case of  $E_6$ and $E_7$, it is more efficient to
perform this computation directly from the solve-list description
of the standard coordinate functions on $\Def (S)$
which was given in section 9.
(That is, we use the explicit description of $\psi^*(s_i)$ to substitute for
$s_i$ in the ingredients of the solve-list, and then solve the
resulting equations.)
When this is done, formulas $e_{12}(\l_1, \l_3, \l_4)$
and $e_{18}(\l_1, \l_3, \l_4,\l_5, \l_7)$ are obtained which express
the pullbacks of the respective ``constant terms" $\e_{12}$ and $\e_{18}$
in terms of the parameters.  We calculated these formulas using
{\sc maple} and
{\sc reduce}, but as they are a bit long, we have not displayed them
in the table.

{\renewcommand{\arraystretch}{1.2}

\begin{table}[p]
\begin{center}
\begin{tabular}{|c|c|l|c|} \hline
 & & & \\
$S$  &  Vanishing  &
                              \multicolumn{1}{c|}{Restricted Polynomial}
& ``Constant Term'' \\
  &    Coordinates  & \multicolumn{1}{c|}{$r_S(U)$} & $\psi^*(\f_N)$ \\
 & & & \\ \hline
 & & & \\
$A_{n-1}$
 & $\a_j$, $_{2 \le j \le n-1}$
& \multicolumn{1}{c|}{$U^n+\l_n$\ \ \ \text{or}\ \ \ $U^n - \l_1^n$}
& $\l_n$\ \ \ or\ \ \ $- \l_1^n$ \\
 & &  & \\
$D_n$ &
     $\c_n$, $\d_{2j}$,  & \multicolumn{1}{c|}{$U^n - \l_{n-1} U$}
& $\l_{n-1}^2$ \\
$_{n\ \text{even}}$ & $_{1 \le j \le n-2}$  & & \\
 & & & \\
$D_7$
& $\d_2, \d_4, \d_6, \c_7$ &
$U^7 + \l_1 U^6 + \frac12 \l_1^2 U^5 + \l_3 U^4 +
(\l_1 \l_3
$ & $(\l_1 \l_5 - \frac12 \l_3 \l_1^3 + \frac{1}{16} \l_1^6
+ \frac12 \l_3^2)^2$ \\
 & & $ -\ \frac18 \l_1^4) U^3 + \l_5 U^2
+
(\l_1 \l_5 - \frac12 \l_3 \l_1^3 $
&  \\
 & & $ +\ \frac{1}{16} \l_1^6
+ \frac12 \l_3^2) U$
&  \\
 & & & \\
$E_4$ & $\e_2$, $\e_3$, $\e_4$
& $U^4 + \l_1U^3 + \frac35 \l_1^2U^2 + \frac{1}{25}\l_1^3U + \frac{11}{125}
\l_1^4$  & $\frac{243}{3125}\l_1^5$ \\
 & & & \\
$E_5$ & $\e_2, \e_4, \e_5$ &
$ U^5 + \l_1U^4 + \frac58 \l_1^2U^3 + \l_3U^2 + (
\frac{15}{128} \l_1^4  $
& ${ { \frac{ 2601 }{16384}}{ \l_{1}^{8}} +  { \frac{ 9 }{4}}{ \l
_{3}^{2}} { \l_{1}^{2}} { -  \frac{ 153 }{128}}{ \l_{1}^{5}} \l_{3
}}$ \\
 & & $  -\ \frac12 \l_1\l_3)U + ( \frac{27}{256}\l_1^5 - \frac12 \l_1^2\l_3)$
&  \\
 & & & \\
$E_6$ & $\e_2, \e_5, \e_6$ &
$U^6 + \l_1 U^5 + \frac23 \l_1^2 U^4 + \l_3 U^3 + \l_4 U^2
$ & $e_{12}(\l_1, \l_3, \l_4)$ \\
 & & $
+\ ( \frac13 \l_1 \l_4 - \frac13 \l_3 \l_1^2 + \frac{2}{27} \l_1^5) U
$ & \\
 & & $
+\ ( \frac{5}{18} \l_1^2 \l_4 - \frac19 \l_3 \l_1^3 +
\frac{11}{486} \l_1^6 - \frac18 \l_3^2 )$
&  \\
 & & & \\
$E_7$ & $\e_2, \e_6$ &
$U^7 + \l_1 U^6 + \frac34 \l_1^2 U^5 + \l_3 U^4 + \l_4 U^3
$ & $e_{18}(\l_1, \l_3, \l_4,\l_5, \l_7)$ \\
  &  & $
+\ \l_5 U^2 + ( - \frac18 \l_3^2 + \frac{3}{64} \l_1^6 - \frac{3}{16} \l_3
\l_1^3
$ & \\
  & & $
-\ \frac14 \l_1 \l_5 + \frac38 \l_1^2 \l_4 ) U
+ \l_7 $ & \\
 & & & \\ \hline
\end{tabular}
\end{center}

\medskip

\caption{}
\label{tableABC}
\end{table}

}

We are now ready to explain how the congruences in table~\ref{table-key1} are
derived.
Let $(E_n,v_k)$ be one of the pairs considered in table~\ref{table-key1}.
The ``constant term" of highest weight $\widetilde{\f}_N$ is induced from the
projection onto a subspace $\widetilde{V}/\widetilde{\W}$ corresponding
to an irreducible subsystem $\widetilde{R}$ of the root system
$R_{E_n}$.  (In all cases except $(E_7,v_2)$ and $(E_8,v_2)$, $R$
is the ``left part" $R'$; in those two cases,
$R$ is the ``right part" $R''$.)
As in section 8, let $I$ be the ideal
in $\C [V]^{\w}$ which is generated by all the standard coordinate
functions on $\PRes (E_n,v_k)$ other than $\widetilde{\f}_N$.
Let $J$ be the ideal in $\C [V]^{\w}$ which is generated by
$J_{\widetilde{R}}$, $\m_1$, and all the standard coordinate functions
on  $\PRes (E_n,v_k)$ which come from $R - \{v_k\} - \widetilde{R}$.
Since the ``constant term" of highest weight is associated to $\widetilde{R}$
but does not belong to $J_{\widetilde{R}}$, we have $J \subset I$.
We extend the map
$\widetilde{\psi}\colon\widetilde{W} \to \widetilde{V}/\Sym_{\widetilde{n}}$
to a map $\psi\colon \widetilde{W} \to V/(\Sym_{n'} \times \Sym_{n''})$
by simply composing it with the natural inclusion
$\widetilde{V}/\Sym_{\widetilde{n}} \subset
 V/(\Sym_{n'} \times \Sym_{n''})$.

We need to compute $\e_i$ modulo $J$.  Since $J$ vanishes on the
image of $\psi$, it suffices to compute $\psi^*(\e_i)$ in terms of
the pullbacks of the standard coordinate functions on $\PRes (E_n,v_k)$
via $\psi^*$.
The first step is to use proposition~\ref{prop71} and table~\ref{tableAA},
which relates the distinguished polynomial for $R$ to those for $R'$
and $R''$.  Now $\psi^*(f_{S'}(U;t'))$ and $\psi^*(f_{S''}(U;t''))$
can be computed immediately:
one of them
is the restricted polynomial for $\widetilde{R}$, and the other one
is just a power of $U$ (since all of the corresponding standard coordinate
functions vanish when
pulled back via $\psi$).  Table~\ref{tableAA} can then be used to compute
$\psi^*(f_S(U;t))$; we carry out this computation below.

We first consider the case $(E_n,v_0)$, in which the complementary root
system has type $A_{n-1}$.  We have
$\psi^*(\m_1)=0$ and
\[\psi^*(f_{A_{n-1}}(U;t'))=r_{A_{n-1}}(U)=U^n+\l_n.\]
(We use the first form of the restricted polynomial since we do not
need a root.)  Thus, by table~\ref{tableAA},
$\psi^*(f_{E_n}(U;t))=U^n+\l_n$.

We next consider the case $(E_n,v_1)$, in which the complementary root
system has type $D_{n-1}$.  Since
$\psi^*(f_{D_{n-1}}(U;t'))=r_{D_{n-1}}(U)$,
the coefficient $\psi^*(\rho_1)$ of $U^{n-2}$ in this polynomial is
0 when $n-1$ is even, and $\l_1$ when $n-1=7$.  We also have
$\psi^*(\m_1)=0$.  Thus by table~\ref{tableAA},
 in the case $(E_7,v_1)$
we get
\begin{align*}
 \psi^*(f_{E_7}(U;t))
&=   (-1)^7 \cdot (-U) \cdot \widetilde{\psi}^*(f_{D_6}(-U;t')) \\
 &=  U^7 + \l_5 U^2
\end{align*}
while in the case $(E_8,v_1)$ we get
\begin{align*}
 \psi^*(f_{E_8}(U;t))
&=  (-1)^8 \cdot (-U+\frac13\l_1)
\cdot \widetilde{\psi}^*(f_{D_7}(-U-\frac16\l_1;t')) \\
 &=  (-U+\frac13\l_1)
\cdot r_{D_7}(-U-\frac16\l_1).
\end{align*}

We next consider the case $(E_n,v_2)$, in which the complementary root
system has components of type $A_1$ and $A_{n-2}$.  In this case,
$\psi^*(f_{A_{n-2}}(U;t''))=r_{A_{n-2}}(U)$,
which we write this time in the form $U^{n-1} - \l_1^{n-1}$ so that
$\psi^*(\s_1)=\l_1$ is a root of this polynomial.
Now
\[\frac{r_{A_{n-2}}(U)}{U-\l_1} = \sum_{i=0}^{n-2}\l_1^i\ U^{n-2-i}.\]
Moreover, $\psi^*(\m_1)=0$, and
$\psi^*(f_{A_1}(U+\frac23\s_1;t'))=(U+\frac23\l_1)^2$.
Table~\ref{tableAA} then implies
\begin{align*}
 \psi^*(f_{E_n}(U;t))
&=   (U+\frac23\l_1)^2 \cdot
\sum_{i=0}^{n-2}\l_1^i\ (U-\frac13\l_1)^{n-2-i}
\end{align*}

Finally, we consider the case $(E_n,v_k)$ with $k \ge 4$,
in which the complementary root
system has components of type $E_k$ and $A_{n-k-1}$.  In this case,
we have
$\psi^*(f_{E_k}(U;t'))=r_{E_k}(U)=U^k+\l_1U^{k-1}+\cdots$,
which implies that the coefficient $\psi^*(\tau_1)$ of
$U^{k-1}$ in this polynomial
is $\l_1$.
Since $\psi^*(\m_1)=0$ and $\psi^*(f_{A_{n-k-1}}(U;t''))=U^{n-k}$,
table~\ref{tableAA}  implies that
\begin{align*}
 \psi^*(f_{E_n}(U;t))
&=
(U-\frac{1}{9-k}\l_1)^{n-k} \cdot   r_{E_k}(U) .
\end{align*}

The second step in the computation of $\e_i$ modulo $J$ is to
 use the coefficients of
the pulled-back distinguished polynomial $\psi^*(f_{E_n}(U;t))$
as ingredients for the solve-lists in section 9, and obtain (using
{\sc maple} and {\sc reduce}) a formula for $\psi^*(\e_i)$ in terms
of the parameters.  Now we have already computed
$\psi^*(\widetilde{\f}_N)$ in terms of the parameters, as indicated
in table~\ref{tableABC}.  So we simply need to compare the formulas
for $\psi^*(\e_i)$ and $\psi^*(\widetilde{\f}_N)^d$.

We illustrate this comparison
with an example which can be carried out by hand.
The case we consider is $(E_7,v_1)$, in which the complementary root
system has type $D_6$.  We have
 $\psi^*(f_{E_7}(U;t)) = U^7 + \l_5 U^2$, which implies that
 $\psi^*(\e_{10})$ is computed by setting $s_5=\l_5$ and all
other $s_j=0$ in the formula for $\e_{10}$.  The only term that then
remains is the term coming from $s_5^2$ in the original formula for
$\e_{10}$.  Now inspection of Appendix 2 shows that the coefficient of
$s_5^2$ in the formula for $16\e_{10}$ is 256.  Thus,
\[\psi^*(\e_{10}) = 16 \l_5^2 = 16 \psi^*({s_5})^2\]
which implies that
\[\e_{10} \equiv 16 s_5^2 \mod J,\]
as required.

To return to the general argument:
in all cases from table~\ref{table-key1} except $(E_8,v_1)$ and
$(E_8,v_7)$, the only monomial in the standard coordinate functions on
$\PRes (E_n,v_k)$ which
has weight $i$ and which does not pull back to zero under $\psi^*$
is $(\widetilde{\f}_N)^d$.  Thus, in those cases it follows that
$\psi^*(\e_i)/\psi^*(\widetilde{\f}_N)^d$
is a constant.  We calculated these constants
 using {\sc maple} and {\sc reduce},
obtaining the values indicated in table~\ref{table-key1}.
This verifies the congruences
$\e_i \equiv c \cdot (\widetilde{\f}_N)^d \mod J$,
which suffices since $J \subset I$ in each case.

In the two remaining cases $(E_8,v_1)$ and
$(E_8,v_7)$, the congruences which hold modulo $J$ are
\begin{align}
\label{eqlast1}
\e_{24} &\equiv
0 \cdot (\d'_8)^3 + -\tfrac{1}{16}\cdot (\d'_{12})^2
\mod J, \quad \text{and}\\
\label{eqlast2}
\e_{18} &\equiv
-\tfrac{1}{3072} \cdot (\e'_8 \e'_{10})
+ \tfrac{1}{64}\cdot (\e'_{18})
\mod J,
\end{align}
respectively, and these imply the desired congruences modulo $I$.
  To verify these congruences, we also need to calculate
\begin{gather*}
\psi^*(\d'_8)=
{ { { {\l_{1}}^{3}} {\l_{5}}} { -  \tfrac{ 3 }{4}
}{ {\l_{1}}^{5}} {\l_{3}} +  { \tfrac{ 3 }{2}}{ {\l_{1}}^{2}} { {\l_{3}}^{2}}
+  { \tfrac{ 5 }{64}}{
{\l_{1}}^{8}} { -  2 {\l_{3}} {\l_{5}}}},
\\
\psi^*(\e'_8)=e_8(\l_1, \l_3, \l_4,\l_5, \l_7),\\
\psi^*(\e'_{10})=e_{10}(\l_1, \l_3, \l_4,\l_5, \l_7).
\end{gather*}
The first of these formulas is obtained from equation (\ref{eqD}),
while the second and third lines refer to formulas which we have calculated
explicitly with {\sc maple} and {\sc reduce}
using the solve-list method, but do not display here.
(Notice that the calculations of $\d'_{12}$ and $\e'_{18}$
are indicated in table~\ref{tableABC}.)
Now the coefficients in equations (\ref{eqlast1}) and (\ref{eqlast2})
can be calculated with the method of undetermined coefficients.
That is, there will be some relation of the form
\begin{align*}
\psi^*(\e_{24}) &=  c_1 \cdot \psi^*(\d'_8)^3
+ c_2 \cdot \psi^*(\d'_{12})^2
\mod J, \quad \text{or} \\
\psi^*(\e_{18}) &=  c_1 \cdot \psi^*(\e'_8) \cdot \psi^*(\e'_{10})
+ c_2 \cdot \psi^*(\e'_{18})
\mod J,
\end{align*}
respectively.  Substituting the calculated values of $\psi^*$
allowed us to solve (using {\sc maple} and {\sc reduce}) for
the undetermined coefficients $c_1$, $c_2$.

This completes the verification of table~\ref{table-key1}, and the
proof of the main theorem.

\bigskip

We would like to offer two pieces of advice to the ambitious reader who
wishes to duplicate our symbolic calculations.  First, it is essential
when computing with solve-lists to keep them unexpanded as long
as possible.  Even when a solve-list must be expanded, it may be
that all relevant information can be extracted by only {\em partially\/}
expanding the solve list, solving for a proper subset of the variables.

Second, the absence of $z$ from the monomials in the solve-lists
(\ref{eq17}c), (\ref{eq41}b), and (\ref{eq41}c)
means that $z$ can be set equal to $0$ before the expansion of
these solve-lists begins.  This cuts down the size of the
computation tremendously.

For the less ambitious reader, {\sc maple} source files for all
calculations described in the paper are available upon request
(directed to the second author).

\newpage

\section*{Appendix 0. A good generating set in the case of $E_8$.}
\begin{align*}
\bar{W} &=  x^3 - y z^2 \\
\bar{Z} &= y^ 3+( s_ 2- s_ 1^ 2)xy^ 2-( s_ 3-
 s_ 1 s_ 2)y^ 2z+( s_ 4- s_ 1 s_ 3)x^ 2y
- ( s_ 5- s_ 1 s_ 4)xyz \\
   &\quad
+ ( s_ 6- s_ 1 s_ 5)
x^ 3 - (s_7 - s_ 1 s_ 6) x^ 2z + (s_8 - s_1 s_7) x z^2 + s_1 s_8 z^3 \\
\bar{Y} &=  x  y^5
 -  2 s_1  y^5 z
+   \left(    s_1^2+  2 s_2 \right)   y^4  x^2
+   \left(   - 2 s_3  -  2 s_1 s_2 \right)  z  y^4 x
+   \left(    s_2^2 +   2 s_1 s_3
+  2 s_4 \right)   z^2  y^4 \\
 &\quad
+   \left(    -  2 s_1 s_4 - 2 s_5  -  2 s_3 s_2 \right)  z  y^3  x^2
+   \left(   s_6  -  2 s_3  s_1^3 +   3  s_1^2 s_4
-  s_2^3 -  s_1^6 +
       3 s_2  s_1^4 \right)   y^2  x^4 \\
 &\quad
+   \left(    2 s_1 s_5+ s_6 +   2 s_3  s_1^3+  s_2^3+  s_3^2
     -  3  s_1^2 s_4 +   2 s_2 s_4
+  s_1^6  -  3 s_2  s_1^4 \right)   z^2  y^3 x \\
&\quad
+   \left(    2 s_1 s_2 s_4
-  2 s_3 s_4 +   s_5  s_1^2
     -  2  s_1^2 s_2 s_3  -   s_1^3 s_4  -  s_2  s_1^5
- s_7
     +    s_2^3 s_1  -  s_1 s_6 +   2  s_2^2  s_1^3 \right. \\
 &\quad \left.
+   s_3  s_1^4
     -  s_3  s_2^2 \right)   z^3  y^3
+   \left(    -  2 s_5 s_2 - s_7 +   2  s_1^2 s_2 s_3
+   s_3  s_2^2
     -  s_1 s_6  -   s_2^3 s_1  -  2  s_2^2  s_1^3  \right. \\
 &\quad \left.
+    s_1^3 s_4
     -  2 s_1 s_2 s_4 +   s_2  s_1^5  -  s_3  s_1^4  -  s_5  s_1^2 \right)
     z  y^2  x^3
+   \left(    -  s_4  s_2^2  -  s_1 s_7 +   2 s_3 s_1 s_4 -  s_4^2
\right. \\
 &\quad \left.
     -   s_3^2  s_1^2 +   2 s_3  s_1^3 s_2 - s_8  -  s_5  s_1^3
     +   s_6 s_2
-  s_4 s_2  s_1^2 +   s_6  s_1^2  -  s_3  s_1^5
     +   s_4  s_1^4 +   s_3 s_1  s_2^2 \right)  y  x^5 \\
&\quad
+   \left(    3 s_1 s_7 +   s_4  s_2^2+  3 s_8 +   2 s_3 s_5
+    s_3^2  s_1^2
     +   s_6 s_2  -  2 s_3  s_1^3 s_2+  2  s_4^2 +   s_3  s_1^5
     -  s_4  s_1^4  \right. \\
 &\quad \left.
-  s_3 s_1  s_2^2 +   s_4 s_2  s_1^2 +   s_5  s_1^3
-  s_6  s_1^2  -  2 s_3 s_1 s_4 \right)   z^2  y^2  x^2
+   \left(    s_5 s_2  s_1^2  -  s_4 s_1  s_2^2 +    s_1^2 s_4 s_3 \right. \\
&\quad \left.
-  2 s_1  s_4^2 +   s_4  s_1^5
+   s_6  s_1^3  -  s_6 s_3
     +   s_1 s_8  -   s_1^2 s_7 +   s_5  s_2^2  -  s_2 s_7
     -  2 s_4 s_2  s_1^3  -  s_5  s_1^4 \right)  z y  x^4  \\
&\quad
+   \left(    s_4 s_1  s_2^2  -  s_6  s_1^3  -  s_5 s_2  s_1^2  -  3 s_1 s_8
     -  s_4  s_1^5 +   2 s_1  s_4^2 +    s_1^2 s_7  -  2 s_5 s_4
     -  s_2 s_7  \right. \\
 &\quad \left.
+   2 s_4 s_2  s_1^3
-  s_6 s_3  -  s_5  s_2^2
     -   s_1^2 s_4 s_3 +   s_5  s_1^4 \right)   z^3  y^2 x
+   \left(    s_2 s_8  -  2 s_5  s_1^3 s_2  -  s_1 s_5  s_2^2  \right. \\
 &\quad \left.
+   s_5  s_1^5
+   s_7 s_3 +   s_2 s_6  s_1^2  -  s_6  s_1^4 +    s_1^3 s_7
     +    s_2^2 s_6 +   s_6 s_4 +   s_3 s_5  s_1^2  -  2 s_1 s_4 s_5 \right. \\
 &\quad \left.
-  s_8  s_1^2 \right)   z^4  y^2
+   \left(    s_7 s_3  -   s_2^2 s_6 +   s_6 s_4  -   s_1^3 s_7  -  s_5  s_1^5
     +  s_5^2 +   s_2 s_8 +   s_6  s_1^4
-  s_3 s_5  s_1^2 \right. \\
 &\quad \left.
+   2 s_1 s_4 s_5 +   2 s_5  s_1^3 s_2 +   s_1 s_5  s_2^2
     -  s_2 s_6  s_1^2 +   s_8  s_1^2 \right)   z^2 y  x^3
+   \left(    s_6  s_1^2 s_3
+   s_4 s_7  -  2 s_6  s_1^3 s_2  \right. \\
 &\quad \left.
-  s_7  s_1^4
     -  s_5 s_6 +   s_6  s_1^5  -  s_6 s_1  s_2^2 +   s_8  s_1^3
     +   s_7  s_2^2  -  2 s_6 s_1 s_4
+   s_3 s_8
     +   s_2 s_7  s_1^2 \right)  z  x^5 \\
 &\quad
+   \left(    -  s_5 s_6  -  3 s_4 s_7 +   s_6 s_1  s_2^2 +   2 s_6  s_1^3 s_2
     -  s_8  s_1^3  -  3 s_3 s_8
-  s_6  s_1^2 s_3  -  s_6  s_1^5
     +   2 s_6 s_1 s_4 \right. \\
 &\quad \left.
+   s_7  s_1^4  -  s_2 s_7  s_1^2
-  s_7  s_2^2 \right)   z^3 y  x^2
+   \left(    2 s_7 s_2  s_1^3
-   s_1^2 s_7 s_3 +   s_8  s_1^4+  s_6^2
     +   s_5 s_7  -  s_8 s_2  s_1^2 \right. \\
 &\quad \left.
+    s_2^2 s_7 s_1  -  s_4 s_8
     +   2 s_7 s_1 s_4  -  s_7  s_1^5
-  s_8  s_2^2 \right)   z^2  x^4
+   \left(     s_1^2 s_7 s_3  -   s_2^2 s_7 s_1 +   s_8 s_2  s_1^2 \right. \\
 &\quad \left.
+   3 s_4 s_8 +   s_5 s_7  -  2 s_7 s_1 s_4
-  2 s_7 s_2  s_1^3
     -  s_8  s_1^4 +   s_8  s_2^2 +   s_7  s_1^5 \right)   z^4 y x
+  \left(    s_8 s_1  s_2^2 \right. \\
 &\quad \left.
+   2 s_1 s_8 s_4 +   2 s_8  s_1^3 s_2
     -  s_8  s_1^5
-  s_3 s_8  s_1^2  -  s_5 s_8 \right)   z^5 y
+   \left(    s_8  s_1^5  -  2 s_6 s_7  -  s_8 s_1  s_2^2  \right. \\
 &\quad \left.
-  2 s_1 s_8 s_4
     -  2 s_8  s_1^3 s_2  -  s_5 s_8
+   s_3 s_8  s_1^2 \right)   z^3  x^3
+   \left(    s_7^2 +   2 s_6 s_8 \right)   z^4  x^2
-  2 s_7 s_8 x  z^5
+   s_8^2  z^6 \\
\bar{X} &= - \frac16 \
\frac{\partial(\bar{Y},\bar{Z},\bar{W})}{\partial(x,y,z)}
\end{align*}

\section*{Appendix 1. Standard coordinates for $E_6$.}
\begin{align*}
 6\,{\e_{2}} &= -2\,{s_{1}}^2 + 3\,s_{2}\\
 81\,{\e_{5}} &= 4\,{s_{1}}^5 - 15\,{s_{1}}^3\,s_{2} + 27\,{s_{1}}^2\,
s_{3} - 27\,s_{1}\,s_{4} + 81\,s_{5}\\
 1944\,{\e_{6}} &= -16\,{s_{1}}^6 + 72\,{s_{1}}^4\,s_{2} -
216\,{s_{1}}^3\,s_{3} + 27\,{s_{1}}^2\,{s_{2}}^2 + 216\,{s_{1}}^2\,s_{
4} \\&\qquad + 162\,s_{1}\,s_{2}\,s_{3} + 324\,s_{1}\,s_{5} - 81\,{s_{
2}}^3 + 324\,s_{2}\,s_{4} - 243\,{s_{3}}^2 - 1944\,s_{6} \\
 34992\,{\e_{8}} &= -64\,{s_{1}}^8 + 384\,{s_{1}}^6\,s_{2} -
864\,{s_{1}}^5\,s_{3} - 324\,{s_{1}}^4\,{s_{2}}^2 + 864\,{s_{1}}^4\,s
_{4} \\&\qquad + 1944\,{s_{1}}^3\,s_{2}\,s_{3} - 2592\,{s_{1}}^3\,s_{5
} - 486\,{s_{1}}^2\,{s_{2}}^3 - 1944\,{s_{1}}^2\,s_{2}\,s_{4} - 2916\,
{s_{1}}^2\,{s_{3}}^2 \\&\qquad + 34992\,{s_{1}}^2\,s_{6} + 2916\,s_{1}
\,{s_{2}}^2\,s_{3} - 11664\,s_{1}\,s_{2}\,s_{5} + 5832\,s_{1}\,s_{3}\,
s_{4} - 729\,{s_{2}}^4 \\&\qquad + 5832\,{s_{2}}^2\,s_{4} - 4374\,s_{2
}\,{s_{3}}^2 + 17496\,s_{3}\,s_{5} - 11664\,{s_{4}}^2 \\
 78732\,{\e_{9}} &= 64\,{s_{1}}^9 - 432\,{s_{1}}^7\,s_{2} +
1296\,{s_{1}}^6\,s_{3} + 324\,{s_{1}}^5\,{s_{2}}^2 - 1296\,{s_{1}}^5\,
s_{4} \\&\qquad - 3888\,{s_{1}}^4\,s_{2}\,s_{3} - 1944\,{s_{1}}^4\,s_{
5} + 1215\,{s_{1}}^3\,{s_{2}}^3 + 972\,{s_{1}}^3\,s_{2}\,s_{4} + 5832
\,{s_{1}}^3\,{s_{3}}^2 \\&\qquad - 14580\,{s_{1}}^3\,s_{6} - 2187\,{s
_{1}}^2\,{s_{2}}^2\,s_{3} + 17496\,{s_{1}}^2\,s_{2}\,s_{5} - 8748\,{s
_{1}}^2\,s_{3}\,s_{4} \\&\qquad + 2187\,s_{1}\,{s_{2}}^2\,s_{4} +
52488\,s_{1}\,s_{2}\,s_{6} - 8748\,s_{1}\,{s_{4}}^2 - 6561\,{s_{2}}^2
\,s_{5} - 78732\,s_{3}\,s_{6} \\&\qquad + 26244\,s_{4}\,s_{5}
\\
 11337408\,{\e_{12}} &= -256\,{s_{1}}^{12} + 2304\,{s_{1}}^{
10}\,s_{2} - 6912\,{s_{1}}^9\,s_{3} - 4320\,{s_{1}}^8\,{s_{2}}^2
 \\&\qquad + 6912\,{s_{1}}^8\,s_{4} + 36288\,{s_{1}}^7\,s_{2}\,s_{3}
 + 10368\,{s_{1}}^7\,s_{5} - 6480\,{s_{1}}^6\,{s_{2}}^3 \\&\qquad -
20736\,{s_{1}}^6\,s_{2}\,s_{4} - 54432\,{s_{1}}^6\,{s_{3}}^2 + 217728
\,{s_{1}}^6\,s_{6} - 11664\,{s_{1}}^5\,{s_{2}}^2\,s_{3} \\&\qquad -
186624\,{s_{1}}^5\,s_{2}\,s_{5} + 93312\,{s_{1}}^5\,s_{3}\,s_{4} +
10935\,{s_{1}}^4\,{s_{2}}^4 + 11664\,{s_{1}}^4\,{s_{2}}^2\,s_{4}
 \\&\qquad + 104976\,{s_{1}}^4\,s_{2}\,{s_{3}}^2 - 1189728\,{s_{1}}^4
\,s_{2}\,s_{6} + 93312\,{s_{1}}^4\,{s_{4}}^2 - 78732\,{s_{1}}^3\,{s_{2
}}^3\,s_{3} \\&\qquad + 437400\,{s_{1}}^3\,{s_{2}}^2\,s_{5} - 209952\,
{s_{1}}^3\,s_{2}\,s_{3}\,s_{4} - 104976\,{s_{1}}^3\,{s_{3}}^3
 \\&\qquad + 2729376\,{s_{1}}^3\,s_{3}\,s_{6} - 279936\,{s_{1}}^3\,s_{
4}\,s_{5} + 13122\,{s_{1}}^2\,{s_{2}}^5 + 196830\,{s_{1}}^2\,{s_{2}}^2
\,{s_{3}}^2 \\&\qquad + 629856\,{s_{1}}^2\,{s_{2}}^2\,s_{6} - 944784\,
{s_{1}}^2\,s_{2}\,s_{3}\,s_{5} - 209952\,{s_{1}}^2\,s_{2}\,{s_{4}}^2
 \\&\qquad + 314928\,{s_{1}}^2\,{s_{3}}^2\,s_{4} - 7558272\,{s_{1}}^2
\,s_{4}\,s_{6} + 2834352\,{s_{1}}^2\,{s_{5}}^2 - 78732\,s_{1}\,{s_{2}}
^4\,s_{3} \\&\qquad + 314928\,s_{1}\,{s_{2}}^3\,s_{5} + 157464\,s_{1}
\,{s_{2}}^2\,s_{3}\,s_{4} - 236196\,s_{1}\,s_{2}\,{s_{3}}^3 \\&\qquad
 + 3779136\,s_{1}\,s_{2}\,s_{3}\,s_{6} - 1259712\,s_{1}\,s_{2}\,s_{4}
\,s_{5} - 472392\,s_{1}\,{s_{3}}^2\,s_{5} \\&\qquad + 629856\,s_{1}\,s
_{3}\,{s_{4}}^2 + 13122\,{s_{2}}^6 - 157464\,{s_{2}}^4\,s_{4} + 118098
\,{s_{2}}^3\,{s_{3}}^2 \\&\qquad - 472392\,{s_{2}}^2\,s_{3}\,s_{5} +
629856\,{s_{2}}^2\,{s_{4}}^2 - 472392\,s_{2}\,{s_{3}}^2\,s_{4} +
177147\,{s_{3}}^4 \\&\qquad - 2834352\,{s_{3}}^2\,s_{6} + 1889568\,s_{
3}\,s_{4}\,s_{5} - 839808\,{s_{4}}^3 \end{align*}

\newpage

{\samepage
\section*{Appendix 2. Standard coordinates for $E_7$.}
This Appendix gives the standard coordinate functions
$\e_i$ for $E_7$, and can also
serve as a correction to the formulas of Bramble \cite[pp.\ 358-360]{[Bra]}.
The $A_i$ which we calculate here are
integer multiples of the $\e_i$ which clear denominators; Bramble's paper
contains the same multiples.  (Note that our $\e_i$ correspond to his
$\a_{ijk\ell}$.)
It is very impressive to observe
that Bramble, calculating by hand, was correct in the calculation of $A_2,
\ A_6,\ A_8$, and
$A_{10}$, and had only two incorrect coefficients for $A_{14}$.  However,
the formulas for $A_{12}$ and $A_{18}$ from \cite{[Bra]} are mostly wrong.
\begin{align*}
A_2 &=
{\e_{2}} = 3\,{s_{1}}^2 - 4\,s_{2}\\
A_6 &=
48\,{\e_{6}} = 18\,{s_{1}}^6 - 72\,{s_{1}}^4\,s_{2} + 96\,{s
_{1}}^3\,s_{3} + 32\,{s_{1}}^2\,{s_{2}}^2 - 96\,{s_{1}}^2\,s_{4} -32
\,s_{1}\,s_{2}\,s_{3} \\&\qquad + 96\,s_{1}\,s_{5} - 64\,s_{2}\,s_{4}
 + 48\,{s_{3}}^2 + 384\,s_{6} \\
A_8 &=
 48\,{\e_{8}} = -27\,{s_{1}}^8 + 144\,{s_{1}}^6\,s_{2} - 192
\,{s_{1}}^5\,s_{3} - 160\,{s_{1}}^4\,{s_{2}}^2 + 192\,{s_{1}}^4\,s_{4}
 \\&\qquad + 320\,{s_{1}}^3\,s_{2}\,s_{3} - 192\,{s_{1}}^3\,s_{5} -
128\,{s_{1}}^2\,s_{2}\,s_{4} - 160\,{s_{1}}^2\,{s_{3}}^2 + 128\,s_{1}
\,s_{3}\,s_{4} \\&\qquad - 2304\,s_{1}\,s_{7} + 768\,s_{2}\,s_{6} +
384\,s_{3}\,s_{5} - 256\,{s_{4}}^2 \\
A_{10} &=
 16\,{\e_{10}} = 3\,{s_{1}}^{10} - 20\,{s_{1}}^8\,s_{2} + 32
\,{s_{1}}^7\,s_{3} + 32\,{s_{1}}^6\,{s_{2}}^2 - 32\,{s_{1}}^6\,s_{4}
 - 96\,{s_{1}}^5\,s_{2}\,s_{3} \\&\qquad + 32\,{s_{1}}^5\,s_{5} + 64\,
{s_{1}}^4\,s_{2}\,s_{4} + 80\,{s_{1}}^4\,{s_{3}}^2 - 128\,{s_{1}}^4\,s
_{6} - 128\,{s_{1}}^3\,s_{3}\,s_{4} \\&\qquad - 256\,{s_{1}}^3\,s_{7}
 + 256\,{s_{1}}^2\,s_{2}\,s_{6} + 128\,{s_{1}}^2\,s_{3}\,s_{5} + 512\,
s_{1}\,s_{2}\,s_{7} - 512\,s_{1}\,s_{3}\,s_{6} \\&\qquad - 1024\,s_{3}
\,s_{7} + 256\,{s_{5}}^2 \\
A_{12} &=
 6912\,{\e_{12}} = -297\,{s_{1}}^{12} + 2376\,{s_{1}}^{10}\,
s_{2} - 3456\,{s_{1}}^9\,s_{3} - 5616\,{s_{1}}^8\,{s_{2}}^2 + 3456\,{s
_{1}}^8\,s_{4} \\&\qquad + 14976\,{s_{1}}^7\,s_{2}\,s_{3} - 3456\,{s_{
1}}^7\,s_{5} + 3328\,{s_{1}}^6\,{s_{2}}^3 - 11520\,{s_{1}}^6\,s_{2}\,s
_{4} \\&\qquad - 11520\,{s_{1}}^6\,{s_{3}}^2 - 6912\,{s_{1}}^6\,s_{6}
 - 9984\,{s_{1}}^5\,{s_{2}}^2\,s_{3} + 11520\,{s_{1}}^5\,s_{2}\,s_{5}
 \\&\qquad + 19584\,{s_{1}}^5\,s_{3}\,s_{4} - 20736\,{s_{1}}^5\,s_{7}
 - 1536\,{s_{1}}^4\,{s_{2}}^2\,s_{4} + 13440\,{s_{1}}^4\,s_{2}\,{s_{3}
}^2 \\&\qquad + 34560\,{s_{1}}^4\,s_{2}\,s_{6} - 14976\,{s_{1}}^4\,s_{
3}\,s_{5} - 11520\,{s_{1}}^4\,{s_{4}}^2 + 3072\,{s_{1}}^3\,s_{2}\,s_{3
}\,s_{4} \\&\qquad + 55296\,{s_{1}}^3\,s_{2}\,s_{7} - 10240\,{s_{1}}^3
\,{s_{3}}^3 - 55296\,{s_{1}}^3\,s_{3}\,s_{6} + 18432\,{s_{1}}^3\,s_{4}
\,s_{5} \\&\qquad - 18432\,{s_{1}}^2\,{s_{2}}^2\,s_{6} - 9216\,{s_{1}}
^2\,s_{2}\,s_{3}\,s_{5} - 6144\,{s_{1}}^2\,s_{2}\,{s_{4}}^2 + 12288\,{
s_{1}}^2\,{s_{3}}^2\,s_{4} \\&\qquad - 55296\,{s_{1}}^2\,s_{3}\,s_{7}
 + 27648\,{s_{1}}^2\,{s_{5}}^2 - 110592\,s_{1}\,{s_{2}}^2\,s_{7} +
73728\,s_{1}\,s_{2}\,s_{3}\,s_{6} \\&\qquad - 18432\,s_{1}\,{s_{3}}^2
\,s_{5} + 6144\,s_{1}\,s_{3}\,{s_{4}}^2 + 221184\,s_{1}\,s_{4}\,s_{7}
 - 110592\,s_{1}\,s_{5}\,s_{6} \\&\qquad + 55296\,s_{2}\,s_{3}\,s_{7}
 + 36864\,s_{2}\,s_{4}\,s_{6} - 55296\,{s_{3}}^2\,s_{6} + 18432\,s_{3}
\,s_{4}\,s_{5} \\&\qquad - 8192\,{s_{4}}^3 - 110592\,s_{5}\,s_{7} -
110592\,{s_{6}}^2 \\
A_{14} &=
 768\,{\e_{14}} = 27\,{s_{1}}^{14} - 252\,{s_{1}}^{12}\,s_{2
} + 384\,{s_{1}}^{11}\,s_{3} + 752\,{s_{1}}^{10}\,{s_{2}}^2 - 384\,{s
_{1}}^{10}\,s_{4} \\&\qquad - 2176\,{s_{1}}^9\,s_{2}\,s_{3} + 384\,{s
_{1}}^9\,s_{5} - 704\,{s_{1}}^8\,{s_{2}}^3 + 1792\,{s_{1}}^8\,s_{2}\,s
_{4} + 1664\,{s_{1}}^8\,{s_{3}}^2
\displaybreak[0]
\\&\qquad - 768\,{s_{1}}^8\,s_{6} +
2816\,{s_{1}}^7\,{s_{2}}^2\,s_{3} - 1280\,{s_{1}}^7\,s_{2}\,s_{5} -
2944\,{s_{1}}^7\,s_{3}\,s_{4} \\&\qquad + 768\,{s_{1}}^7\,s_{7} - 1536
\,{s_{1}}^6\,{s_{2}}^2\,s_{4} - 3968\,{s_{1}}^6\,s_{2}\,{s_{3}}^2 +
2816\,{s_{1}}^6\,s_{2}\,s_{6} \\&\qquad + 2432\,{s_{1}}^6\,s_{3}\,s_{5
} + 1280\,{s_{1}}^6\,{s_{4}}^2 + 4608\,{s_{1}}^5\,s_{2}\,s_{3}\,s_{4}
 - 2048\,{s_{1}}^5\,s_{2}\,s_{7} \\&\qquad + 2048\,{s_{1}}^5\,{s_{3}}^
3 - 5120\,{s_{1}}^5\,s_{3}\,s_{6} - 2048\,{s_{1}}^5\,s_{4}\,s_{5} -
1024\,{s_{1}}^4\,{s_{2}}^2\,s_{6} \\&\qquad - 512\,{s_{1}}^4\,s_{2}\,s
_{3}\,s_{5} - 1024\,{s_{1}}^4\,s_{2}\,{s_{4}}^2 - 4096\,{s_{1}}^4\,{s
_{3}}^2\,s_{4} + 8192\,{s_{1}}^4\,s_{3}\,s_{7} \\&\qquad + 8192\,{s_{1
}}^4\,s_{4}\,s_{6} - 1536\,{s_{1}}^4\,{s_{5}}^2 + 4096\,{s_{1}}^3\,{s
_{2}}^2\,s_{7} + 2048\,{s_{1}}^3\,{s_{3}}^2\,s_{5} \\&\qquad + 2048\,{
s_{1}}^3\,s_{3}\,{s_{4}}^2 + 16384\,{s_{1}}^3\,s_{4}\,s_{7} - 12288\,{
s_{1}}^3\,s_{5}\,s_{6} - 30720\,{s_{1}}^2\,s_{2}\,s_{3}\,s_{7}
 \\&\qquad - 4096\,{s_{1}}^2\,s_{2}\,s_{4}\,s_{6} + 8192\,{s_{1}}^2\,s
_{2}\,{s_{5}}^2 - 2048\,{s_{1}}^2\,{s_{3}}^2\,s_{6} - 2048\,{s_{1}}^2
\,s_{3}\,s_{4}\,s_{5} \\&\qquad - 12288\,{s_{1}}^2\,s_{5}\,s_{7} +
12288\,{s_{1}}^2\,{s_{6}}^2 - 8192\,s_{1}\,s_{2}\,s_{4}\,s_{7} + 32768
\,s_{1}\,{s_{3}}^2\,s_{7} \\&\qquad + 8192\,s_{1}\,s_{3}\,s_{4}\,s_{6}
 - 8192\,s_{1}\,s_{3}\,{s_{5}}^2 + 49152\,s_{1}\,s_{6}\,s_{7} - 24576
\,s_{2}\,s_{5}\,s_{7} \\&\qquad + 16384\,s_{3}\,s_{4}\,s_{7} - 24576\,
s_{3}\,s_{5}\,s_{6} + 8192\,s_{4}\,{s_{5}}^2 + 49152\,{s_{7}}^2
\\
A_{18} &=
 9\cdot16^3\,{\e_{18}} = 63\,{s_{1}}^{18} - 756\,{s_{1}}^{16}\,s
_{2} + 1152\,{s_{1}}^{15}\,s_{3} + 3264\,{s_{1}}^{14}\,{s_{2}}^2 -
1152\,{s_{1}}^{14}\,s_{4} \\&\qquad - 9600\,{s_{1}}^{13}\,s_{2}\,s_{3}
 + 1152\,{s_{1}}^{13}\,s_{5} - 5888\,{s_{1}}^{12}\,{s_{2}}^3 + 8448\,{
s_{1}}^{12}\,s_{2}\,s_{4} \\&\qquad + 7488\,{s_{1}}^{12}\,{s_{3}}^2 +
24576\,{s_{1}}^{11}\,{s_{2}}^2\,s_{3} - 7680\,{s_{1}}^{11}\,s_{2}\,s_{
5} - 13824\,{s_{1}}^{11}\,s_{3}\,s_{4} \\&\qquad + 4608\,{s_{1}}^{11}
\,s_{7} + 3584\,{s_{1}}^{10}\,{s_{2}}^4 - 16896\,{s_{1}}^{10}\,{s_{2}}
^2\,s_{4} - 36864\,{s_{1}}^{10}\,s_{2}\,{s_{3}}^2 \\&\qquad - 3072\,{s
_{1}}^{10}\,s_{2}\,s_{6} + 12288\,{s_{1}}^{10}\,s_{3}\,s_{5} + 6912\,{
s_{1}}^{10}\,{s_{4}}^2 - 17920\,{s_{1}}^9\,{s_{2}}^3\,s_{3} \\&\qquad
 + 12800\,{s_{1}}^9\,{s_{2}}^2\,s_{5} + 53760\,{s_{1}}^9\,s_{2}\,s_{3}
\,s_{4} - 29184\,{s_{1}}^9\,s_{2}\,s_{7} + 20480\,{s_{1}}^9\,{s_{3}}^3
 \\&\qquad + 7680\,{s_{1}}^9\,s_{3}\,s_{6} - 12288\,{s_{1}}^9\,s_{4}\,
s_{5} + 5120\,{s_{1}}^8\,{s_{2}}^3\,s_{4} + 37120\,{s_{1}}^8\,{s_{2}}^
2\,{s_{3}}^2 \\&\qquad + 13312\,{s_{1}}^8\,{s_{2}}^2\,s_{6} - 39424\,{
s_{1}}^8\,s_{2}\,s_{3}\,s_{5} - 21504\,{s_{1}}^8\,s_{2}\,{s_{4}}^2 -
49152\,{s_{1}}^8\,{s_{3}}^2\,s_{4} \\&\qquad + 52224\,{s_{1}}^8\,s_{3}
\,s_{7} + 24576\,{s_{1}}^8\,s_{4}\,s_{6} - 16128\,{s_{1}}^8\,{s_{5}}^2
 - 20480\,{s_{1}}^7\,{s_{2}}^2\,s_{3}\,s_{4} \\&\qquad + 106496\,{s_{1
}}^7\,{s_{2}}^2\,s_{7} - 40960\,{s_{1}}^7\,s_{2}\,{s_{3}}^3 - 81920\,{
s_{1}}^7\,s_{2}\,s_{3}\,s_{6} \\&\qquad + 40960\,{s_{1}}^7\,s_{2}\,s_{
4}\,s_{5} + 45056\,{s_{1}}^7\,{s_{3}}^2\,s_{5} + 36864\,{s_{1}}^7\,s_{
3}\,{s_{4}}^2 - 24576\,{s_{1}}^7\,s_{4}\,s_{7} \\&\qquad - 8192\,{s_{1
}}^6\,{s_{2}}^3\,s_{6} - 4096\,{s_{1}}^6\,{s_{2}}^2\,s_{3}\,s_{5} -
8192\,{s_{1}}^6\,{s_{2}}^2\,{s_{4}}^2 \\&\qquad + 49152\,{s_{1}}^6\,s
_{2}\,{s_{3}}^2\,s_{4} - 299008\,{s_{1}}^6\,s_{2}\,s_{3}\,s_{7} -
90112\,{s_{1}}^6\,s_{2}\,s_{4}\,s_{6} \\&\qquad + 86016\,{s_{1}}^6\,s
_{2}\,{s_{5}}^2 + 20480\,{s_{1}}^6\,{s_{3}}^4 + 77824\,{s_{1}}^6\,{s_{
3}}^2\,s_{6} - 77824\,{s_{1}}^6\,s_{3}\,s_{4}\,s_{5} \\&\qquad - 8192
\,{s_{1}}^6\,{s_{4}}^3 + 73728\,{s_{1}}^6\,{s_{6}}^2 - 212992\,{s_{1}}
^5\,{s_{2}}^3\,s_{7} + 131072\,{s_{1}}^5\,{s_{2}}^2\,s_{3}\,s_{6}
 \\&\qquad - 40960\,{s_{1}}^5\,s_{2}\,{s_{3}}^2\,s_{5} + 24576\,{s_{1}
}^5\,s_{2}\,s_{3}\,{s_{4}}^2 + 139264\,{s_{1}}^5\,s_{2}\,s_{4}\,s_{7}
 \\&\qquad - 49152\,{s_{1}}^5\,s_{2}\,s_{5}\,s_{6} - 49152\,{s_{1}}^5
\,{s_{3}}^3\,s_{4} + 229376\,{s_{1}}^5\,{s_{3}}^2\,s_{7} \\&\qquad +
90112\,{s_{1}}^5\,s_{3}\,s_{4}\,s_{6} - 122880\,{s_{1}}^5\,s_{3}\,{s_{
5}}^2 + 32768\,{s_{1}}^5\,{s_{4}}^2\,s_{5}
\displaybreak[0]
\\&\qquad + 147456\,{s_{1}}
^5\,s_{6}\,s_{7} + 778240\,{s_{1}}^4\,{s_{2}}^2\,s_{3}\,s_{7} + 81920
\,{s_{1}}^4\,{s_{2}}^2\,s_{4}\,s_{6} \\&\qquad - 81920\,{s_{1}}^4\,{s
_{2}}^2\,{s_{5}}^2 - 286720\,{s_{1}}^4\,s_{2}\,{s_{3}}^2\,s_{6} +
40960\,{s_{1}}^4\,s_{2}\,s_{3}\,s_{4}\,s_{5} \\&\qquad - 32768\,{s_{1}
}^4\,s_{2}\,{s_{4}}^3 - 122880\,{s_{1}}^4\,s_{2}\,s_{5}\,s_{7} -
245760\,{s_{1}}^4\,s_{2}\,{s_{6}}^2 \\&\qquad + 81920\,{s_{1}}^4\,{s_{
3}}^3\,s_{5} - 409600\,{s_{1}}^4\,s_{3}\,s_{4}\,s_{7} + 122880\,{s_{1}
}^4\,s_{3}\,s_{5}\,s_{6} \\&\qquad - 131072\,{s_{1}}^4\,{s_{4}}^2\,s_{
6} + 122880\,{s_{1}}^4\,s_{4}\,{s_{5}}^2 + 147456\,{s_{1}}^4\,{s_{7}}^
2 \\&\qquad - 327680\,{s_{1}}^3\,{s_{2}}^2\,s_{4}\,s_{7} - 819200\,{s
_{1}}^3\,s_{2}\,{s_{3}}^2\,s_{7} + 163840\,{s_{1}}^3\,s_{2}\,s_{3}\,{s
_{5}}^2 \\&\qquad - 393216\,{s_{1}}^3\,s_{2}\,s_{6}\,s_{7} + 163840\,{
s_{1}}^3\,{s_{3}}^3\,s_{6} - 163840\,{s_{1}}^3\,{s_{3}}^2\,s_{4}\,s_{5
} \\&\qquad + 65536\,{s_{1}}^3\,s_{3}\,{s_{4}}^3 + 196608\,{s_{1}}^3\,
s_{3}\,s_{5}\,s_{7} + 393216\,{s_{1}}^3\,s_{3}\,{s_{6}}^2 \\&\qquad -
262144\,{s_{1}}^3\,{s_{4}}^2\,s_{7} + 393216\,{s_{1}}^3\,s_{4}\,s_{5}
\,s_{6} - 294912\,{s_{1}}^3\,{s_{5}}^3 \\&\qquad - 983040\,{s_{1}}^2\,
{s_{2}}^2\,s_{5}\,s_{7} + 589824\,{s_{1}}^2\,{s_{2}}^2\,{s_{6}}^2 +
1572864\,{s_{1}}^2\,s_{2}\,s_{3}\,s_{4}\,s_{7} \\&\qquad - 393216\,{s
_{1}}^2\,s_{2}\,s_{3}\,s_{5}\,s_{6} - 131072\,{s_{1}}^2\,s_{2}\,{s_{4}
}^2\,s_{6} + 131072\,{s_{1}}^2\,s_{2}\,s_{4}\,{s_{5}}^2 \\&\qquad -
393216\,{s_{1}}^2\,s_{2}\,{s_{7}}^2 + 327680\,{s_{1}}^2\,{s_{3}}^3\,s
_{7} - 131072\,{s_{1}}^2\,{s_{3}}^2\,s_{4}\,s_{6} \\&\qquad + 65536\,{
s_{1}}^2\,{s_{3}}^2\,{s_{5}}^2 - 65536\,{s_{1}}^2\,s_{3}\,{s_{4}}^2\,s
_{5} + 393216\,{s_{1}}^2\,s_{3}\,s_{6}\,s_{7} \\&\qquad + 393216\,{s_{
1}}^2\,s_{4}\,s_{5}\,s_{7} - 393216\,{s_{1}}^2\,s_{4}\,{s_{6}}^2 +
1179648\,s_{1}\,{s_{2}}^2\,s_{6}\,s_{7} \\&\qquad + 1572864\,s_{1}\,s
_{2}\,s_{3}\,s_{5}\,s_{7} - 1179648\,s_{1}\,s_{2}\,s_{3}\,{s_{6}}^2 -
262144\,s_{1}\,s_{2}\,{s_{4}}^2\,s_{7} \\&\qquad - 1441792\,s_{1}\,{s
_{3}}^2\,s_{4}\,s_{7} + 393216\,s_{1}\,{s_{3}}^2\,s_{5}\,s_{6} +
262144\,s_{1}\,s_{3}\,{s_{4}}^2\,s_{6} \\&\qquad - 131072\,s_{1}\,s_{3
}\,s_{4}\,{s_{5}}^2 + 393216\,s_{1}\,s_{3}\,{s_{7}}^2 - 1572864\,s_{1}
\,s_{4}\,s_{6}\,s_{7} \\&\qquad - 1179648\,s_{1}\,{s_{5}}^2\,s_{7} +
1179648\,s_{1}\,s_{5}\,{s_{6}}^2 + 589824\,{s_{2}}^2\,{s_{7}}^2
 \\&\qquad - 1179648\,s_{2}\,s_{3}\,s_{6}\,s_{7} - 393216\,s_{2}\,s_{4
}\,s_{5}\,s_{7} + 589824\,{s_{3}}^2\,{s_{6}}^2 \\&\qquad + 524288\,s_{
3}\,{s_{4}}^2\,s_{7} - 393216\,s_{3}\,s_{4}\,s_{5}\,s_{6} + 65536\,{s
_{4}}^2\,{s_{5}}^2 - 1572864\,s_{4}\,{s_{7}}^2 \\&\qquad + 2359296\,s
_{5}\,s_{6}\,s_{7}
\end{align*}
}

\end{document}